\ifx\rokhloaded\undefined\input rokh.def \fi
\input epsf

\def\mnote#1{}

\MakeToc{toc}

\topmatter
\author
A.~Degtyarev and V.~Kharlamov
\endauthor
\title
Topological properties of real algebraic varieties: 
du cote de chez Rokhlin
\endtitle
\rightheadtext{Topological properties of real algebraic varieties}

\dedicatory
 To 
  Vladimir Abramovich Rokhlin
\bigskip\bigskip
\begingroup
\def\\{\noalign{\smallskip}}
\font\eightcyr=wncyr8

\chardef\ya="1F
\chardef\yo="1B
\chardef\ssig="7E
\chardef\MZ="5E

\def\epigraph{\vbox{\def\\{\kern0pt}\eightpoint\eightcyr
 \halign{&##\hss\cr
  Inye vospominani\ya{} napolovinu st\yo{}rlis\ssig\cr
  v e\yo{} pam\ya{}ti, drugie ne sotrut\\s\ya{} nikogda.\footnotemark\cr}}}

\dimendef\bxsize=0
\bxsize=3in
\dimendef\unit=2
\unit=1pt
\catcode`\@=11
\def\over[#1,#2]#3{\kern#1\unit\raise#2\unit\vbox{%
 \hbox to\z@{\hss#3\hss}}\kern-#1\unit}

\let\lstyle\scriptscriptstyle
\def\lup#1{\eightpoint$\m@th \setbox\z@\hbox{$\lstyle\braceld$}\lstyle
  \bracelu\leaders\vrule height\ht\z@ depth\z@\hskip#1\unit\braceru$}
\def\ldown#1{$\m@th \setbox\z@\hbox{$\lstyle\braceld$}\lstyle
  \braceld\leaders\vrule height\ht\z@ depth\z@\hskip#1\unit\bracerd$}
\def\ligeup[#1,#2,#3]{\over[#1,#2]{\lup{#3}}}
\def\ligedown[#1,#2,#3]{\over[#1,#2]{\ldown{#3}}}

\def\wv[#1,#2,#3]{\over[#1,#2]{\vbox{%
 \setboxz@h{$\m@th\scriptstyle\lhook$}\dp\z@-1.5\p@%
 \cleaders\box\z@\vskip#3\unit}}}

\def\staff{\vbox{\hbox{\epsfxsize=\bxsize$
  \over[30,32]{$\pmb{pp}$}
  \over[60,0]{\sevenit con Ped.}
  \over[45,29].\over[45,40].
  \wv[40,43,12]\wv[40,15,12]
  \over[81,27].\over[81,38].
  \over[98.5,25.5].\over[98.5,38].
  \over[121,27].\over[121,39.5].
  \wv[115,39,8]\wv[115,11,16]
  \over[139.5,25.5].
  \wv[134.5,11,13]
  \over[157,25].\over[157,39.5].
  \wv[152.5,11,13]
  \over[175,38].
  \wv[170,11,13]
  \over[196,29].\over[196,40].
  \wv[192,43,12]\wv[192,15,12]
  \ligeup[71.5,35,13]\ligedown[71.5,30,13]
  \ligeup[148,35,13]
  \ligedown[187,36.5,14]
  \epsffile{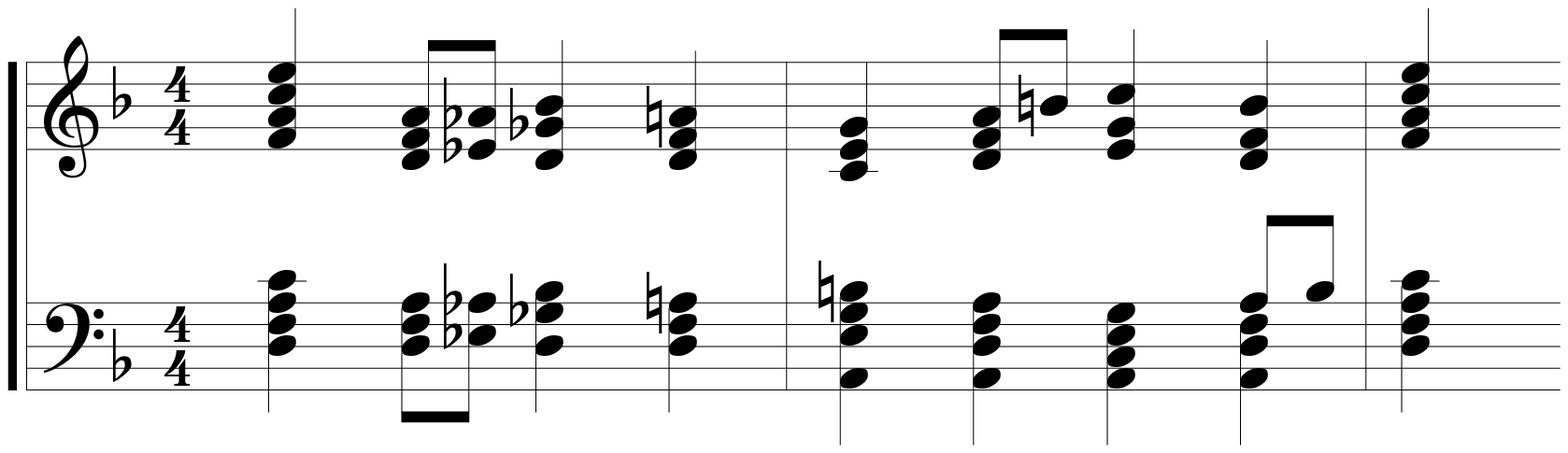}
  $}}}

\def\epititle{\bf Contes de la Vieille Grand' M\`ere}

\def\episign{\eightcyr S.~PROKOF\MZ{}EV. Soq.~31 (1918)}

\tenpoint
\line{\hss\vbox{\halign{#\cr
 \hss\epititle\hss\cr\\\\\\
 \hss\epigraph\cr\\\\\\\\
 \staff\cr
 \hss\episign\cr
}}}
\endgroup

\enddedicatory
\contents
 \nosubsection
 \nosubsubsection
\endcontents
\abstract The survey gives an overview of the achievements in topology of
real algebraic varieties in the direction initiated in the early 70th by
V.~I.~Arnold and V.~A.~Rokhlin. We make an attempt to systematize the
principal results in the subject. After an exposition of general tools
and results, special attention is paid to surfaces and curves on surfaces.
\endabstract

\address
Bilkent University, Ankara, Turkey
\endaddress
\email degt\,\@\,fen.bilkent.edu.edu
\endemail

\address
Institut de Recherche Math\'ematique Avanc\'ee\newline
\indent Universit\'e Louis Pasteur et CNRS\newline
\indent 7 rue Ren\'e-Descartes\newline
\indent 67084, Strasbourg, France
\endaddress
\email kharlam\,\@\,math.u-strasbg.fr
\endemail

\subjclass 
14P25, 14J99, and 57S25
\endsubjclass

\keywords
Real algebraic variety, real algebraic surface, real algebraic
curve, Hilbert's 16th problem
\endkeywords

\thanks
The title imitates that of M.~Proust's novel, usually translated as {\it
Swan's way} (cf\. also similar play on words in L.~Guillou, A.~Marin,
{\it A la recherche de la topologie perdue\/}).
\endgraf
\footnotemark[1] Some of her memories are half gone; others will stay
forever. {\it S.~Prokofiev,} Op.~31 (1918)
\endthanks

\endtopmatter

\begin
\document

\begingroup
\def\include#1{\input #1\relax}
\let\end\endinput


\ifx\rokhloaded\undefined\input rokh.def \fi

\def\Arf{\operator{Arf}}
\def\+{\mathbin{\scriptstyle\sqcup}}
\def\PD{\operator{D\null}}
\def\CC{\Cal C}

\begin
\document

\section{Introduction}\label{intro}

\counterstyle\thm{\thesection.\arabicnum\thmno}

\subsubsection**{The break-through}
To great extent, the interest shown by Vladimir Abramovich Rokhlin to
topology of real algebraic varieties was motivated by\mnote{toZv: paragraf
edited}
the results obtained in the late 60th by D.~Gudkov~\cite{Gudkov69} and
the subsequent paper by V.~Arnol$'$d \cite{Arnold},
which made a considerable contribution
to the solution of Hilbert's 16th problem.
Gudkov disproved one of Hilbert's conjectures on the arrangement of ovals
(i.e., two-sided components) of plane $M$-sextics (i.e., sextics with the
maximal number of ovals).
He corrected the conjecture, proved it for degree~$6$, and
suggested~\cite{Gudkov71}, as a new conjecture, an extension of his
result to $M$-curves of any even degree.  (Recall that an $M$-curve of
genus~$g$ is a curve with the maximal number of connected components of
the real part; due to Harnack's bound the maximal number  is $M=g+1$.)

To state Gudkov's conjecture, recall that an oval of a plane curve is
called \emph{even} (\emph{odd}) if it lies inside an even (respectively,
odd) number of other ovals. The number of even (odd) ovals is denoted
by~$p$ (respectively,~$n$). In this notation Gudkov's conjecture claims
that $p-n=k^2\bmod8$ for any $M$-curve of degree $2k$.

V.~Arnol$'$d in his remarkable paper~\cite{Arnold} related the study of
real plane curves to topology of $4$-manifolds and arithmetics of
integral quadratic forms and, besides other results, proved a weaker
version of Gudkov's conjecture ($p-n=k^2\bmod 4$). After that the events
were developing swiftly. First, Rokhlin suggested a proof of Gudkov's
conjecture based on his formula relating the signature of a $4$-manifold
with the $\Arf$-invariant of a characteristic surface (see
Appendix~\ref{RGMform}). Then, after Kharlamov's~\cite{Kh:10}
generalization of Arnol$'$d's results to surfaces, Rokhlin~\cite{Ro3}
found another proof, not using any specific tools of $4$-dimensional
topology, and extended the statement to varieties of any dimension. The
revolutionary break-through originated by Arnol$'$d was over and a period
of systematic study started. (More details and a brief account of the
further history of the subject can be found in~\cite{Gudkov:survey},
\cite{Wilson},\mnote{toZv:additional reference} and~\cite{ArnOlei}.)

\subsubsection**{Rokhlin's heritage}
Rokhlin published six papers on topology of real algebraic varieties (see
\cite{Ro}--\cite{Ro5}). This number is not very large, but each of these
papers originated a whole new direction in the subject.\mnote{toZv:
paragraf edited and extended}
(The only exception is probably the short note \cite{Ro4}, which extended
the range of applications of some previous results from complete
intersections to arbitrary real algebraic varieties.) In
this survey we discuss, in more or less details, all papers
except~\cite{Ro5}, where nonalgebraic coverings are used to study
algebraic curves; for other examples of this approach see~\cite{Fie:ro}.

In\mnote{toZv: par} Rokhlin's first paper \cite{Ro} there is a mistake in
the proof of Gudkov's conjecture. However, the approach of the paper,
namely, using characteristic surfaces in a $4$-manifold to evaluate its
signature $\bmod\,16$, became a powerful method in the study of real
algebraic curves. It was used by A.~Marin, who, together with L.~Guillou
(see~\cite{GM}), extended Rokhlin's signature formula to nonorientable
characteristic surfaces, and by means of it corrected the mistake. (The
Rokhlin-Guillou-Marin formula and related quadratic form are discussed in
Appendix~\ref{RGMform}.)

Another fundamental result, difficult to overestimate, is Rokhlin's
formula of complex orientations. As soon as
observed,\mnote{toZv: edited; a piece on Klein moved ahead} the notion of
complex orientation of a dividing real curve (see below), as well as
Rokhlin's formula and its proof, seem incredibly transparent.
The formula settles, for example, two of\mnote{toZv:sentence changed}
D.~Hilbert's conjectures on 11~ovals of plane sextics, which the latter
tried to prove in a very sophisticated way and included into his famous
problem list (as $16$-th problem).  Nowadays Rokhlin's formula is a major
tool in the study of the topology of real curves on surfaces. It is one
of the few phenomena related to curves that still do not have a
satisfactory generalization to varieties of higher dimension. (In fact,
even for curves on surfaces some basic questions have not been clarified
completely yet, cf.~\ref{C3.2}.)

It is worth mentioning that it\mnote{toZV: a piece moved and edited}
was F.~Klein
who
first studied dividing curves. He
discovered some of their remarkable properties
and pointed out
a way to improve and generalize Harnack's bound on the number of connected
components for such curves. However, he did not
notice the complex orientations.

\subsubsection**{Rokhlin's formula of complex orientation and Hilbert's 16th problem}
An irreducible\mnote{\Dg: edited} over~$\R$ nonsingular real curve $A$
(i.e., either an irreducible nonsingular complex curve with an
antiholomorphic involution $\conj\:A\to A$ or a nonsingular complex curve
consisting of two irreducible components transposed by an
anti-holomorphic involution) is called \emph{dividing} or \emph{of
type~$\I{}$}, if its real part $\R A=\Fix\conj$ divides $A$ into two
halves: two connected $2$-manifolds $A_+$ and $A_-$ having $\R A$ as
their common boundary. (Note that if $A$ is irreducible over~$\R$ but
reducible over~$\C$, it is dividing.)\mnote{\Dg: remark added} The complex
conjugation $\conj\:A\to A$ interchanges $A_\pm$ and the complex
orientation of $A_\pm$ induces two opposite orientations on $\R A$,
called its \emph{complex orientations}. If the curve is reducible or
singular, it is called dividing if so is the normalization of each real
component. Clearly, $A$ is dividing if and only if the quotient
$A/\!\conj$ is orientable, and if this is the case, the complex
orientations\mnote{toZv: plurial} are induced by orientations of
$A/\!\conj$. The principal example is provided by $M$-curves, which are
all dividing.

Let $A$ be a nonsingular dividing plane curve. Two ovals are said to form
an \emph{injective pair} if they bound an annulus in~$\Rp2$. An injective
pair is called \emph{positive} if the complex orientations of the ovals
are induced from an orientation of the annulus, and \emph{negative}
otherwise. Denote by~$\Pi^+$ and~$\Pi^-$ the numbers of, respectively,
positive and negative injective pairs, and let $\Pi=\Pi^++\Pi^-$.

\claim[Rokhlin's formula]\label{Rokh.form}
For a nonsingular dividing plane curve of even degree $2k$ with $l$ real
components one has $2(\Pi^+-\Pi^-)=l-k^2$.
\endclaim

Rokhlin's formula has numerous applications. For example, the Arnol$'$d
and Slepian congruences are its straightforward consequences. Given a
dividing curve of degree $2k$, the Arnol$'$d congruence states that
$p-n=k^2\bmod4$, and the Slepian congruence states that $p-n= k^2\bmod8$
provided that each odd oval contains immediately inside it an odd number
of even ovals (so that it bounds from outside a component of the
complement of the curve with even Euler
characteristic).\mnote{toZv:sentence edited}
Another consequence of Rokhlin's formula is the fact that for any
dividing curve of degree $2k$ one has $\Pi\ge\ls|l-k^2|$ and $l\ge k$.
The latter, together with the Klein congruence $l=k\bmod2$, which also
follows from Rokhlin's formula, is the complete set of restrictions on
the number of components of a nonsingular dividing plane curve of
degree~$2k$; for curves of degree $2k-1$ the corresponding version of
Rokhlin's formula (Rokhlin-Mishachev formula, see,  e.g., \cite{Ro2})
implies $l\ge k$ and a similar statement holds.\footnote{As we learned
from A.~Gabard, the problem on the number of components of a dividing
curve of a given degree was first posed by F.~Klein.}

Given a dividing curve~$A$, its \emph{real scheme} (i.e., the topological
type of $(\Rp2,\RA)$) and the complex orientation determine $A\cup\Rp2$
up to homeomorphism. (If the curve is not dividing, $A\cup\Rp2$ is
determined by the real scheme of~$A$.) One of the related problems,
suggested by Rokhlin, is the study of complex orientations appearing in a
given degree; it occupies an intermediate position between the stronger
question on the equivariant topology of the pair $(\Cp2,A)$ and the
weaker one on $(\Rp2,\RA)$ (the original setting of Hilbert's $16$th
problem). Two other related problems are worth mentioning here: the study
of the fundamental group $\pi_1(\Cp2\sminus(\Rp2\cup\RA))$ and of the
topology of the subspace $(\Rp2\cup A)\subset\Cp2/\!\conj=S^2$.

\subsubsection**{Rokhlin's school}
Rokhlin's influence on topology of real algebraic varieties extends far
beyond his own results. He directed into the subject a number of his
students (Cheponkus, Fiedler, Finashin, Kharlamov, Mischachev, Slepian,
Viro, Zvonilov) and generously shared his ideas and broad knowledge.
Unfortunately, Vladimir Abramovich could not realize all his plans: his
untimely death
interrupted his work on quite a number of projects.

As early as in
the middle
seventies\mnote{toZv:sentence edited} Rokhlin
had an intention to write two detailed surveys of the results achieved by
that time: one on real plane curves, and another one on higher dimensional
varieties. At the very beginning the project was joined by V.~Kharlamov
and later, by O.~Viro. Detailed plans were elaborated, Rokhlin wrote
summaries of the first chapters of the book on plane curves, and
Kharlamov, summaries of the first chapters of the other book. Rokhlin was
planning to deliver graduate courses in Leningrad State University;
unfortunately, he was forced to retire, and the course on plane curves
was given first by Viro (in a shorten version) and then by Viro and
Kharlamov. Kharlamov and Viro continued the work on the first project
after Rokhlin had passed away. The work was interrupted several times,
due to various political and geographical reasons. However, a draft of
the book has been written and a considerable part of it has been
polished. A preliminary version of the book has been used by a number of
graduate students as an introduction to the subject. Pitifully, it seems
that the book will never appear (partially due to the fact that the
subject is developing faster than the text is being written). The
introductory section of Viro's survey~\cite{Viro:survey} can give one an
idea about the first chapters of the book. (The survey, which
is mainly devoted to
other, construction aspects of the subject,\mnote{ toZv: a few
changements}
gives
a broader coverage of the principal notions and results, while the book
treats the subtle details that are only appropriate for a graduate text.)
The present survey
reflects the plan\mnote{toZv: ``in its turn'' removed} of the second
project; in its preliminary form it was designed as an appendix to the
book on curves (and, naturally, much inspired by it; many results
in~\ref{C3.3} and~\ref{s16.1} are straightforward generalizations of the
corresponding results for curves treated in the book).

\subsubsection**{What is and what is not contained in the survey}
We do not intend to give a complete overview of all known results on
topology of real algebraic varieties: the subject is so developed and
divers that it seems impossible to cover it within a reasonable volume.
Instead, we tried to select the results more closely related to the
phenomena discovered and studied by Rokhlin
and\mnote{toZv: ``himself'' and ``by'' removed}
his followers (certainly, even this choice is not exhaustive and is
mainly due to our personal taste) and made an attempt to illustrate how
the approaches and techniques developed work. Most proofs are either
omitted or just sketched. The citation list is also incomplete; in many
cases we refer to papers and textbooks where the concepts and proofs are
presented in a suitable way rather than to the original works.

Historically a great deal of attention was paid to the case of plane
curves; we address the reader to the excellent
surveys~\cite{Gudkov:survey;Wilson;Viro:survey} and concentrate on the
results in higher dimensions. However, at the end we return to surfaces
and curves on surfaces and discuss a few later developments, especially
those related to complex orientations. (Note that the topological
properties of abstract, not embedded, real curves are simple and have
been understood completely since F.~Klein, see,
e.g.,~\cite{Ro2;Natanzon:survey}.) A separate topic are real $3$-folds.
For a long time they stayed outside the main scope of the subject;
recently, due to J.~Koll\'ar, the situation started changing, and now the
$3$-folds should deserve a separate paper.

It turns out that most known results for plane curves are derived from
appropriate results for certain auxiliary surfaces; moreover, in many
cases the more general higher dimensional setting is, in fact, more
suitable as it gives proper understanding of the relation between the
topology of the real part of a variety and its complexification. However,
there are few exceptions, which are still not extended to higher
dimensions; the most remarkable of them is Rokhlin's formula of complex
orientations, mentioned above, and, among newer developments, results by
S.~Orevkov, based on a systematic study of pencils of lines. (Orevkov's
results, except his formula of complex orientations~\ref{orev:COF},
mainly deal with curves of low degrees; their general meaning and relation
to other approaches has not been completely revealed yet, even for plane
curves.)

Topology of real algebraic varieties is developing in two directions:
prohibitions and constructions. In this survey we completely ignore the
latter (see, e.g., \cite{Viro:survey;VIten;Risler} for an overview of
corresponding methods and results) and concentrate on the prohibition
type results, i.e., the restrictions on the topology of the real point
set of a real algebraic variety imposed by the topology of its complex
point set.\mnote{toZv: edited}
The latter is usually assumed known;
a typical example is considering varieties within a fixed complex
deformation family.

Some other important topics, currently developing but\mnote{toZv:edited}
ignored in this survey, are: special polynomials, fewnomials, complexity,
singularities and singular varieties, approximations, metric properties,
Ax principle, toric varieties, algebraic cycles, moduli spaces, minimal
models, relations to symplectic geometry.

Although our main subject are real algebraic varieties, many results
extend to much wider categories. In particular, instead of algebraic
varieties defined over~$\R$ we often consider closed complex manifolds
supplied with an anti-holomorphic involution. In many cases the complex
structure does not need to be integrable. (The most intriguing exception
is the generalized Comessatti-Petrovsky-Oleinik inequality~\ref{C1.3.A},
which has a topological proof for surfaces,
while the known proofs for higher dimensional varieties use the
integrability of the complex structure and even the existence of a
K\"ahler form.)\mnote{\Dg: edited} Moreover, many prohibition results are
topological in their nature and, thus, hold for arbitrary smooth
manifolds with involution; sometimes one needs to assume, in addition,
that locally the fixed point set of the involution behaves as the real
point set of a real algebraic variety, i.e., its normal\mnote{\Dg: edited}
and cotangent bundles are isomorphic. These phenomena are not extremely
surprising, considering the modern reconciliation of differential and
algebraic geometries. One can also speculate on the unity of complex and
real algebraic geometries, based, e.g., on the correspondence $X\mapsto
X\times\bar X$ ($\bar X$ standing for the complex conjugate variety) or
$X\mapsto X\+\bar X$;
in the latter case there is a bijection between the real structures
$f\:X\to\bar X$ and involutive automorphisms
of $X\+\bar X$ of the form $f\+\bar f$.\mnote{toZv:modified slightly}

\subsection**{Notation}
A \emph{real structure} on a complex analytic manifold~$X$ is an
anti-holom\-orphic involution. (In the case of algebraic varieties defined
over $\R$ it is the Galois involution.) Usually, we denote it by
$\conj\:X\to X$. The fixed point set of the real structure is called the
\emph{real part} of the variety and denoted by $\RX$. The quotient space
$X/\!\conj$ is usually denoted by~$X'$.

As it has been mentioned above, many results are topological in their
nature and, thus, hold for wider classes of varieties. Although it is
difficult to incorporate all the necessary hypotheses in a single
statement, in most cases the following definition is quite suitable: a
\emph{flexible real variety} is a closed smooth manifold~$X$ supplied
with a smooth involution $\conj\:X\to X$ and a (not necessarily
integrable) quasi-complex structure~$J$ in a neighborhood of the
\emph{real part} $\Fix\conj$, so that $\conj$ is anti-holomorphic in
respect to~$J$. (When the Bezout theorem is concerned, it may also be
useful to consider a symplectic structure compatible with both $\conj$
and~$J$.)\mnote{\Dg: new remark}

A closed complex (not necessarily irreducible or reduced) submanifold
$A\subset X$ is called \emph{real} if it is $\conj$-invariant. To avoid
the confusion between $\conj$-invariant complex submanifolds and real
submanifolds in the ordinary sense, we sometimes call the former
\emph{real smooth cycles} (or \emph{real smooth divisors}, if the complex
codimension is~$1$). A \emph{curve} in a complex manifold~$X$ is a
reduced effective cycle (divisor, if $X$~is a surface).

Many notions in algebraic geometry have different meaning over~$\C$ and
over~$\R$; in order to designate the real version we use the prefix~$\R$.
Thus, a submanifold~$A$ is \emph{$\R$-irreducible} if either it is
irreducible over~$\C$ or it has two components permuted by~$\conj$; an
\emph{$\R$-component} of~$A$ is either a $\conj$-invariant component or a
pair of components permuted by~$\conj$. (Realizing the awkwardness of
this terminology, we
adopt
it as we could not find anything better.
Designating these notions as `real' might cause confusion with their
topological counterparts, while `over~$\R$' is not always applicable,
grammatically or, even worse, semantically, as in general our manifolds
are not algebraic varieties defined over~$\R$.)

Given a linear system~$\ls|D|$ (see~\ref{C1.0}) on a complex manifold~$X$
with real structure, we denote by $\R\ls|D|$ its real part, which
consists of the real divisors linearly equivalent to~$D$, and by
$\Delta_{\ls|D|}$ and $\R\Delta_{\ls|D|}$, the discriminant of~$\ls|D|$
and its real part, respectively. In the particular case $X=\Cp{n}$ we use
the more common notation $\CC_q$, $\Delta_q$ for the space of
hypersurfaces of degree~$q$ and its discriminant, respectively, and
$\R\CC_q$, $\R\Delta_q$ for their real parts.

We use the standard bracket notation to encode finite partially ordered
sets. Typically this notation is applied to the set of ovals of a real
curve on a real surface. The particular partial order used is to be
specified explicitly for each class of curves; e.g., for two ovals $C_1$,
$C_2$ of a plane curve one has $C_1\succ C_2$ if $C_2$ belongs to the
disk bounded by~$C_1$.

Unless stated explicitly, all the homology and cohomology groups have
coefficients~$\ZZ$. We use $b_i(\,\cdot\,)$ and $\Gb_i(\,\cdot\,)$ for
the Betti numbers over~$\Q$ and over~$\ZZ$, respectively, and
$b_*(\,\cdot\,)$ and $\Gb_*(\,\cdot\,)$, for the corresponding total
Betti numbers. Given a topological space~$X$ with involution~$c$, we
denote by $(\rH_*(X),\rd_*)$ (or just $(\rH_*,\rd_*)$\,) Kalinin's
spectral sequence of~$X$, by~$\bv_*$, the Viro homomorphisms, and by
$\CF^*(X)=\CF^*$,
Kalinin's filtration on $H_*(\Fix c)$ (see~\ref{C1.7}).

If $X$ is a manifold, $w_i=w_i(X)$ stand for its Stiefel-Whitney classes,
$u_i=u_i(X)$, for its Wu classes, and, if $X$ is complex, $c_i=c_i(X)$,
for its Chern classes. If $X$ is closed, we denote by
$\PD=\PD_X\:H^*(X)\to H_*(X)$ the Poincar\'e duality isomorphism and often
use the same notation~$w_i$, $u_i$, $c_i$ for the Poincar\'e dual homology
characteristic classes.

\subsection**{Acknowledgements}
We are grateful to our numerous colleagues with whom we worked on various
topics in the subject and who shared with us their knowledge. The list of
persons who helped us to clarify certain details or communicated to us
their unpublished results during our work on this survey includes, but is
not limited to, B.~Chevallier, S.~Finashin, I.~Itenberg, A.~Marin,
S.~Orevkov, O.~Viro, J.-Y.~Welschinger. Our particular gratitude is to
A.~M.~Vershik, who encouraged us to write this survey.


\counterstyle\thm{\thesubsection.\arabicnum\thmno}

\end



\ifx\rokhloaded\undefined\input rokh.def \fi

{\catcode`\@=11 \gdef\proclaimfont@{\sl}
\global\let\proclaimfont\proclaimfont@}

\def\wt#1{\smash{\widetilde{#1}}}
\def\Pic{\operator{Pic}}

\def\BO{\text{\sl BO}}
\def\BSO{\text{\sl BSO}}


\begin

\section{General tools and results}\label{C1}

\subsection{Divisors and linear equivalence}\label{C1.0}
Recall that a \emph{divisor} on a nonsingular compact complex
manifold~$X$ is a formal finite integral linear combination $A=\sum
m_iD_i$, where $m_i\in\Z$ and $D_i$ are irreducible (possibly singular)
compact codimension 1 subvarieties of~$X$. If all $m_i\ge0$, $A$ is called
\emph{effective}. Divisors form a group in respect to the formal addition
of linear combinations. Two divisors~$D_1$, $D_2$ are called
\emph{linearly equivalent}, $D_1\sim D_2$, if their difference is a
\emph{principal divisor}, i.e., there is a meromorphic function on~$X$
whose zeros and poles are the components of $D-D'$, considered with their
multiplicities. All effective divisors linearly equivalent to a given
divisor~$D$ form a projective space; it is called a \emph{linear system}
and denoted by~$\ls|D|$.  Clearly, linearly equivalent divisors realize
the same class in $H^2(X)$; if $X$ is simply connected (or, more
generally, $H^1(X)=0$), this condition is also sufficient (see, e.g.,
\cite{Hirzebruch} or~\cite{Hartshorne}).

A divisor is called \emph{very ample} if it can be realized as a
hyperplane section under an appropriate embedding of~$X$ into a
projective space. A divisor~$D$ is called \emph{ample} if some positive
multiple~$mD$ is very ample. (Some useful criteria of ampleness can be
found in~\cite{Hartshorne}.)

If $A$ is a real divisor on a real variety $(X,\conj)$, the class of $A$
in $H^2(X;\Z)$ is $\conj^*$-skew-invariant. If $H^1(X)=0$, every divisor
whose class is $\conj^*$-skew-invariant is linearly equivalent to a real
one. Real divisors equivalent over $\C$ are equivalent over $\R$,
provided that $\RX\ne\varnothing$.

\claim[Lefschetz theorem on hyperplane section \ROM{(see, e.g.,~\cite{Milnor})}]\label{C1.0.A}
If $A$ is an ample divisor on a compact complex $n$-dimensional
manifold~$X$, then $X\sminus A$ has homotopy type of
$n$-dimensional {\sl CW\/}-complex. In particular, for any abelian
group~$G$ and $r<n$ one has $H_r(X,A;G)=H^r(X,A;G)=0$.
\qed
\endclaim

\subsection{Double coverings}\label{C1.1}
Double coverings play a special r\^ole in real algebraic geometry: the
position of a subvariety $A\subset X$ is reflected, to great extent, in
the topological properties of the double covering of~$X$ ramified in~$A$
(according to V.~Arnol$'$d's principle, the notion of double covering is
the complexification of that of manifold with boundary; Arnol$'$d claims
that it is due to this observation that he found a proof of relaxed
Gudkov's congruence and the other remarkable results of~\cite{Arnold:break}).
Thus, once an absolute result is obtained, it may be applied to branched
double coverings and thus produce some relative prohibitions on the
position of a subvariety.

A divisor $A\subset X$ is called \emph{even} if the fundamental
class~$[A]$ vanishes in $H_{2n-2}(X)$, where $n=\dim X$. In this case
there exists a double covering~$Y$ of~$X$ branched over~$A$, which is
also a complex variety, nonsingular if so is~$A$. Clearly, an isomorphism
class of the covering is determined by a class $\omega\in H^1(X\sminus
A)$ whose image under $H^1(X\sminus A)\to H_{2n-1}(X,A)\to H_{2n-2}(A)$
is~$[A]$. Alternatively, the homology vanishing condition is equivalent
to $A\sim2E$ for some divisor~$E$ on~$X$, and the isomorphism classes of
double coverings are in a canonical one-to-one correspondence with the
classes $E\in\Pic X$ such that $2E\sim A$, see~\cite{Hirz53}. To indicate
a particular choice of the covering we use the notation $Y(E)$.

The following is an immediate consequence of~\ref{C1.0.A}:

\claim[Lefschetz theorem for double covering]\label{C1.1.A}
Let $A$ be an even ample divisor on a compact complex manifold~$X$ and $Y$
a double covering of~$X$ branched over~$A$. Then for any abelian
group~$G$ one has
$$
H_r(Y,A;G)=H^r(Y,A;G)=0\quad\text{for $r<n=\dim X$}.\rlap\qed
$$
\endclaim

\claim[Corollary of \ref{C1.0.A} and \ref{C1.1.A}]\label{C1.1.B}
Let $A\subset X$ and\/ $Y$ be as in~\ref{C1.1.A} and $p\:Y\to X$ the
covering projection. Then for any abelian group~$G$ the induced map
$p_*\: H_r(Y;G)\to H_r(X;G)$ is an isomorphism for $r<n$ and an
epimorphism for $r=n$\rom; the induced map $p^*\:  H^r(X;G)\to H^r(Y;G)$
is an isomorphism for $r<n$ and a monomorphism for $r=n$. Furthermore, if
$G$ is a field, $\char G\ne2$, then $p_*$ induces an isomorphism between
$H_*(X;G)$ and the $(+1)$-eigenspace $H_*^+(Y;G)$ of the deck translation
\rom(respectively, $p^*$ induces an isomorphism between $H^*(X;G)$ and
$H^*_+(Y;G)$\,\rom).
\qed
\endclaim

Note that \ref{C1.1.A} and~\ref{C1.1.B} still hold if $A$ and, hence, $Y$
are singular. If $A$ is smooth but not even, one may try to find another
smooth divisor~$H$ on~$X$, intersecting~$A$ transversally, so that
$[H]+[A]$ vanishes in $H_{2n-2}(X)$. Then there is a double covering~$Y_H$
of~$X$ branched along~$H\cup A$, and one can study its blow-up~$\wt Y_H$
at~$S=H\cap A$. (A typical example is when $X\subset\Cp{N}$ and $A$~is cut
on~$X$ by a hypersurface of odd degree: taking for~$H$ a generic
hyperplane section one obtains information about an affine part of
$(X,A)$.) Equivalently, one can first blow up~$S$ at~$X$ to obtain a
manifold~$\wt X$; then $\wt Y$ is the double covering of~$\wt X$ branched
along the proper transform of $A+H$. If~$H$ is ample, so is $aH+A$ for
$a\gg0$, and since $H_r(\wt X)=H_r(X)\oplus H_{r-2}(S)$,
Corollary~\ref{C1.1.B} takes the following form:

\claim[Proposition]\label{C1.1.C}
Let $A\subset X$ be a smooth divisor, $H\subset X$ an ample divisor
transversal to~$A$, and $\wt Y_H$ as above. Then for any field~$F$,
$\char F\ne2$, and any $r\in\Z$ there are canonical isomorphisms
$H_r^+(\wt Y_H;F)=H_r(X;F)\oplus H_{r-2}(H\cap A;F)$ and $H^r_+(\wt
Y_H;F)=H^r(X;F)\oplus H^{r-2}(H\cap A;F)$.  Furthermore, $H_r^-(\wt
Y_H;F)=H^r_-(\wt Y_H;F)=0$ for $r\ne n=\dim X$. \rom(Here $H_*^\pm$ and
$H^*_\pm$ are the $(\pm1)$-eigenspaces of the deck translation
involution.\rom)
\qed
\endclaim

\subsubsection{}\label{real.cov}
Assume now that $X$ is supplied with a real structure and $A$ is an even
real divisor. Assume, further, that $E$ with $2E\sim A$ is also chosen
real and let $Y=Y(E)$. (Such a choice of~$E$ exists, e.g., if $H^1(X)=0$;
in general the existence depends on whether the class of~$A$ in
$H^2(X;\Z)$ is the double of a skew-invariant class.)
If $\RX\ne\varnothing$, the
involution $\conj$ on~$X$ lifts to two involutions $T_\pm\:Y\to Y$, which
are antiholomorphic, commute with each other and with the deck translation
$\tau\:Y\to Y$, and satisfy $\tau=T_+\circ\nomathbreak T_-$
(see~\ref{apdx.lifts}). The corresponding real parts of~$Y$ are denoted
by $\RY_\pm$, and their projections to~$\RX$, by~$\RX_\pm$. Obviously,
$\RX_+\cup\RX_-=\RX$ and $\RX_+\cap\RX_-=\RA$. The characteristic classes
$\omega_\pm\in H^1(\RX_\pm)$ are Poincar\'e dual to the corresponding
restrictions $\rel[\RE]\in H_{n-1}(\RX_\pm,\RA)$.

\subsection{Smith inequalities}\label{C1.2}
This group of results, obtained from the Smith inequality and its
relative version (see \iref{1.1.2}2 and~\eqref{1.1.12}, respectively),
generalizes the Harnack inequality for curves.

\claim[Theorem]\label{C1.2.A}
For\mnote{toZv:new title} a manifold~$X$ with real structure one has
$\beta_*(\RX)\le\beta_*(X)$ and $\beta(\RX)=\beta(X)\bmod2$. If,
furthermore, $A\subset X$ is a real submanifold, then
$\beta_*(\RX,\RA)\le\beta_*(X,A)$ and $\beta_*(\RX,\RA)=\beta(X,A)\bmod2$.
\qed
\endclaim

In~\ref{C1.2.A} one can easily recognize the familiar relations between
the number of real roots of a real polynomial and its degree, if $\dim
X=0$, and the classical Harnack inequality, which states that the number
of real components of $\RX$ does not exceed $\operatorname{genus}(X)+1$,
if $\dim X=1$. The bound given by~\ref{C1.2.A} is sharp for hypersurfaces
of any degree in projective spaces of any dimension (see \cite{ViroM} and
the forthcoming~\cite{ItenViro});\mnote{toZv:additional information and
reference} for a hypersurface $X\subset\Cp{q}$ of degree~$m$ the bound
takes the form $\beta_*(X) \le\frac1m((m-1)^{q+1}+(-1)^q)+q-(-1)^q$,
see~\cite{ThomSmith}.

\claim[Complementary inequality for even pairs]\label{C1.2.C}
Let $A$ be an even ample real submanifold on a manifold~$X$ with real
structure, $\dim X=n$. Let, further, $E$ be a real divisor on~$X$ with
$2E\sim A$ and $\omega_\pm\in H^1(\RX_\pm)$ the class Poincar\'e dual to
$\rel[\RE]\in H_{n-1}(\RX_\pm,\RA)$, where $\RX_\pm$ are the two halves
of $\RX\smallsetminus\RA$ \rom(see.~\ref{real.cov}\rom). Then
$$
2\beta_*(\RX_\pm)-
  2\bigl[\dim\Ker\partial-\dim\Ker(\omega_\pm\oplus\partial)\bigr]\le
\beta_*(X)+(-1)^n[\chi(X)-\chi(A)],
$$
where $\omega_\pm\oplus\partial\:H_*(\RX_\pm,\RA)\to
H_{*-1}(\RX_\pm,\RA)\oplus H_{*-1}(\RA)$ is the boundary homomorphism in
the Smith exact sequence~\ref{Smith.seq}.
\endclaim

\proof
The statement is the Smith inequality applied to one of the
involutions~$T_\pm$ in the double covering $Y\to X$ branched over~$A$
(see~\ref{real.cov}). The left hand side here is $\beta_*(\RY_\pm)$, as
it follows immediately from the Smith exact sequence of the deck
translation involution, exact sequence of~$(\RX,\RA)$, and Poincar\'e
duality. The right hand side equals $\beta_*(Y)$. Indeed,
from~\ref{C1.1.B} and the symmetry of the Betti numbers it follows that
$\beta_r(Y)=\beta_r(X)$ for all $r\ne n$, and to find the remaining Betti
number $\beta_n(Y)$ it suffices to compare the Euler characteristics
(using, e.g., the Riemann-Hurwitz formula \iref{1.1.2}6).
\endproof

\remark{Remark}
The left hand side of the inequality in~\ref{C1.2.C} above is equal to
$\beta_*(\RA)+2\dim\Ker(\omega_\pm\oplus\partial)$. In the case of
projective hypersurfaces the inequality can be simplified, see
\ref{l-prop}; in this case $\omega_\pm$ is $k$-times the generator of
$H^1(\RP^n)$, where $2k=\deg A$. (By a strange mistake, in the early
papers on the subject another class was indicated.)
\endremark

Unlike the case of curves, in general it is a difficult problem to
estimate~$\beta_i(\RX)$ (and even~$\beta_0(\RX)$\,) separately; the sharp
bound is not known even for surfaces of degree~5 in~$\Cp3$,
see~\ref{C2.3}. Note also that~\ref{C1.2.A} and~\ref{C1.2.C}, when they
both apply, give, in general, different restrictions to the topology of
$(\RX,\RA)$.  Certainly, the ampleness condition in~\ref{C1.2.C} is only
needed to evaluate the Betti numbers of~$Y$; without it the right hand
side of the inequality should be replaced with~$\beta_*(Y)$.

According to~\ref{C1.2.A}, the difference $\beta_*(X)-\beta_*(\RX)$ is a
nonnegative even integer. If this difference is~$2d$, one calls $X$ is an
\emph{$(M-d)$-manifold} (or $\conj$ an \emph{$(M-\nomathbreak
d)$-involution}). Similar to the case of curves and the classical Harnack
inequality it is $M$- (or close to~$M$-) manifolds that satisfy certain
additional congruence type prohibitions (see~\ref{C1.5}).

\subsection{Comessatti-Petrovsky-Oleinik inequalities}\label{C1.3}
Recall that, if a compact complex manifold $X$ admits a K\"ahler metric
(and this is the case for projective manifolds, since they inherit a
K\"ahler metric from the ambient projective space), there is a canonical
\emph{Hodge decomposition} (\emph{Hodge structure})
$$
H^r(X,\C)=\bigoplus_{p+q=r}H^{p,q}(X), \quad 0\le p,q\le\dim_{\C}X,
$$
where $H^{p,q}(X)=H^q(X;\Omega^p(X))$
can be interpreted as the subspace of classes realized by
$(p,q)$-forms.\mnote{\Dg: edited} The \emph{Hodge numbers}
$h^{p,q}(X)=\dim H^{p,q}(X)$ are deformation invariants of~$X$. A
particular K\"ahler metric yields a further decomposition
$H^{p,q}(X)=\bigoplus\Omega^k\wedge P^{p-k,q-k}(X)$, where $\Omega\in
H^{1,1}(X)$ is the fundamental class of the metric and $P^{a,b}\subset
H^{a,b}$ are the subspaces of so called \emph{primitive classes} (see,
e.g.,~\cite{Chern}). If $X$ is real, it admits an invariant K\"ahler
metric (e.g., obtained by the averaging) and for such a metric one has
$\conj^*\Omega=-\Omega$.

The following result incorporates the Petrovsky~\cite{Petrovsky},
Petrovsky-Oleinik~\cite{PetrovskyOleinik}, and
Comessatti~\cite{Comessatti} inequalities (for double planes, projective
hypersurfaces, and surfaces, respectively). For other generalizations of
Petrovsky-Oleinik inequalities see \cite{ArnoldPO;KhoPO;KhoPO2}.

\claim
[Theorem]\label{C1.3.A}
If\mnote{toZv:new title} $X$ is a compact complex K\"ahler manifold with
real structure, $\dim X=n=2k$ even, then
$$
\bigl|\chi(\RX)-1\bigr|\le h^{k,k}(X)-1.
$$
If, besides, $A\subset X$ is an even real divisor, $A\sim2E$ for a real
divisor~$E$, then
$$
\gather
\bigl|\chi(\RX_-)-\chi(\RX_+)\bigr|\le h_-^{k,k}(Y),\rlap{\quad and}\\
 \noalign{\vglue1\jot}
\bigl|2\chi(\RX_\pm)-1\bigr|\le h^{k,k}(Y)-1,
\endgather
$$
where $Y=Y(E)$, $\tau$ is the deck translation, and
$h_-^{pq}(Y)=h^{p,q}(Y)-h^{p,q}(X)$ is the dimension of the
$\tau$-skew-invar\-iant part $H_-^{pq}(Y)\subset H^{p,q}(Y)$.
\endclaim

\proof
We prove the first assertion; the two others are similar. According to
the Lefschetz fixed point theorem,
$$
\chi(\RX)=\sum_r(-1)^r\Trace(\conj^*,H^r(X)).
$$
Since $\conj_*H^{p,q}(X)=H^{q,p}(X)$ and $\conj^*\Omega=-\Omega$ (for an
invariant K\"ahler metric), in the decomposition of $H^*(X;\C)$ into
primitive classes only the terms present in~$H^{k,k}(X)=\sum
\Omega^j\wedge P^{k-j,k-j}$ may contribute to~$\chi(\RX)$:
$$
\chi(\RX)=\sum_{j=0}^{k}\Trace (\conj^*,P^{k-j,k-j}(X)).\eqtag\label{eqC.1}
$$
It remains to observe that $\conj^*=\id$ on $P^{0,0}(X)$ and that
$\ls|\Trace{}|\le\dim$.
\endproof

The
above inequalities are sharp for curves in $\Cp2$,
see~\cite{Petrovsky},
and for surfaces in $\Cp3$,
see~\cite{ViroM};
in these cases they are, respectively, the Petrovsky and
Petrovsky-Oleinik inequalities.
To our knowledge, in higher dimensions the question is still
open.\mnote{toZv:moved and modified}

\subsection{Hodge numbers}
If\mnote{toZv:a separate section} $A$ is ample, the last two inequalities
in~\ref{C1.3.A} can be made effective. Following F.~Hirzebruch, consider
the \emph{generalized Todd genus} (or just \emph{$T_y$-genus})
$T_y(X)=\sum(-1)^qh^{p,q}(X)y^p$. By the Hirzebruch-Riemann-Roch
theorem~\cite{Hirzebruch} $T_y(X)$ is equal to the value on $[X]$ of a
certain degree~$n$ polynomial $T_n(y;c_1,\dots,c_n)$ in the Chern classes
$c_i=c_i(X)$. Furthermore, given $u_1,\dots,u_r\in H^2(X)$, one can
define the \emph{virtual Todd genus} $T^r_y(u_1,\dots,u_r)$, which is a
polynomial in~$c_i$ and~$u_j$ with the following property: if the classes
Poincar\'e dual to~$u_i$ are realized by codimension~$1$
submanifolds~$U_i$ intersecting transversally, then
$T_y^r(u_1,\dots,u_r)=T_y(U_1\cap\ldots\cap U_r)$.

The following statement is a consequence of~\ref{C1.1.B} (with $G=\C$):

\claim[Proposition \ROM{(see~\cite{Kh:petr1})}]\label{C1.3.B}
Let $A\subset X$ and $Y=Y(E)$ be as in~\ref{C1.3.A}. If $A$ is ample, then
$h^{p,q}(Y)=h^{p,q}(X)$ for all $p$,~$q$ with $p+q\ne n$ and
$$
\sum_{i=0}^n(-1)^i(h^{n-i,i}(Y)-h^{n-i,i}(X))y^i=
\sum_{r\ge0}a_r(y)T_y^r(e,\dots,e)_X-T_y(X),
$$
where $a_r(y)$ is the coefficient of~$x^r$ in the formal power series
expansion
$$
\frac{(1+yx)^2-(1-x)^2}{(1+yx)^2+y(1-x)^2}\cdot\frac1x= \sum_{r=0}^\infty
a_r(y)x^r  \eqtag\label{eqC.2}
$$
and $e=c_1(E)$. \rom(The result does not depend on the choice of~$E$ with
$2E\sim A$, as the calculation is done in the rational homology.\rom)
\qed
\endclaim

\proof
From~\ref{C1.1.B}, the naturallity of the Hodge decomposition, and Serre
duality it follows that $h^{p,q}(Y)=h^{p,q}(X)$ for $p+q\ne n$, and to
find the remaining Hodge numbers in the middle dimension it suffices to
know the Todd genus of~$Y$, which, due to~\cite{Kh:petr1}, is given by
$T_y(Y)=\sum_{r\ge0}a_r(y)T_y^r(e,\dots,e)_X$.
\endproof

Explicit calculation for small dimensions gives the following values for
$T_y(X)=\sum_{i=0}^nT_{n,i}y^i$ (where $n=\dim X$):

\claim[Corollary]\label{C1.3.C}
Let $c_i=c_i(X)$. If $\dim X=2$, then
$$
\gather
T_{2,0}=T_{2,2}=\tfrac1{12}(c_1^2+c_2)[X],\qquad
T_{2,1}=\tfrac1{6}(c_1^2-5c_2)[X],\\
h^{1,1}(Y)=h^{1,1}(X)-T_{2,1}-(ec_1-3e^2)[X].
\endgather
$$
If $\dim X=4$, then
$$
\gather
T_{4,0}=T_{4,4}=
\tfrac1{720}(-c_1^4+4c_1^2c_2+3c_2^2+c_1c_3-c_4)[X],\\
T_{4,1}=T_{4,3}=
\tfrac1{180}(-c_1^4+4c_1^2c_2+3c_2^2-14c_1c_3-31c_4)[X],\\
T_{4,2}=
\tfrac1{120}(-c_1^4+4c_1^2c_2+3c_2^2-19c_1c_3+79c_4)[X],\\
\hbox to\hsize{\qquad $h^{2,2}(Y)=h^{2,2}(X)+T_{4,2}+{}$\hfil}\\
\hbox to\hsize{\hfil
 $\tfrac1{12}(115e^4-46e^3c_1-5e^2c_1^2+31e^2c_2+ec_1c_2-12ec_3)[X].$\rlap\qed\qquad}
\endgather
$$
\endclaim

If $X$ is a regular complete intersection in~$\Cp{N}$, its Hodge numbers
and, hence, the bounds in~\ref{C1.3.A} can be found recursively. They
depend only on the polydegree of~$X$.\mnote{toZv:modified} Define
polynomials $\chi_s^q(m_1,\dots,m_s;y)$ as follows:
$$
{\aligned
&\chi_0^0(y)=1,\quad
\chi_s^0(m_1,\dots,m_s;y)=0\quad\text{for $s>0$,\quad and}\\
&\chi_s^q(\dots,m_s;y)=m_s\chi_{s-1}^{q-1}(\dots;y)\\
&\qquad\qquad+\sum_{\mu=1}^{m_s}\bigl[(y-1)\chi_s^{q-1}(\dots,\mu-1;y)-
y\chi_{s+1}^{q-1}(\dots,\mu-1,\mu;y)\bigr].
\endaligned}\eqtag\label{eqC.3}
$$
(Note that $\chi_0^q(y)=\sum_{r=0}^q(-1)^ry^r\!$, \
$\chi_q^q(m_1,\dots,m_q;y)=m_1\ldots m_q$, and $\chi_s^q=0$ for~$s>q$.)

\claim[Proposition]\label{C1.3.D}
If a manifold~$X$ of even dimension~$n=2l$ is a regular complete
intersection in~$\Cp{n+s}$ of polydegree $(m_1,\dots,m_s)$, then
$$
T_y(X)=\chi_s^{n+s}(m_1,\dots,m_s;y).
$$
In particular, $(-1)^lh^{l,l}(X)$ is equal to the coefficient of~$y^l$ in
$\chi_s^{n+s}(m_1,\dots,m_s;y)$.

If, further, $A\subset X$ is a submanifold cut on~$X$ by a hypersurface of
degree~$2k$ and $Y$ is the double covering of~$X$ branched over~$A$, then
$$
T_y(Y)= \sum_{r=0}^{n}a_r(y)\chi_{r+s}^{n+s}(m_1,\dots,m_s,k,\dots,k;y).
$$
\rom(See~\eqref{eqC.2} for the definition of~$a_r(y)$.\rom) In particular,
$(-1)^lh^{l,l}(Y)$ is equal to the coefficient of~$y^l$ in the above
polynomial.
\endclaim

\proof[Proof\/ \ROM(see~\cite{Kh:petr1})]
The statement is an immediate consequence of the functional equation for
$T_y$-genus~\cite{Hirzebruch, Theorem~11.3.1} and the well known fact
that $h^{p,q}(X)=0$ unless $p=q$ or $p+q=n$.
\endproof

In small dimensions (surfaces in~$\Cp{q}$ and hypersurfaces in~$\Cp5$)
Proposition~\ref{C1.3.D} gives the following (to simplify the formulas we
denote by~$\mu_i$ the $i$-th elementary symmetric polynomial in
$m_1,\dots,m_s$; certainly, $\mu_0=1$ and $\mu_i=0$ for $i>s$).

\claim[Corollary]\label{C1.3.E}
If $X\subset\Cp{q}$, $\dim X=2$, then
$$
\gather
h^{1,1}(X)=
\tfrac1{12}\mu_{q-2}\bigl(8\mu_1^2-10\mu_2-6(q+1)\mu_1+(q+1)(3q-2)\bigr),\\
h^{1,1}(Y)=\mu_{q-2}k\bigl(3k+\mu_1-(q+1)\bigr)+2h^{1,1}(X);
\endgather
$$
If $X\subset\Cp5$, $\dim X=4$, then
$$
\align
\text{$s=0$:}\quad&
 h^{2,2}(Y)=\tfrac{115}{12}k^4-\tfrac{115}6k^3+
  \tfrac{185}{12}k^2-\tfrac{35}6k+2,\\ 
\text{$s=1$:}\quad&
 h^{2,2}(X)=\tfrac{11}{20}m_1^5-\tfrac{11}4m_1^4+
  \tfrac{23}4m_1^3-\tfrac{25}4m_1^2+\tfrac{37}{10}m_1,\\ 
&h^{2,2}(Y)=
\tfrac{115}{12}m_1k^4+\left(\tfrac{23}6m_1-23\right)m_1k^3\\
&\qquad\qquad+\left(\tfrac{13}6m_1^2-\tfrac{21}2m_1+
\tfrac{95}4\right)m_1k^2\\
&\qquad\qquad+\left(\tfrac{11}{12}m_1^3-5m_1^2+
\tfrac{43}4m_1-\tfrac{25}2\right)m_1k+2h^{2,2}(X)\rlap\qed.
\endalign
$$
\endclaim

\remark{Remark}
Originally the bound in the Petrovsky-Oleinik inequality was expressed
via the number of certain monomials rather than in terms of the Hodge
structure. The coincidence of the two bounds was proved by
V.~Zvonilov~\cite{ZvonPetr} using~\ref{C1.3.B}.
It is closely related to\mnote{toZv: slightly modified} the following
beautiful rule: if $X\subset\Cp{2k+1}$ is a hypersurface, $\deg X=d$, then
$h^{k,k}(X)-1$ is the number of the integral points in the layer
$kd<\sum(x_i+1)<(k+1)d$ of the cube $[0,d-2]^{2k+1}\subset\R^{2k+1}$.
Other layers of the cube give the other Hodge numbers. The formula
extends to hypersurfaces in spaces of even dimension. Relations of this
type were first observed by, probably, J.~Steenbrink~\cite{Steen}, who
found an explicit monomial basis for the vanishing Hodge structure of a
semiquasihomogeneous singularity. For further generalizations
see~\cite{KhoDanilov} and~\cite{VarHodge}.
\endremark


\subsection{Types of involutions}\label{C1.4}
In the theory of real curves the notion of separating (or type~$\I{}$,
see Introduction) curve is crucial. In higher dimensions, when $\RX$ has
(real) codimension greater than one, it cannot separate~$X$, but the
homology vanishing condition, $[\RX]=0\in H_n(X)$, still makes sense. In
various applications the properties of~$\RX$ may depend on whether it is
homologous to a certain distinguished class $u\in H_n(X)$ which is
usually natural (say, a characteristic class of~$X$), but does not need to
be zero.  Thus, given a class $u\in H_n(X)$, we say that the real
structure $\conj$ on~$X$ (or just $X$ itself, or~$\RX$) is of
\emph{type~$\I{u}$} if $[\RX]$ is homologous to~$u$. Here are the special
cases used most commonly:
\roster
\item
$\I0$ or $\Iabs$, with $u=0$. Type~$\I0$ is sometimes called the
\emph{absolute type~$\I{}$};
\item
$\Iwu$, where $u=u_n(X)$ is the  $n$-th Wu class of~$X$. (Recall that
$u_n(X)$ is the characteristic element of the intersection form of~$X$.
Due to the Wu formula, $u_n(X)$ is a certain polynomial in the
Stiefel-Whitney classes of~$X$\, \footnote{Since $X$ is a complex
manifold, $w_{2i+1}(X)=0$ and $w_{2i}(X)=c_i(X)\bmod2$. Note that
$\conj_*c_i=(-1)^ic_i$; in particular, all the $\!\bmod2$ characteristic
classes of~$X$ are $\conj_*$-invariant.});
\item
$\Ihp$, with $u$ equal to the $(n/2)$-th power of the hyperplane
section.  (Certainly, this only makes sense for a fixed embedding
$X\hookrightarrow\Cp{N}$.)
\endroster

In the case of projective varieties the notion of types~$\I0$ and~$\Ihp$
appeared in~\cite{Kh75}. In~\cite{Viro:survey} types~$\I0$ and~$\Ihp$ are
called~$\Iabs$ and~$\Irel$ respectively, and a hypersurface which is not
of one of these types is said to have type~$\II$. However, this
terminology does not seem commonly accepted. Usually, depending on a
particular problem, it is reasonable to distinguish a class $u\in H_n(X)$
and consider manifolds of types~$\Iabs$ and $\I{u}=\Irel$, regarding the
rest as type~$\II$. E.g., in most results below it is type~$\Iwu$ that
plays essential r\^ole.

\subsection{Congruences}\label{C1.5}
Next two theorems are direct generalizations of Gudkov, Ar\-nol$'$d,
Rokhlin, Kharlamov, and Krakhnov congruences.

\claim[Extremal congruences]\label{C1.5.A}
Let $X$ be an even dimensional $(M-d)$-manifold with real structure.
Then\rom:
\roster
\item\local1
if $d=0$, then $\chi(\RX)=\sigma(X)\bmod{16}$\rom;
\item\local2
if $d=1$, then $\chi(\RX)=\sigma(X)\pm2\bmod{16}$\rom;
\item\local3
if $d=2$ and $\chi(\RX)=\sigma(X)+8\bmod{16}$, then $X$ is of type~$\Iwu$.
\endroster
\endclaim

\claim[Generalized Arnol$'$d congruence]\label{C1.5.B}
Let $X$ be a manifold with real structure of type\/~$\Iwu$. If\/ $\dim X$
is even, then $\chi(\RX)=\sigma(X)\bmod8$.
\endclaim

\claim[Proposition]\label{C1.5.C}
An $M$-manifold with real structure is of type~$\Iwu$.
\endclaim

The signature of a complex analytic manifold~$X$ equals~$T_1(X)$
(see~\ref{C1.3}). If $X$ is a surface, this gives
$\Gs(X)=\13(c_1^2-2c_2)[X]$. If $X$ is a regular complete intersection
in~$\Cp{q}$, its signature can be found using Proposition~\ref{C1.3.D}
and~\eqref{eqC.3}. Here are some partial results:

\claim[Signature of a complete intersection]\label{C1.5.D}
Let $X$ be a regular complete intersection in~$\Cp{q}$ of polydegree
$(m_1,\dots,m_s)$ and even dimension $n=q-s$.
\roster
\item\local1
If $n=2$, then $\Gs(X)=-\tfrac13\mu_{q-2}(\mu_1^2-2\mu_2-q-1)$, where
$\mu_i$ is the $i$-th elementary symmetric polynomial in
$(m_1,\dots,m_s)$.
\item\local2
If $s=q+2\bmod4$ and $m_1+\dots+m_s=q+1\bmod2$, then $\Gs(X)=0\bmod{16}$.
\item\local3
If $s=q\bmod4$ and all~$m_i$ but, maybe, one \rom(say,~$m_s$\rom) are odd,
then $\Gs(X)=m_1\ldots m_s\bmod{16}$.
\item\local4
If $X$ is a hypersurface \rom(i.e., $s=1$\rom), then
\item""
$\Gs(X)=m_1\pmod{16}$ if $n=0\bmod4$\rom;
\item""
$\Gs(X)=0\bmod{16}$ if $n=2\bmod4$ and $m_1$ is even\rom;
\item""
$\Gs(X)=1-m_1(m_1-1)\bmod{16}$ if $n=2\bmod4$ and $m_1$ is odd.
\endroster
\endclaim

\proof
We use the identity $\Gs(X)=\chi_s^q(m_1,\dots,m_s;1)$. Statement~\loccit1
is proved by induction using~\eqref{eqC.3}. Statement~\loccit2 follows
from the fact that the signature of a $(8k+4)$-dimen\-sional
$\Spin$-manifold is divisible by~16 (see~\cite{Ochanine}). To
prove~\loccit3, we proceed by induction in~$s$ and~$m_s$ and use the
functional equation for $T_y$-genus \cite{Hirzebruch, Theorem~11.3.1)},
which implies
$$
\chi_s^q(\dots,m_s;1)=\chi_s^q(\dots,m_s-1;1)+
\chi_{s-1}^{q-1}(\dots;1)-\chi_{s+1}^{q-1}(\dots,m_s-1,m_s;1).
\eqtag\label{eqC.4}
$$
Under the hypotheses of~\loccit3 the last term in~\eqref{eqC.4} vanishes
due to~\loccit2.

Statements~\loccit2 and~\loccit3 cover the first two cases in~\loccit4. In
the last case we use~\eqref{eqC.4} again; now
$\chi_{1}^{q}(m-1;1)=0\bmod{16}$ due to~\loccit2,
$\chi_{0}^{q-1}(1)=\Gs(\Cp{q-1})=1$, and
$\chi_{2}^{q-1}(m-1,m;1)=(m-1)m\bmod{16}$ due to~\loccit3.
\endproof

The proofs of~\ref{C1.5.A}--\ref{C1.5.C} are similar to each other; the
key ingredient is Lem\-ma~\ref{C1.5.E} below.

Let $n=\dim X$.  Denote $H=H_n(X;\Z)/\Tors$ and let $H^{\pm1}\subset H$ be
the eigensubgroups of~$\conj_*$ and $J=H/(H^{+1}\oplus H^{-1})$.  Denote
by~$B\:H\otimes H\to\Z$ the intersection form $(x,y)\mapsto x\circ y$
of~$X$, and by~$B^{\conj}$, the twisted intersection form $(x,y)\mapsto
x\circ\conj_*y$, see~\ref{involutions}. Let $q^\pm$ be the discriminant
quadratic space associated with $B|_{H^{\pm}}$. As is known
(see~\ref{inv.lattice}), $J$ and $q^\pm$ are $\Z_2$-vector spaces
isomorphic as groups.

\claim[Lemma]\label{C1.5.E}
Let,\mnote{toZv:new title} as above, $X$ be an $(M-d)$-manifold with real
structure, $\dim X=n$. Then\rom:
\roster
\item\local1
the image of $[\RX]$ in $H_2(X)$ is the characteristic element of the
twisted intersection form~$B^{\conj}$\rom;
\item\local2
both $u_n(X)$ and $[\RX]\in H_n(X)$ are integral classes, i.e., they
belong to the image of $H_n(X;\Z)$\rom;
\item\local3
$d\ge\dim\Coker[(1+\conj_*)\:H_*(X)\to H_*(X)]\ge\dim J$.
\endroster
Let, besides, $n=\dim X$ be even, $n=2k$. Denote $\epsilon={(-1)^{k+1}}$.
Then
\roster
\item[4]\local4
$\chi(\RX)=(-1)^{k+1}\sigma(B^{\conj})=\sigma(X)-2\sigma(H^\epsilon)$\rom;
\item\local5
$H^\epsilon$ is an even lattice\rom;
\item\local6
$\dim J=d\bmod2$.
\endroster
\endclaim

\proof
Statement~\loccit1 is the Arnol$'$d lemma~\ref{Arnold.lemma}. \loccit4 is
an immediate consequence of \ref{Hst}, and \loccit5 follows from the fact
that $u_n(X)$, like any characteristic class of~$X$, is realized by a
formal linear combination of real cycles
of dimension~$k$ (see, e.g., \cite{Ro4}) and, hence, belongs
to~$H^{-\epsilon}$.

To prove~\loccit2 note that both the standard and twisted $\Z_2$-valued
intersection forms vanish on the image of $\Tors H_n(X;\Z)$ (as so do
their $\Z$-valued counterparts). Hence, their characteristic classes
annihilate $\Tors H_n(X;\Z)$ and thus are integral.

Statement~\loccit3 follows from~\iref{1.1.2}1 and the construction of the
Smith exact sequence given in~\ref{1.1.9}. Indeed, $\Ker\pr_*$
in~\ref{Smith.seq} consists of geometrically invariant, and hence
$\conj_*$-invariant cycles, which must belong to $\Ker(1+\conj_*)$.  This
gives the first inequality in~\loccit3; the second one is obvious, as
$\dim J$ equals $\dim\Ker(1+\conj_*)$ restricted to $H\otimes\Z_2$.

The proof of~\loccit6 is based on the following two facts:

\claim[Lemma]\label{C1.5.F}
For any closed manifold~$M$ of even real dimension $n=2k$ one has
$\beta_*(M)=(-1)^{k}\chi(M)\bmod4$.
\endclaim

\proof
If $k=2l$~is even, by Poincar\'e duality
$\beta_*(M)-\chi(M)=2\sum_{r=0}^{k-1}\beta_{2r+1}(M)=
4\sum_{r=0}^{l-1}\beta_{2r+1}(M)$. The proof for $k$~odd is similar.
\endproof

\claim[Lemma]\label{C1.5.G}
For any closed complex manifold~$X$ of even complex dimension $n=2k$ one
has $\sigma(X)=(-1)^k\chi(X)\bmod4$.
\endclaim

\proof
The proof is similar to the previous one, using the $T_y$-characteristic
(see~\ref{C1.3}), the identities $\sigma(X)=T_1(X)$ and
$\chi(X)=T_{-1}(X)$, and the symmetry $T_{n,i}=(-1)^nT_{n,n-i}$ of the
coefficients~$T_{n,i}$ in $T_y(X)=\sum T_{n,i}y^i$ (see \cite{Hirzebruch,
Section~1.8}; note that $T_{n,i}$ take integral values on any, not
necessarily K\"ahler, complex manifold.)
\endproof

Lemmas~\ref{C1.5.F} and~\ref{C1.5.G} applied to~$X$ give
$\beta_*(X)=(-1)^k\sigma(X)\bmod4$. Then \loccit4 turns into
$(-1)^k\beta_*(\RX)=(-1)^k\beta_*(X)-2\sigma(H^\epsilon)\bmod4$, and
\loccit6 follows from the congruence $\sigma(H^\epsilon)=\dim
q^\epsilon\bmod2$ (the lattice is even due to~\loccit5) and the identity
$\dim J=\dim q^\pm$.
\endproof

\remark{Remark}
Lemma \ref{C1.5.G} extends to almost complex manifolds; one can use either
Atiyah-Dupont formula~\cite{AD} (see G.~Wilson~\cite{Wilson}) or
cobordism arguments.  \endremark

\proof[Proof of Theorem \ref{C1.5.A}]
Under the hypotheses of \iref{C1.5.A}1 or~\ditto2 from \iref{C1.5.E}3
and~\ditto6 it follows that $\dim q^\epsilon=\dim J= 0$ or~$1$,
respectively (where still $\epsilon=(-1)^{k+1}$), and Theorem~\ref{10.2.B}
applied to $H^\epsilon$ in~\iref{C1.5.E}4 gives \iref{C1.5.A}1
and~\ditto2.

Similar arguments show that under the hypotheses of~\ref{C1.5.A}(3) the
discriminant forms $q^{\pm}$, whose dimensions are at most~2, are even
and, due to~\ref{10.4.C}, the difference $[\RX]-u_n(X)$ annihilates all
the integral classes in $H_n(X)$, i.e., belongs to the image of $\Tors
H_n(X;\Z)$. Now, using again \iref{C1.5.E}3, where both the inequalities
turn into equalities, one concludes that each element of
$H_n(X)/H_n(X;\Z)\otimes\Z_2$ has a $\conj_*$-invariant representative.
Hence, $B$ and $B^{\conj}$ coincide on such classes and $u_n(X)-[\RX]$
vanish.
\endproof

\proof[Proof of Theorem~\ref{C1.5.B}]
Since $[\RX]$ and $u_n(X)$ coincide, from~\ref{10.4.C} it follows that
$q^\pm$ are even discriminant lattices. Hence, $\Br q^\pm=0\bmod4$, and
\iref{C1.5.E}4 applies.
\endproof

\proof[Proof of Proposition~\ref{C1.5.C}]
The statement follows from~\iref{C1.5.E}1 and~\iref{C1.5.E}3, which
implies that $\conj_*$ acts as identity on $H_*(X)$.
\endproof

\subsection{Additional prohibitions}\label{C1.6}

\claim[Generalized Arnol$'$d inequality]\label{C1.6.A}
Let~$X$ be
a closed complex K\"ahler manifold with real structure and
$\dim X=n=2k$ even. Then the number~$p_-$ of the orientable components
of~$\RX$ with negative Euler characteristic satisfies the inequality
$$
p_-\le\frac14\bigl({b_n(X)-(-1)^k\sigma(X)}\bigr)+
\frac12\sum_{j=1}^k(-1)^jh^{k-j,k-j}(X)+\frac12\bigl({(-1)^k-1}\bigr).
$$
\endclaim

\proof
We proceed exactly as in~\cite{Arnold:break} and derive the inequality from
$p_-\le\sigma_{\epsilon}^{+1}$, where $\epsilon=(-1)^{k+1}$ and
$\sigma_\pm^{+1}$ are the inertia indices of the intersection form of~$X$
restricted to~$H^{+1}_n(X;\Z)$. In order to find~$\sigma_\epsilon^{+1}$,
we use the relations
$$
\align
\Gs_+^{+1}+\Gs_-^{+1}+\Gs_+^{-1}+\Gs_-^{-1}&=b^n(X),\\
\Gs_+^{+1}-\Gs_-^{+1}+\Gs_+^{-1}-\Gs_-^{-1}&=\Gs(X),\\
\Gs_+^{+1}-\Gs_-^{+1}-\Gs_+^{-1}+\Gs_-^{-1}&=(-1)^k\chi(\RX),\\
\Gs_+^{+1}+\Gs_-^{+1}-\Gs_+^{-1}-\Gs_-^{-1}
&=\Trace\bigl(\conj^*,H^n(X;\Z)\bigr).
\endalign
$$
In the right hand side of the last equation $H^n(X;\Z)$ can be replaced
with $H^{k,k}(X)$, and an estimate on~$\sigma_\epsilon^{+1}$ follows from
comparing
$$
\Trace\bigl(\conj^*,H^{k,k}(X)\bigr)=
\sum_{j=0}^{k}(-1)^j\Trace(\conj^*,P^{k-j,k-j}(X))
$$
with~\eqref{eqC.1} and the obvious inequality
$$
-\Trace(\conj^*,P^{i,i}(X))\le \dim P^{i,i}(X)=h^{i,i}(X)-h^{i-1,i-1}(X).
$$
(In the case of $k$ odd one can also use the fact that
$-\Trace(\conj^*,P^{0,0}(X))$ contributes to~$\sigma_\epsilon^{+1}$ and,
on the other hand, $\conj^*=\id$ on $P^{0,0}(X)$.)
\endproof

\subsubsection{}
For a subset $S\subset\Rp{N}$ let
$$
\ell(S)=\max\bigl\{-1,i\bigm|\text{$\inc^*\: H^i(\Rp{N})\to H^i(S)$ is
nontrivial}\bigr\},
$$
where $\inc\:S\hookrightarrow\Rp{N}$ is the inclusion. If
$A\subset\Cp{N}$ is a real projective variety, we let $\ell(A)=\ell(\RA)$.
It is clear that $\inc^*\: H^i(\Rp{N})\to H^i(\RA)$ is a nontrivial
homomorphism for all $i\le\ell(A)$ and that $\ell(A)=n=\dim A$ if $\deg
A$ is odd. (Recall that the degree of~$A$ is the number of intersection
points of~$A$ and a generic $(N-n)$-plane in~$\Cp{N}$.) The following
simple statement is a direct consequence of the Poincar\'e duality and the
standard exact sequences.

\claim[Proposition]\label{l-prop}
If $\deg A$ is even, then $\ell(\RP_+)=\ell(A)\le\ell(\RP_-)=n-l(A)$ and,
in particular, $\ell(A)\le\frac{n}2$.
\qed
\endclaim

\remark{Remark}
As it follows from~\ref{l-prop} and the remark after~\ref{C1.2.C}, if
$\deg A$ is even, the double covering of $\Cp{n+1}$ is an $M$-variety
(under a proper choice of the covering real structure) if and only if $A$
is an $(M-d)$-variety, where for $n$ even
$$
d=\cases
 n-2\ell(\RA),& \text{if $m=0\bmod 4$},\\
 1,&         \text{if $m=2\bmod 4$ and $n>2\ell(\RA)$, or}\\
 0,&         \text{if $m=2\bmod 4$ and $n=2\ell(\RA)$.}
\endcases
$$
and for $n$ odd
$$
d=\cases
 n-2\ell(\RA)-1,& \text{if $m=0\bmod 4$},\\
 0,&         \text{if $m=2\bmod 4$.}
\endcases
$$
Another, related, consequence of the same calculation is an improvement
of the Smith-Thom bound: {\proclaimfont $A$ can not be an $(M-l)$-variety
with $l<d$ given above.}
\endremark

The invariant $\ell(\RA)$ was introduced in~\cite{Kh75}; the results
of~\cite{Kh75} were later improved by I.~Kalinin~\cite{Kalinin}.

\claim[Additional extremal
congruences for $(M-d)$-hypersurfaces]\label{C1.6.B}
Let
$A\subset\Cp{n+1}$ be a real $(M-d)$-hypersurface of even dimension~$n$
and degree~$m$.  If $m=\pm2\bmod8$ and $n/2-\ell(A)$ is odd, then
\roster
\item if $d=n/2-\ell(A)$, then
$\chi(\RA)=\Gs(A)\mp2\bmod{16}$\rom;
\item
If $d=n/2-\ell(A)+1$, then $\chi(\RA)=\Gs(A),\,\Gs(A)\mp4\bmod{16}$.
\endroster
If $m=0\bmod8$ and $\ell(A)<n/2$, then
\roster
\item[2]
if $d=n-2\ell(A)$, then $\chi(\RA)=\Gs(A)\bmod{16}$\rom;
\item
if $d=n-2\ell(A)+1$, then $\chi(\RA)=\Gs(A)\pm2\bmod{16}$.
\endroster
\endclaim

\proof[Proof of Theorem~\ref{C1.6.B}\/ \ROM(see~\cite{Kalinin})]
Without going too deep into the details, we just outline the principal
ideas.  Like the other extremal congruences, \ref{C1.6.B}~is derived
from~\iref{C1.5.E}4 and~\ditto5 using Theorem~\ref{10.2.B}. The crucial
point, which replaces~\iref{C1.5.E}3 and~\ditto6 and gives an estimate on
$\dim J=\dim q^\pm$, is the following lemma:

\claim[Lemma\/ \ROM(see~\cite{Kalinin})]\label{C1.6.C}
Let $A\subset\Cp{n+1}$ be a real $(M-d)$-hypersurface of dimension~$n$ and
degree~$m$.
\roster
\item
If $\ell(A)\ge\intpart{(n-1)/2}$, then $\dim J=d$.
\item
If $\ell(A)<\intpart{(n-1)/2}$, then
\endroster
$$
\dim J=d-
\cases
2\bigl(\intpart{(n-1)/2}-\ell(A)\bigr)&\text{if $m=0\bmod4$,}\\
2\intpart{(n/2-\ell(A))/2}&\text{if $m=2\bmod4$ and $n$ is even,}\\
(n-1)/2-\ell(A)&\text{if $m=2\bmod4$ and $n$ is odd.}
\endcases
$$
\rom(Here $\intpart{x}$ denotes the integral part of~$x$.\rom)
\endclaim

Lemma~\ref{C1.6.C} is proved using
Kalinin's\mnote{\Dg: edited after Zv} spectral sequence $(\rH^*,\rd^*)$,
see~\ref{C1.7}. If $A$ is a regular complete intersection in~$\Cp{N}$,
the difference $d-\dim J$ is equal to the number of nontrivial
differentials~$\rd^i$ with $r>1$. Lemma~\ref{C1.6.C} is obtained
in~\cite{Kalinin} from an explicit calculation of $\rH^*(\Cp{N})$ and
$\rH^*(\Cp{n+1}\sminus A)$. This calculation also gives the other key
ingredient of the proof of~\ref{C1.6.B}, which is stated below. Denote
by~$h$ the generator of $H^2(\Cp{N})$ and let $\Gd_h(A)=0$ if $n=\dim A$
is odd or $n$~is even and $\inc^*h^{n/2}\in\Im(1+\conj^*)$ in
$H^n(A;\Z_2)$, and $\Gd_h(A)=1$ otherwise.

\claim[Lemma\/ \ROM(see~\cite{Kalinin})]\label{C1.6.D}
Let $A\subset\Cp{n+1}$ be a real hypersurface of even dimension~$n$ and
degree~$m$. If $m=2\bmod4$, then
$$
\align
\Gd_h(A)=&
\cases
0&\text{if $\ell(A)<n/2$ and $\ell(A)=n/2-1\bmod2$,}\\
1&\text{otherwise\rom;}
\endcases\\
\noalign{\medskip\noindent if $m=0\bmod4$, then
\medskip}
\Gd_h(A)=&
\cases
0&\text{if $\ell(A)<n/2$,}\\ 1&\text{if $\ell(A)\ge n/2$.}\\
\endcases
\endalign
$$
\endclaim

From the lemma it follows that under the hypotheses of~\ref{C1.6.B} one
has $\inc^*h^{n/2}\in\Im(1+\conj^*)$. This gives additional information
on~$\Br q^\pm$ when $\dim q^\pm=\dim J$ is small (see~\ref{smallBr}).
\endproof

\subsection{Inherited structures on the real part}
In
the rest of this section we briefly discuss a few
constructions generalizing, to an extent, the notion of complex
orientation of a dividing real curve.

\claim[$\Spin$-orientations]\label{Spin.orient}
Let $X$ be a manifold with real structure~$\conj$. Than any
$\conj$-invariant $\Spin$-structure on~$X$ defines, in a natural way, a
semi-orientation \rom(i.e., pair of opposite orientations\rom) on~$\RX$.
\endclaim

\remark{Remark}
If $X$ is $\Spin$ and $H_1(X)=0$, the only $\Spin$-structure on~$X$ is
obviously $\conj$-invariant and, hence, $\RX$ has a canonical
semi-orientation. In particular, it is orientable (cf\.~\ref{C2.2.H}).
The orientability statement extends to involutions on arbitrary manifolds
and is known as Edmonds theorem \cite{Edmonds}.
\endremark

\proof[Construction]
One needs to compare orientations at two points $x_1,x_2\in\RX$.
Connect~$x_1$, $x_2$ by a path~$\gamma$ in~$X$, represent the two
orientations by tangent $n$-frames $\Xi_1$, $\Xi_2$ (where
$n=\dim_{\C}X=\dim_{\R}\RX$), extend $(\Xi_1,\sqrt{-1}\,\Xi_1)$ and
$(\Xi_2,\sqrt{-1}\,\Xi_2)$ to a $2n$-frame field $(\xi_1,\dots,\xi_{2n})$
on~$\gamma$, and evaluate the chosen $\Spin$-structure on the loop
$\gamma\circ\conj_*\gamma$, where $\conj_*\gamma$ is framed with
$$
(\conj^*\xi_1,\dots,\conj^*\xi_n,-\conj^*\xi_{n+1},\dots,-\conj^*\xi_{2n}).
$$
The two orientations are regarded coherent if the resulting value is~$0$.
It is straightforward to check that, if the $\Spin$-structure is
$\conj$-invariant, the result does not depend on the choices
made.
\endproof

\subsubsection{Stiefel orientations}\label{Stiefel}
Recall that an orientation of a smooth manifold~$Y$ can be defined as a
homotopy class of lifts to~$\BSO$ of a classifying map $f_Y\:Y\to\BO$ of
the tangent bundle of~$Y$. The fibration $\BSO\to\BO$ can, in turn, be
regarded as the $K(\ZZ,0)$ fibration killing $w_1\in H^1(\BO)$.
Generalizing this approach one can fix a characteristic class $\omega\in
H^{i+1}(\BO)$ and define $\omega$-structures on~$Y$ as homotopy classes of
lifts of~$f_Y$ to a $K(\ZZ,i)$-fibration $\BO_\omega\to\BO$
killing~$\omega$. It is easy to see that $Y$ admits an $\omega$-structure
if and only if $\omega(Y)=0$ and, if nonempty, the set of
$\omega$-structures on~$Y$ is an affine space over $H^{i}(Y)$ (see,
e.g.,~\cite{Degt:ks}). Thus, orientations and $\Spin$- (more precisely,
$\Pin^+$-) structures are, respectively, $w_1$- and $w_2$-structures.

If $\omega=w_{i+1}$ is a Stiefel-Whitney class, one can replace
$\BO_\omega$ with the corresponding associated Stiefel bundle and show
that a $w_i$-structure can be regarded as a $\ZZ$-valued function on the
homotopy classes of $(n-i)$-framed $i$-cycles (where $n=\dim Y$).  This
gives rise to the following generalization of the notion of complex
orientation of a dividing curve (see~\cite{Viro.???} or, for a more formal
approach,~\cite{Degt:Stiefel}). Let $X$ be a complex manifold with real
structure, $\dim X=n$, and $[\RX]$ vanishes in $H_n(X)$. Then $\RX$
possesses a canonical $w_n$-structure, whose value on a $1$-framed
$(n-1)$-cycle~$\gamma$ is defined as the linking coefficient of~$\RX$ and
a shift of~$\gamma$ along the framing multiplied by~$\sqrt{-1}$. (More
precisely, this is a partial structure, defined on the kernel of the
inclusion homomorphism $H_{n-1}(\RX)\to H_{n-1}(X)$.) The assumption
$[\RX]=0$ assures that the linking coefficient is well defined. The
construction admits a further generalization. Let $\gamma$ be an
$(n-i)$-cycle and $(\xi_1,\dots,\xi_i)$ its framing. Proceed as in the
geometric construction of Viro homomorphisms (see~\ref{C1.7}) and
include~$\gamma$ into an $(n-1)$-cycle~$\gamma'$ in~$X$ tangent to
$\sqrt{-1}\,\xi_1$, \dots, $\sqrt{-1}\,\xi_{i-1}$; then shift~$\gamma'$
along an extension of $\sqrt{-1}\,\xi_i$ and evaluate the linking
coefficient of the shift and~$\RX$. A detailed analysis shows that this
construction gives rise to a partial $w_{n-i+1}$-structure on~$\RX$,
defined on the kernel of $\vr_{n-1}\:H_{n-i}(\RX)\to\RH{i}_{n-1}(X)$,
provided that $\vr_{n+i-1}[\RX]$ vanishes in $\RH{i}_{n+i-1}(X)$. In
particular, the last condition implies that
$w_n(X)=\ldots=w_{n-i+1}(X)=0$.

\end



\ifx\rokhloaded\undefined\input rokh.def \fi

{\catcode`\@=11
\gdef\proclaimfont@{\sl}
\global\let\proclaimfont\proclaimfont@}

\def\PD{\mathop{\roman D\null}}
\def\PV{\Cal P}
\def\1#1{\mathchoice{\smash{\displaystyle\frac1{#1}}}%
  {\smash{\textstyle\frac1{#1}}}{\smash{\scriptstyle\frac1{#1}}}%
  {\smash{\scriptscriptstyle\frac1{#1}}}}
\def\CC{\Cal C}
\def\CS{\Cal S}

\let\+\dplus

\def\tX{\tilde X}

\let\ix\index

\document


\section{Surfaces}\label{C2}

\subsection{Basic results}\label{C2.1}
In
this
section by \emph{complex surface} we mean a closed complex analytic
manifold of complex dimension~2. First, we restate some results of
Section~\ref{C1} in more topological terms.

\claim
[Theorem]\label{C2.1.A}
Let~$X$ be\mnote{toZv:title changed similar to r-tools} a complex surface
with real structure. Then
$$
 \beta_*(\RX)\le\beta_*(X)\quad\text{and}\quad
 \beta_*(\RX)=\beta_*(X)\bmod2.
$$
Let, further, $A$ be a nonsingular even ample real divisor on~$X$ and
$Y\to X$ a double covering branched over~$A$. Then
$$
\gather
\Gb_*(\RX_-)+\Gb_*(\RX_+)\le\Gb_2(X)+\Gb_*(A)\rlap{ and}\\\vspace{1\jot}
2\Gb_*(\RX_\pm)-4c_0(\RX_\pm)\le2\Gb_2(X)+\Gb_*(A),
\endgather
$$
where $\RX_\pm$ are the two halves of $\RX\smallsetminus\RA$
\rom(see~\ref{C1.1}\rom) and $c_0(\RX_\pm)$ is the number of closed
components of~$\RX_\pm$ covered nontrivially in~$\RY_\pm$.
\qed
\endclaim

According to the Nakai-Moishezon criterion (see, e.g.,~\cite{BPV}) an
effective divisor~$A$ on a nonsingular K\"ahler surface~$X$ is ample if and
only if $A\circ D>0$ for any irreducible curve $D\subset X$. In general,
the ampleness condition in~\ref{C2.1.A} can be replaced with a weaker
requirement that $A$ should be connected and the inclusion homomorphism
$H_1(A)\to H_1(X)$ should be onto. Under these hypotheses from the Smith
exact sequence it follows that $\Gb_1(Y)=\Gb_1(X)$.

\claim
[Theorem]\label{C2.1.B}
Let\mnote{toZv: title changed similar to r-tools} $X$ be a complex
K\"ahler surface with real structure. Then
$$
\bigl|\chi(\RX)-1\bigr|\le\Gs^-(X),
$$
where $\Gs^-(X)=\frac12[b_2(X)-\Gs(X)]$ is the negative inertia index of
the intersection form of~$X$. Let, further, $A\subset X$ be a nonsingular
ample even real divisor on~$X$ so that there is a real divisor~$E$ on~$X$
with $2E\sim A$. Then
$$ \gather
 \bigl|\chi(\RX_-)-\chi(\RX_+)\bigr|\le
  \Gs^-(X)-b_1(X)+1+\frac14{[A]^2}-\frac12{\chi(A)},\\ \noalign{\vglue2\jot}
 \bigl|2\chi(\RX_\pm)-1\bigr|\le
  2\Gs^-(X)-b_1(X)+1+\frac14{[A]^2}-\frac12{\chi(A)}.\rlap\qed
\endgather
$$
\endclaim

Note that here, as well as in many other statements about surfaces (at
least in this survey) the K\"ahler property is only used to assert that
$\conj_*$ is traceless on $H_1(X;\R)$. If $X$ is not K\"ahler (but still
complex analytic), there still is a
canonical subspace $H^{1,0}(X)\subset H^1(X;\C)$ (generated by the
classes of closed holomorphic forms) and one has
$H^{1,0}\cap\conj^* H^{1,0}= H^{1,0}\cap\bar H^{1,0}=0$ and either
$b_1(X)=2h^{1,0}(X)$ or $b_1(X)=2h^{1,0}(X)+1$. Thus, the worst that can
happen is that the trace of $\conj^*$ on $H^1(X;\R)$ is $\pm1$, and the
first inequality should be replaced with
$\bigl|\chi(\RX)-1\bigr|\le\sigma^-(X)+1$.
Alternatively, the K\"ahler property can be replaced, e.g., with the
requirement $b_1(X)=0$. (Then the result holds for flexible varieties as
well.)\mnote{\Dg: edited} A stronger condition $\Gb_1(X)=0$ would also
imply the existence of~$E$ as in the statement. Note also that the
inequalities in~\ref{C2.1.B} do not appeal to any covering; however, we
do not know whether the statement still holds if there is no {\bf real}
double covering branched in~$A$.

The ampleness condition in the last two inequalities can be
replaced with $b_1(Y)=b_1(X)$; otherwise $b_1(Y)-b_1(X)$ should
be added to the right hand sides of both the inequalities. Note
that, unlike~\ref{C1.3.A}, all the statements can be proved
topologically following the lines of~\cite{Arnold:break}. (In
order to prove the second inequality one should consider both the
involutions~$T_\pm$ on~$Y$ and evaluate the inertia
indices~$\smash{\Gs_\pm^{\pm1,\pm1}}$ of the restrictions of the
intersection form to the bi-eigenspaces~$H_2^{\pm1,\pm1}$. The
inequality would then follow from $\Gs_-^{\pm1,\mp1}\ge0$.)
Moreover, taking into account various classes realized by the
orientable components of~$\RX$ and~$\RY$, one can obtain refined
inequalities and their extremal properties. Here is an example:

\claim
[Theorem]\label{C2.1.C}
Let\mnote{toZv:title changed for uniformity} $X$ be a complex K\"ahler
surface with real structure. Then
$$
-\Gs^-(X)\le1-\chi(\RX)\le\Gs^-(X)-2p_+,
$$
where $p_+$ is the number of orientable components of~$\RX$ of positive
Euler characteristic. If $A\subset X$ is a nonsingular ample even real
divisor on~$X$ so that there is a real divisor~$E$ on~$X$ with $2E\sim A$,
then
$$
\chi(\RX_-)-\chi(\RX_+)\le \Gs^-(X)-b_1(X)+1+\14{[A]^2}-\12{\chi(A)}-2q_+,
$$
where $q_+$ is the number of orientable components of~$\RY_+$ of positive
Euler characteristic \rom(which, clearly, correspond to the components
of~$\RX_+$ with positive Euler characteristic for which the restriction of
the covering projection $\RY_+\to\RX_+$ is the orientation double
covering.\rom)
\qed
\endclaim

The extremal congruences~\ref{C1.5.A}, Arnol$'$d congruence~\ref{C1.5.B},
and Proposition~\ref{C1.5.C} (as well as Lemma~\ref{C1.5.E}) are
transferred to surfaces without changes. The generalized Arnol$'$d
inequality~\ref{C1.6.A} takes the following simpler form:

\claim[Arnol$'$d inequality for surfaces]\label{C2.1.D}
Let~$X$ be a K\"ahler surface with real structure. Then the number~$p_-$ of
the orientable components of~$\RX$ with negative Euler characteristic
satisfies the inequality
$$
p_-\le\12(\Gs^+(X)-1),
$$
where $\Gs^+(X)=\frac12[b_2(X)+\Gs(X)]$.
\qed
\endclaim

Below is another example of a refined statement specific for surfaces
(see~\cite{Nikulin:condition};
in fact, one proves that the class $[\RX]\in
H_2(X;\Z)$ is divisible by $2^{r+1}$). More examples of refined extremal
congruences taking into account particular additional properties of the
surface are found in~\ref{C2.2}.

\claim[Nikulin's congruence]\label{C2.2.F}
Let $X$ be an $M$-sur\-face with real structure and $H_1(X)=0$. Suppose
that $\RX$ is orientable and that the Euler characteristic of each
component of~$\RX$ is divisible by~$2^r$ for some~$r\ge\nobreak1$. Then
$\chi(\RX)=0\bmod{2^{r+3}}$.
\endclaim

If~$X$ is a regular complete intersection in~$\Cp{q}$, then the complex
ingredients of~\ref{C2.1.A}--\ref{C2.1.D} (i.e., the topological
invariants of~$X$ and~$A$) can easily be found by induction
using~\ref{C1.3.D}.

\claim\label{C2.1.E}
Let a surface~$X$ be a regular complete intersection in~$\Cp{q}$ of
polydegree $(m_1,\dots,m_{q-2})$. Then
$$
\gather
 b_1(X)=0,\qquad b_2(X)=\chi(X)-2,\\ \vspace{1\jot}
 \chi(X)=\mu_{q-2}\left(\mu_1^2-\mu_2-(q+1)\mu_1+\tfrac12q(q+1)\right),\\
 \vspace{1\jot} \Gs(X)=-\tfrac13\mu_{q-2}(\mu_1^2-2\mu_2-q-1),
\endgather
$$
where $\mu_i$ is the $i$-th elementary symmetric polynomial in
$(m_1,\dots,m_{q-2})$. If, besides, $A$ is a nonsingular curve cut on~$X$
by a hypersurface of degree~$m$, then
$$
\chi(A)=-m\mu_{q-2}(m+\mu_1-q-1).\qed
$$
\endclaim

\subsection{The orbit space of the complex conjugation}
Another phenomenon specific for surfaces is the fact that the fixed point
set of the complex conjugation has codimension~$2$. Hence, the quotient
$X/\!\conj$ is a manifold; moreover, one can easily see that, up to
isotopy, there is a unique smooth structure on $X/\!\conj$ such that the
projection $X\to X/\!\conj$ is a double covering branched over $\RX$.

The resulting $4$-manifolds~$X/\!\conj$ form a very interesting class. On
one hand, they are closely related to algebraic surfaces; on the other
hand, in many respects their properties are just opposite to those of
algebraic surfaces. S.~Akbulut conjectured that the Seiberg-Witten
invariants of the quotient vanish whenever the positive inertia index of
the quotient is at least $2$. There is a strong evidence: if $\RX$ has a
component of genus $\ge2$, the vanishing of the invariants follows
immediately from the `smooth version' of the adjunction inequality applied
to the image of the component in the quotient (for the appropriate version
of the inequality, going back to Kronheimer and Mrowka, see~\cite{Szabo});
if $\RX=\varnothing$, it was proved by Sh.~Wang~\cite{Wang}. (If
$\RX=\varnothing$, the projection $X\to X/\!\conj$ is an honest double
covering and one can easily control the behaviour of solutions to the
Seiberg-Witten equation.) Furthermore, in many cases (see below) the
quotient $X/\!\conj$ is \emph{completely decomposable}, i.e., splits into
connected sum of copies of $\Cp2$, $\overline{\Cp2\!}\,$, and $S^2\times
S^2$. (Recall that minimal algebraic surfaces are irreducible.) Thus, in
many cases one can assert that $X/\!\conj$ admits neither complex nor
symplectic structure. There is a remarkable exception: if $X$ is a
$K3$-surface, $X/\!\conj$ is diffeomorphic to either rational or Enriques
surface, see below.

According to Arnol$'$d, the following
result
should\mnote{toZv: edited} be attributed to Maxwell:

\theorem\label{maxwell}
The orbit  space $\Cp2/\!\conj$ is diffeomorphic to $S^4$.
\endtheorem

The proof indicated below is found, e.g., in~\cite{MarinonKuiperMassey}.
It is a little bit shorter than the other, more direct, proofs published
in~\cite{Kuiper;Massey;Arnold}. (Probably, except the one in
\cite{Massey}; as explained in~\cite{ArOnMaxwell}, this beautiful
explicit proof was essentially known to Maxwell; a generalization of the
statement to higher dimensions is also given in~\cite{ArOnMaxwell}). The
reason is the usage of the following Cerf theorem (see~\cite{Cerf}):
{\proclaimfont the group of diffeomorphisms of the $3$-dimensional sphere
coincides with the group of diffeomorphisms of the $4$-dimensional ball
\rom(in short, $\Gamma_4 = 0$\rom)}. The main advantage of this proof is
that it can be generalized to other real surfaces, often using the
Laundenbach-Poenaru theorem~\cite{LP}, extending the Cerf theorem to
handlebodies.

\proof
Pick a real flag $\Cp0\subset\Cp1\subset\Cp2$. Its quotient by~$\conj$
gives a simple cell decomposition of the quotient space: a $4$-ball is
attached to a $2$-ball. The quotient of a closed round $4$-ball centered
at the origin of $\Cp2\sminus\Cp1$ is a $4$-ball, and the closure of its
complement is a closed regular neighborhood of the $2$-ball $\Cp1/\!\conj$.
As is known, a closed regular neighborhood of a smooth $2$-ball in a
smooth $4$-manifold is a smooth $4$-ball. Thus, two $4$-balls are patched
together and, according to the Cerf theorem, the total space is
diffeomorphic to~$S^4$.
\endproof

Most generalizations of~\ref{maxwell} state that, under certain
assumptions on the surface (usually, for surfaces explicitly constructed
in a certain way so that one can control the topology of the quotient) the
orbit space is completely decomposable. However, there are also a few
examples of surfaces whose quotient is simply connected but not
completely decomposable. Below we give a brief account of the known
results.

For nonempty quadrics and cubics in $\Cp3$ the decomposability result was
obtain by M. Letizia~\cite{Leti}. (Due to~\ref{maxwell}, a blow up at a
real point does not change the orbit space; hence, for quadrics and
connected cubics the result follows from~\ref{maxwell}.) For nonempty
quartics and, more generally, nonempty $K3$-surfaces the decomposability
was proved by S.~Donaldson~\cite{Donald}. (Curiously, some $K3$-surfaces,
so called Fresnel-Kummer surfaces, are related to the Maxwell
electro-magnetism theory.) Donaldson's approach is based on changing the
complex structure in the twistor family and transforming the real
structure to a holomorphic involution reversing the holomorphic forms; the
quotient becomes a rational surface, which is obviously diffeomorphic to
$\Cp2\mathbin{\#}k\overline{\Cp2\!}\,$ or $S^2\times S^2$.

The complete decomposability of the quotient has been proved for double
planes branched over curves of degree~$2k$ whose real part consists of a
single nest of depth~$k$ (S.~Akbulut~\cite{Akbulut}) or over curves
obtained by small perturbation of union of real line
(S.~Finashin~\cite{Finashin1}). Finashin extended these results to double
quadrics branched over small perturbations of unions of generatrices and
to certain regular complete intersections in $\Cp{N}$ obtained
recursively by perturbation of unions. For a wider class of plane curves
he proved that the quotients of the corresponding double planes become
completely decomposable after adding several copies
of~$\overline{\Cp2\!}\,$.

The complete decomposability of the quotient holds for all real
structures on rational surfaces (the proof is based on the classification
of real rational surfaces, see~\ref{C2.4.A}) and for all real Enriques
surfaces (using a modified version of Donaldson's trick, see the end
of~\ref{C2.5}, Enriques surfaces can be reduced to rational ones).
See~\cite{Finashin2} for details.

An example of complex surface with real structure whose quotient is not
completely decomposable is constructed by S.~Finashin and E.~Shustin
\cite{FinShust}. The quotient is a simply connected $\Spin$-manifold
whose signature is not~$0$ (and, hence, the manifold cannot be
diffeomorphic to $k(S^2\times S^2)$). Currently it is not known whether
there are surfaces whose quotient by complex conjugation do not decompose
into connected sum of copies of~$\Cp2$, $\overline{\Cp2\!}\,$, $S^2\times
S^2$, and $K3$-surface (the last option would allow $\Spin$-manifolds with
nonvanishing signature).

\subsection{Complex separation and Pontrjagin-Viro quadratic form}\label{C4.Mikh}
Let $X$ be a complex surface with real structure, $H_1(X)=0$, and $A$ a
real curve on~$X$. (We do not exclude the case $A=\varnothing$.) The pair
$(X,A)$ is said to be of \emph{characteristic type} if $[\RX]+[A]=w_2(X)$
in $H_2(X)$. Denote by $\pr\:X\to X'$ the projection to the quotient space
$X'=X/\!\conj$.\mnote{\Dg: extended after Zv} From the Smith exact
sequence of the deck translation involution and the projection formula
$w_2(X)=\pr^*w_2(X')+\PD^{-1}[\RX]$ it follows that $\rel w_2(X')=[A']$ in
$H_2(X',\RX)$ (where $A'=A/\!\conj$). Hence, $[\RA]=\partial w_2(X')$
vanishes in $H_1(\RX)$ and there is a surface $\RX_+^w\subset\RX$ such that
$\partial\RX_+^w=\RA$ and $\AA_+^w=\RX_+^w\cup A'$\mnote{$^w$ after Zv}
realizes $w_2(X')$. Let $\RX_-^w$ be the closure of $\RX\sminus\RX_+^w$.
Since $\RX$ is homologous to zero in~$X'$, the surface
$\AA_-^w=\RX_-^w\cup A'$\mnote{$A'$ after Zv}
also realizes $w_2(X')$. The decomposition
$\RX=\RX_+^w\cup\RX_-^w$ satisfying these conditions is unique as, again
due to the Smith exact sequence, $\dim\Ker[H_2(\RX)\to H_2(X')]=1$. It is
called the \emph{complex separation} of pair~$(X,A)$.

Since the \emph{Arnol$'$d surfaces}~$\AA_\pm^w$ are characteristic and
$H_1(X)=0$, they possess canonical $\Pin^-$-structures and corresponding
Rokhlin-Guillou-Marin forms, see~\ref{RGMform}. Denote these forms
by~$\qq_\pm^w$ and their restrictions to~$\RX_\pm^w$, by~$q_\pm^w$. It is
easy to see that $q_\pm^w$ vanishes on the kernel of the inclusion
homomorphism $H_1(\RX_\pm^w)\to H_1(\RX)$.

The following results are due to G.~Mikhalkin~\cite{Mikhalkin}.

\claim[Lemma]\label{C4.Mikh.A}
Under the above assumptions one has
$$
 \chi(\RX_-^w)=\14[\chi(\RX)-\Gs(X)]+\12[A]^2+\Br\qq_-^w\bmod8.
$$
\endclaim

\theorem\label{C4.Mikh.B}
Let $X$, $A$ be as above. Assume that $q_-^w$ vanishes on the image of
$H_1(\RA)$ and denote $R=\frac14\bigl(\chi(\RX)-\Gs(X)+[A]^2\bigr)+\Br
q_-^w$.
\roster
\item\local1
If $A$ is an $M$-curve, then $\chi(\RX_-^w)=R\bmod8$.
\item\local2
If $A$ is an $(M-1)$-curve, then $\chi(\RX_-^w)=R\pm1\bmod8$.
\item\local3
If $A$ is an $(M-2)$-curve and $\chi(\RX_-^w)=R+4\bmod8$, then $A$ is of
type~$\I{}$.
\item\local4
If $A$ is of type~$\I{}$, then $\chi(\RX_-^w)=R\bmod4$.
\endroster
\endtheorem

Note that Lemma~\ref{C4.Mikh.A} and Theorem~\ref{C4.Mikh.B} are of
topological nature and, hence, also hold for flexible curves on flexible
surfaces; $A$ may even be nonorientable.

\proof[Proof of~\ref{C4.Mikh.A} and~\Ref{C4.Mikh.B}]
Lemma~\ref{C4.Mikh.A} follows from the Rokhlin-Guillou-Marin
congruence~\cite{GM} and a straightforward calculation, which gives
$[\AA_-^w]\circ[\AA_-^w]=\frac12[A]^2-2\chi(\RX_-^w)$. Proof of
Theorem~\ref{C4.Mikh.B} repeats one of the proofs of the classical
extremal congruences for plane curves (see~\cite{GM.zamecha}); it is based
on~\ref{C4.Mikh.A}, additivity of the Brown invariant, and the following
obvious facts: if $A$ is an~$(M-d)$-curve, then $\dim H_1(A')=d$, and
$A'$ is orientable if and only if $A$ is of type~$\I{}$.
\endproof

The following fact is an immediate consequence of the definition of the
Rokhlin-Guillou-Marin form and~\ref{C4.Mikh.B}:

\claim[Proposition]\label{C4.Mikh.C}
The restrictions of~$\qq_\pm^w$ to $H_1(A')$ differ by~$2\Go$, where $\Go$
is the characteristic class of the covering $A\sminus\RA\to A'\sminus\RA$.
As a consequence, $q_\pm^w$ vanish on the image of $H_1(\RA)$
simultaneously, and in this case
$$
\Br q_+^w+\Br q_-^w=\12(\chi(\RX)+\Gs(X)-[A]^2)\bmod8
$$
and, if $A$ is of type~$\I{}$,
$$
\chi(\RX_+^w)-\chi(\RX_-^w)=\Br q_+^w-\Br q_-^w\bmod8.
$$
\rom(If flexible curves are admitted, in the first congruence $A$ must be
orientable\rom; the second one holds for nonorientable~$A$ as well.\rom)
\qed
\endclaim

The restrictions $q_\pm^w$ can be calculated using Kalinin's spectral
sequence (see~\ref{C1.7}). We will show that $q_\pm^w$ can be patched
together into a quadratic form defined on $H_1(\RX)$, which depends only
on the topology of~$(X,\conj)$ and the class $[A]\in H_2(X;\Z)$.

Fix a class $c\in H_2(X;\Z)$ so that $\conj_*c=-c$ and
$c=[\RX]+w_2(X)\bmod2$. Pick an element $x\in H_1(\RX;\Z)$. There is an
integral chain $x'$ in~$X$ such that $\partial x'=x$. (One may have to
multiply~$x$ by an odd integer, which is irrelevant in what follows.)
Then it is easy to see that $(1-\conj_*)x'$ is an integral cycle and
$$
\PV_c(x)=[(1-\conj_*)x']^2+c\circ(1-\conj_*)x'\bmod4
$$
does not depend on the choice of~$x'$. Furthermore, the function
$x\mapsto\PV_c(x)$ is a quadratic extension of the intersection form
$(x,y)\mapsto x\circ y$ and, since $c$ dies in~$\iH_2$, one has
$\PV_c(-x)=\PV_c(x)+2(c\circ\vr_2x)=\PV_c(x)$, i.e., $\PV_c$ factors
through $H_1(\RX)$. The function $\PV_c\:H_1(\RX)\to\Z_4$ is called the
\emph{Pontrjagin-Viro quadratic form}.

\proposition\label{C4.Mikh.D}
One has
\roster
\item\local1
$\PV_{c+2e}(x)=\PV_c(x)+2(e\circ\vr_2x)$ for any $e\in H_2(X;\Z)$,
$\conj_*e=-e$\rom;
\item\local2
$\Br\PV_c=\12[\chi(\RX)+\Gs(X)-c^2]\bmod8$\rom;
\item\local3
if $A$ is a \rom(flexible\rom) real curve on~$X$ such that $(X,A)$ is of
characteristic type, then $q_\pm^w$ are the restrictions of~$\PV_{[A]}$
to~$\RX_\pm^w$.
\endroster
\endproposition

\proof
\loccit1 follows immediately from the definition; \loccit3 is proved by
direct calculation of indices of membranes. To prove~\loccit2 it suffices
to construct an orientable flexible curve~$A$ on~$X$ such that the pair
$(X,A)$ is of characteristic type and $\PV_{[A]}$ vanishes on~$\RA$ (e.g.,
$\RA=\varnothing$); the statement for~$\PV_{[A]}$ would follow then
from~\loccit3 and~\ref{C4.Mikh.C}, and it would extend to the other
values of~$c$ using~\loccit1 and~\iref{Br}3. Such a curve can easily be
constructed: since $[\RX]+w_2(X)=0$ in $\iH_2$, there is a class $f\in
H_2(X;\Z)$ such that $(1-\conj_*)f=[\RX]+w_2(X)\bmod2H_2(X;\Z)$.
Realize~$f$ by a smooth orientable surface~$F$ transversal to~$\RA$ and
perturb the singularities of $F\cup(-\conj F)$ according to the
orientation.
\endproof

Mention two special cases of the Pontrjagin-Viro form. An integral
lift~$c$ of~$w_2(X)$ corresponds to a $\Spin^{\C}$-structure on~$X$.
Hence, if $[\RX]=0$ in $H_2(X)$, the correspondence $c\mapsto\PV_c$
assigns to any $\Spin^{\C}$-structure on~$X$ a $\Pin^-$-structure
on~$\RX$ (cf@.~\ref{Stiefel}; the correspondence can be extended to all,
not necessarily $\conj_*$-skew-invariant classes~$c$
via~\iref{C4.Mikh.D}1). There also is a similar correspondence in the
case when $\RX$ is orientable, but not necessary homologous to zero; it
is given by $c\mapsto\PV_{c-[\RX]}$, where $[\RX]$ is the integral
fundamental class of~$\RX$ corresponding to some orientation. (In view
of~\iref{C4.Mikh.D}1 his map does not depend on the orientation, as
orientable components of~$\RX$ annihilate $\vr_2H_1(\RX)$.)

The other special case worth mentioning is when $\RX$ itself is
characteristic in~$X$ (and, hence, $[A]=0\bmod2H_2(X;\Z)$). In this case
one can take $c=0$, and the resulting quadratic form $\PV=\PV_0$ is
nothing but the Rokhlin-Guillou-Marin form on~$\RX$. There is an
alternative way to construct~$\PV$, which explains the name: it is easy
to see that, since $[\RX]=w_2(X)$, the Pontrjagin square
$P_2\:H_2(X)\to\Z_4$ descends to $\RH2_2=\iH_2$, and one can let
$\PV=P_2\circ\vr_2\:H_1(\RX)\to\Z_4$. In this form the construction admits
generalizations to nonhomogeneous classes in $H_*(\RX)$ and to varieties
of higher dimensions, not necessarily simply connected. Namely, let $X$ be
a closed complex variety, $\dim X=2n$, with real structure. (In fact, most
statements can be modified to hold for an arbitrary closed manifold with
orientation preserving involution whose dimension is divisible by~$4$;
details can be found in~\cite{Dg:PVF} or~\cite{DIK}.) Assume that the
Pontrjagin square $P_{2n}\:H_{2n}(X)\to\Z_4$ descends to $\iH_{2n}(X)$ and
define the \emph{Pontrjagin-Viro form} via
$\PV=P_{2n}\circ\bv_{2n}\:\CF^{2n}\to\Z_4$. When defined, $\PV$ is a
quadratic extension of Kalinin's intersection pairing
$*\:\CF^{2n}\otimes\CF^{2n}\to\ZZ$.

The following statements are proved in~\cite{Dg:PVF} (see also~\cite{DIK}).

\proposition\label{PVF.def}
The Pontrjagin-Viro form is well defined if and only if the Wu class
$u_{2n}(X'\sminus\RX)$ vanishes. In this case one has
$\Br\PV=\Gs(X)\bmod8$.
\qed
\endproposition

\proposition\label{PVF.sufficient}
The Pontrjagin square $P_{2n}$ descends to $\RH2_{2n}$ if and only if
$[\RX]$ realizes the Wu class $u_{2n}(X)$. Thus, $[\RX]=u_{2n}(X)$ is a
necessary condition for~$\PV$ to be well defined. The following are
sufficient conditions\rom:
\roster
\item
$X$ is an $M$-variety\rom;
\item
$[\RX]=u_{2n}(X)$ and $X$ is $\ZZ$-Galois maximal
\rom(see~\ref{1.1.8}\rom)\rom;\mnote{\Dg: edited after Zv}
\item
$[\RX]=u_{2n}(X)$ and $H_i(X)=0$ for $0<i<2n$.
\qed
\endroster
\endproposition

Assume now that $n=1$, i.e., $X$ is a surface. The condition
$u_2(X'\sminus\RX)=0$ in~\ref{PVF.def} is then equivalent to the
requirement that $u_2(X')=w_2(X')$ belong to the image of the inclusion
homomorphism $H_2(\RX)\to H_2(X')$, i.e., $\RX$ contains a subsurface
characteristic in $X'$. Such subsurfaces $F\subset\RX$ are characterized
by the property $\PV(x)=2([F]\circ x)\bmod4$ for all $x\in\CF^2\cap
H_0(\RX)$. Given such a surface~$F$, denote by $H_F$ the image of~$\CF^2$
under the natural map $\CF^2\hookrightarrow H_*(\RX)\to H_1(\RX)\to
H_1(F)$ and define a quadratic function $\PV_F\:H_F\to\Z_4$ via
$$
x_1\mapsto\PV(x_1+x_0)+2([F]\circ x_0),
$$
where $x_0\in H_0(\RX)$ is any element such that $x_1+x_0\in\CF^2$.

\proposition
The function $\PV_F\:H_F\to\Z_4$ coincides with the Rokhlin-Guillou-Marin
form of the characteristic surface $F\subset X'$. In particular, $H_F$
coincides with the kernel of the inclusion homomorphism $H_1(F)\to
H_1(X')$, the subspace $H_F\subset H_1(F)$ is informative \rom(in respect
to $\PV_F$\rom; see~\ref{forms} for the definition\rom), and one has $$
\chi(F)=\frac14[\chi(\RX)-\sigma(X)]+\Br\PV_F\bmod8.\rlap\qed $$
\endproposition

If, in addition, $H_1(X)=0$, the last statement implies~\iref{C4.Mikh.D}3
(for $A=\varnothing$); furthermore, in this case one can see that
the complex separation defined in~\ref{C4.Mikh} can be expressed in terms
of~$\PV$: two components $C_1$, $C_2$ of $\RX$ belong to the same half
$\RX_+^w$ or $\RX_-^w$ if and only if $\PV(\<C_1-C_2>)=0$.

\subsection{$\Spin$-congruences}\label{C2.2}
Below we reproduce two extremal congruences which take into account
$w_2(X)$ and orientability of~$\RX$. Other similar results, which employ
more subtle information about the homotopy type of~$X$, can be found
in~\ref{C2.5} (see also~\cite{DIK}, devoted to detailed study of real
Enriques surfaces; in fact, the congruences below are a by-product of the
results and tools developed in~\cite{DIK}).

\claim[Extremal congruence for non $\Spin$-surfaces]\label{C2.2.I}
Let $X$ be a complex closed surface with real structure and $w_2(X)\ne0$.
\roster
\item\local1
If $X$ is an $M$-surface, then $\RX$ is nonorientable.
\item\local2
If $w_2(X)$ belongs to the image of $\Tors H_2(X;\Z)$ in $H_2(X)$ and
$\RX$ is orientable, then $X$ is an $(M-d)$-surface, $d\ge2$, and
\subroster
\subdash
if $d=2$, then $\chi(\RX)=\Gs(X)\bmod{16}$\rom;
\subdash
if $d=3$, then $\chi(\RX)=\Gs(X)\pm2\bmod{16}$\rom;
\subdash
if\hbox{ $d=4$ }and $\chi(\RX)=\Gs(X)+8\bmod{16}$, then $\RX$ is of
type~$\Itors$.
\endsubroster
\endroster
\endclaim

\claim[Extremal congruence for $\Spin$-surfaces]\label{C2.2.H}
Let $X$ be a complex closed surface with real structure and $w_2(X)=0$.
\roster
\item\local1
If $H_1(X)=0$, then $\RX$ is orientable.
\item\local2
If $b_1(X)=0$ and $\RX$ is nonorientable, then $X$ is an
$(M-d)$-sur\-face, $d\ge2$, and
\subroster
\subdash
if $d=2$, then $\chi(\RX)=\Gs(X)\bmod{16}$\rom;
\subdash
if $d=3$, then $\chi(\RX)=\Gs(X)\pm2\bmod{16}$\rom;
\subdash
if\hbox{ $d=4$ }and $\chi(\RX)=\Gs(X)+8\bmod{16}$, then $\RX$ is of
type~$\Itors$.
\endsubroster
\endroster
\endclaim

To prove the congruence part (Statement~\loccit2) of~\ref{C2.2.H}
and~\ref{C2.2.I} we show that a part of $T_2(X)=\Tors
H_2(X;\Z)\otimes\ZZ\subset H_2(X)$ dies in~$\iH_2$: then $\dim J<d$ and
the congruences follow (cf.~\ref{C1.5}).

\proof[Proof of~\Ref{C2.2.I}]
The class $w_2(X)\ne0$ dies in $\iH_2$ if and only if $\RX$ is orientable
(see Table~\ref{tabKif} in~\ref{Kalinin.surface}), and the statement
follows.
\endproof

\proof[Proof of~\Ref{C2.2.H}]
\loccit1 If there is $x_1\in H_1(\RX)$ with $x_1^2=1$, then
$(\bv_2x_1)^2=1$. ($\bv_2x_1$ is well defined since $H_1(X)=0$.) This
contradicts to the assumption that the intersection form in~$H_2(X)$ and,
hence, in $\iH_2$ is even.

\leavevmode\loccit2 By the assumption, there is $x_1\in H_1(\RX;\Z)$ with
$x_1^2=1$. Similar to~\loccit1 one concludes that $\bv_1x_1\not\in\bv_1
H_0(\RX)$. Furthermore, any nonorientable component~$C_i$ of~$\RX$ is of
even genus (as otherwise $[C_i]^2=1\bmod2$). Then
$\bv_1x_1\not\in\Sq_1\iH_2$ and, hence, $H_2(X)/H_2(X;\Z)\otimes\Z_2$ is
not covered by~$\Im\bv_2$.
\endproof

\remark{Remark}
Alternatively, \iref{C2.2.H}1 can be derived
from~\ref{Spin.orient}: if $H_1(X)=0$, the only $\Spin$-structure on~$X$
is obviously $\conj$-invariant, and the construction of~\ref{Spin.orient}
produces a canonical semi-orientation on~$\RX$.
\endremark

\subsection{Surfaces in~$\Cp3$}\label{C2.3}
In the next few section we cite the known classification results for
surfaces. In order to describe the homeomorphism type of the real part of
a surface we denote by~$\SS p$, $p\ge0$, an orientable closed surface of
genus~$p$, $\SS p=\#_p(S^1\times S^1)$, and by~$\VV q$, $q\ge1$, a
nonorientable closed surface of genus~$q/2$, $\VV q=\#_q\Rp2$.
(See~\ref{C3.1} for our definition of genus in this case.)\mnote{\Dg:
edited avter Zv}
Let $S=S_0$.

From the topological point of view, surfaces in $\Cp3$ are traditionally
studied up to one of the following equivalence relations: homeomorphism
of the real part of a surface, ambient isotopy of the real parts in
$\Rp3$, rigid isotopy (i.e., isotopy in the class of nonsingular or, more
generally, equisingular in some appropriate sense surfaces of the same
degree), and rough projective equivalence (i.e., projective
transformation and rigid isotopy). The difference between the last two
relations is due to the fact that the group $\text{\sl PGL}(4;\R)$ of
projective transformations of $\Rp3$ has two connected components.
Certainly, the transformations in the component of unity transform a
surface into a rigidly isotopic one.

For any degree the number of classes is finite (for each of the four
equivalence relations); this follows from the fact that for each relation
the universal family of hypersurfaces of a given degree admits a finite
(semi-algebraic) stratification such that, in respect to the chosen
relation, the family is locally trivial over each strata. (For
topological equivalence relations the existence of a stratification was
proved by A.~Wallace \cite{Wallace} and A.~N.~Varchenko \cite{Varchenko}.)

Obviously, a real surface in~$\Cp3$ of degree~$1$ is a real projective
plane.  The classification of surfaces of degree~$2$ is a subject of
classical analytic geometry. Up to rigid isotopy there are eight types of
such surfaces (including singular ones), which are distinguished by the
topology of their real parts. A nonsingular real surface~$X\subset\Cp3$ of
degree~$2$ may be either hyperboloid ($\RX=S_1$), or ellipsoid ($\RX=S$),
or empty surface ($\RX=\varnothing$). Note that a nonsingular quadric~$X$
is isomorphic to $\Cp1\times\Cp1$, and when considered as an abstract
surface, $X$ admits four nonequivalent real structures: there are two
real structures with empty real part.

In accordance with the general setting of this survey we consider
nonsingular surfaces only. For such surfaces a complete classification
(for each of the four equivalence relations) is known in degrees~$\le4$.
In degrees $\ge5$ even the maximal number of connected components of the
real part of a surface is not known. (Note that the maximum can always be
realized by a nonsingular surface.) The best known upper bound is
obtained by adding the Smith and Comessatti (Petrovsky-Oleinik)
inequalities. It is worth mentioning that each inequality separately is
sharp (see \cite{ViroM}, where the refined Comessatti inequality
$1-\chi(\RX)\le \sigma^-(X)-2p_+$ is also shown to be sharp). However,
there is an evidence that the derived bounds for ~$\beta_0$ and~$\beta_1$
are not sharp.  Asymptotically one has $\max\beta_i\sim k_id^3$, $i=0,1$,
where $d$ is the degree, $\frac{13}{36}\le k_0\le\frac5{12}$, and
$\frac{13}{18}\le k_1\le\frac56$ (examples establishing the lower bounds
are constructed by F.~Bihan~\cite{Bihan}).

\subsubsection{Real cubic surfaces}
The classification of real surfaces in~$\Cp3$ of degree~3 is due to
L.~Schl\"afli~\cite{Schlafli}, F.~Klein~\cite{Klein}, and
H.~G.~Zeuthen~\cite{Zeuthen}.

\claim\label{C2.3.A}
There are five rigid isotopy types of nonsingular real cubic surfaces
in~$\Cp3$. They are distinguished by their real parts, which are $\VV7$,
$\VV5$, $\VV3$, $\VV1\+S$, or $\VV1$.
\endclaim

\proof
The topological invariants of the complex part of a nonsingular cubic
surface $X\subset\Cp3$ are given by~\ref{C2.1.E}: $\Gb_*(X)=9$,
$\Gs(X)=-4$, and $\Gs^-(X)=5$. Hence, the five homeomorphism types of the
real part listed above are all not prohibited by the Smith
inequality~\ref{C2.1.A}, extremal congruence~\ref{C1.5.A}, and
Comessatti-Petrovsky inequality~\ref{C2.1.B}. Thus, it remains to show
that all the five topological types are realized and that the space of
nonsingular cubics has at most five connected components. Consider the
space $\R\CC_3=\Rp{19}$ of all real cubic surfaces in~$\Cp3$ and let
$\R\Delta\subset\R\CC$ be the set of singular cubics,
$\Delta_0\subset\R\Delta$, the set of
cubics with a single nondegenerate double point, and
$\CS=\R\Delta\sminus\Delta_0$. Then $\dim\R\Delta=18$ and $\dim\CS\le17$,
and from the Poincar\'e duality and the exact sequence of the triple
$(\R\CC_3,\R\Delta,\CS)$ it follows that $\dim
H^0(\R\CC_3\sminus\R\Delta)=1+\dim\Ker\inc_*$, where
$\inc_*\:H_{18}(\R\Delta,\CS)\to H_{18}(\R\CC_3,\CS)$ is induced by the
inclusion. The group $H_{18}(\R\Delta,\CS)$ is generated by fundamental
classes of the connected components of $\Delta_0$, i.e., rigid isotopy
classes of cubics with a single nondegenerate double point. Such a cubic
(assuming that the singular point is $(1:0:0:0)$) is given by an equation
of the form
$$
x_0p_2(x_1,x_2,x_3)+p_3(x_1,x_2,x_3)=0,
$$
where $p_2$ and~$p_3$ are real homogeneous polynomials of degree~2 and~3
respectively. $(1:0:0:0)$ is the only singular point of the surface and it
is a nondegenerate double point if and only if the quadric $C_2=\{p_2=0\}$
is nonsingular and the cubic $C_3=\{p_3=0\}$ intersects~$C_2$
transversally. (Note that $C_3$ may be singular.) Through any set of six
disjoint points of~$C_2$ one can trace a cubic curve not containing~$C_2$,
and such curves form a 3-dimensional affine space. Hence, $\Delta_0$ has
at most five connected components: either $\R C_2=\varnothing$, or $\R
C_2\cong S^1$ and there are $2i$, $i=0,1,2,3$, real intersection points
of~$C_2$ and~$C_3$. Perturbing the singular surfaces obtained in this way
one can construct all the five topological types of nonsingular cubic
surfaces.

It remains to prove that $\dim\Ker\inc_*\le4=\dim H_{18}(\R\Delta,\CS)-1$,
and for this purpose, again due to the Poincar\'e duality, it suffices to
construct a topological circle in $\R\CC_3\sminus\CS$ which
intersects~$\Delta_0$ transversally at two points in distinct strata. Such
a circle is contained in the versal $2$-parameter perturbation of the
cuspidal surface given by
$$
x_0(x_1^2+x_2^2)+x_3^3+x_2^3=0.\rlap\qed
$$
\endproof\nofrills

\remark{Remark}
Topologically, the $V_q$ component of the real part of a nonsingular cubic
surface is situated in $\Rp3$ as the standard $\Rp2$ with unlinked and
unknotted handles attached. This can be deduced, e.g., from the
construction of the surfaces indicated in the proof above.
\endremark

\subsubsection{Real quartic surfaces}
The systematic study of the topology of the surfaces of degree~$4$ was
launched by K.~Rohn and D.~Hilbert at the beginning of this century.
(Hilbert
posed this question as part of his 16-th problem.)\mnote{toZv: modified}
However, first complete proofs of some of their statements appeared
later and are due to I.~Petrovsky and O.~Oleinik (the bound on the number
of spherical components) and to R.~Thom (the bound on the total Betti
number). The study was continued in the 60th--70th by G.~Utkin, and the
classification was completed by V.~Kharlamov~\cite{Kh:K3}, \cite{Kh:iso}
(the topological and isotopic classification; one isotopy class was
eliminated in collaboration with V.~Nikulin) and
V.~Nikulin~\cite{Nikulin:forms} (the rough projective classification). The
chirality of degree~$4$ surfaces was studied by
V.~Kharlamov~\cite{Kh:rigid}, \cite{Kh:chiral}.

As it is the case for any nonsingular real surface of even degree, the
real part of a nonsingular quartic is orientable. Its position in~$\Rp3$
is relatively simple: it is isotopic to a union of ellipsoids and
hyperboloids with unknotted and unlinked handles. With one exception the
components are outside each other; in the exceptional case the real part
consists of two nested spheres. (In most cases these results can be
derived from the Bezout theorem.)

With the only exception of two nested spheres the position of $\RX$ in
$\Rp3$ is determined up to isotopy by the images of the inclusion
homomorphisms $H_1(C_j)\to H_1(\Rp3)=\ZZ$, where $C_j\subset\RX$ are the
components of $\RX$. (This follows from general theorems of
$3$-dimensional topology and the fact that the handles are unknotted and
unlinked.) For a nonspherical component we use the notation $S_p^0$
(respectively, $S_p^1$) if the image is (respectively, is not) trivial.
Further, $nS$ stands for $n$ unnested spheres and $S(S)$, for a nest of
two spheres.

In the case of degree~$4$ surfaces all the four classifications
(topological, isotopic, rough projective, and rigid) are distinct.

\claim
[Classifications of real quartic surfaces]\label{qua.cla} A\mnote{toZv:
new title} real nonsingular surface of degree~$4$ is determined up to
rigid isotopy by the isotopy type of its real part~$\RX$, its complex type
\rom($\I0$, $\Ihp$, or $\II$, see~\ref{C1.4}\rom) and its chirality, i.e.,
whether the surface is rigidly isotopic to its mirror image. The isotopy
types of the real parts are those obtained by a sequence of Morse
simplifications \rom(see \ref{C2.5}\rom)\mnote{\Dg: extended after Zv}
from the following extremal types\rom:

\topinsert
\input f-4s-nc
\endinsert

\roster
\item\local1
$M$-surfaces\rom:\quad$\SS{10}^1\dplus S$, $\SS6^1\dplus5 S$,
$\SS2^1\dplus9 S$\rom;
\item \local2
$(M-2)$-surfaces\rom:\quad$\SS7^1\dplus2 S$, $\SS3^1\dplus6 S$,
${S}_a^0\dplus(9-a)S$, $0\le a\le 9$\rom;
\item\local3
pairs of tori\/\rom:\quad $\SS1^1\dplus\SS1^1$, $\SS1^0\dplus\SS1^0$\rom;
\item\local4
pair of nested spheres\/\rom:\quad $S(S)$.
\endroster
The following is the complete list of nonchiral isotopy types\rom:
$pS\dplus{S}_q^0$, $p\ge3$, $q\ge3$, and $pS\dplus{S}_q^1$, $p\ge4$. The
complex types are shown in Tables~\ref{4surf.nc} and~\ref{4surf.c}.
\endclaim

\remark{Remark}
As it follows from~\ref{qua.cla}, with the exception of the case
$\RX=S(S)$ the isotopy type of the real part~$\RX$ of a real quartic
surface in $\Cp3$ is determined by its topological type and whether $\RX$
is contractible in $\Rp3$ or not.
\endremark

\subsubsection{Surfaces of degree~$5$}
Starting from degree~5, one observes the following phenomenon: on one
hand, any surface obtained by a small deformation of a surface of a given
degree in $\Cp3$ can be embedded in $\Cp3$ (as a surface of the same
degree); on the other hand, by a big deformation one can obtain a surface
which can not be embedded in $\Cp3$. (It is worth mentioning that in
degree 4 the situation is just the opposite: a nonsingular quartic
surface can be made nonprojective by a small deformation, and, on the
other hand, any surface of the same deformation class, i.e., a
$K3$-surface, see~\ref{C2.5.B} below, can be made a nonsingular quartic
by a small deformation. This explains why the topological classification
of real nonsingular degree 4 surfaces is the same as the topological
classification of real $K3$-surfaces.)

In degree~5 the phenomenon was studied by E.~Horikawa \cite{Horikawa}. He
proved that the base of a versal deformation of a nonsingular quintic
consists of two smooth irreducible components $M_0$, $M_1$ intersecting
normally, $\dim_{\C}M_0=\dim_{\C}M_1=40$, and $\dim_{\C}M_0\cap M_1=39$.
The points of $M_0\sminus(M_0\cap M_1)$ correspond to nonsingular
quintics, those of $M_1\sminus(M_0\cap M_1)$, to double coverings of
$\Cp1\times\Cp1$, and those of $M_0\cap M_1$, to double coverings of the
quadratic cone~$\Sigma_2$. Certainly, the deformation of a real quintic
can be made equivariant and its members can be used to produce quintics
of the same real deformation class. This observation was used to prove
the existence part of~\ref{quintimax} and~\ref{Mquinti}.

\claim[Proposition]\label{quintimax}
The extremal values of the Betti numbers of nonsingular quintics
in~$\Cp3$ are in the range $22\le\max\beta_0\le25$, $\max\beta_1=45$
or~$47$. For the surfaces in the same deformation class one has
$23\le\max\beta_0\le25$ and $\max\beta_1=47$.
\endclaim

\proof
For a nonsingular quintic~$X$ one has $h^{1,1}(X)=45$, $\beta_{*}(X)=55$,
and $\sigma(X)=-35$. Thus, the upper bounds follow from the Smith and
Comessatti inequalities \ref{C2.1.A}, \ref{C2.1.B} and the extremal
congruence~\ref{C1.5.A}. Quintics with $\beta_1=45$ were constructed by
V.~Kharlamov~\cite{KhLen}, quintics with $\beta_0=22$, by I.~Itenberg and
V.~Kharlamov~\cite{ItenKh}, and $\C$-deformation equivalent surfaces with
$\beta_0=23$ and\mnote{toZv:extended} those with $\beta_1=47$, by
F.~Bihan~\cite{Bihan}.
\endproof

\claim[Proposition]\label{Mquinti}
The real part of an $M$-surface $\C$-deformation equivalent to a
nonsingular quintic can only have $5$, $9$, $13$, $17$, $21$, or $25$
connected components. $M$-quintics with $5$, $9$, $13$, $17$, and $21$
connected components of the real part do exist.
\endclaim

\proof
The other values are prohibited by the extremal congruence \ref{C1.5.A}
and the Smith inequality~\ref{C2.1.A}. $M$-surfaces with 5, 9, 13, 17,
and 21 components were constructed by V.~Kharlamov~\cite{KhLen} (the
construction for 17 and 21 components is based on some auxiliary curves
constructed by O.~Viro).
\endproof

The $M$-surfaces constructed in~\cite{KhLen} are
$$
20S\+V_{13},\quad 16S\+V_{21},\quad 12S\+V_{29},\quad 8S\+V_{37},\quad
4S\+V_{45}.
$$
Their \emph{twins}, obtained by switching the sign of the double
covering, are
$$
4S\+V_{43},\quad 8S\+V_{35},\quad 12S\+V_{27},\quad
11S\+S_1\+V_{25},\quad 16S\+V_{19}.
$$
Among other pairs constructed at the moment are $19S\+S_2\+V_{7}$,
$4S\+V_{43}$ and $19S\+S_1\+V_{7}$, $4S\+V_{45}$. As it follows from the
construction of the auxiliary curves, the Morse simplifications of all the
topological types above can be realized by real quintics. The variety of
the examples obtained in this way has never been analyzed systematically.

The surface with $22$ components constructed by V.~Kharlamov and
I.~Itenberg in~\cite{ItenKh} is $21S\+V_7$. Its twin is
$S\+\SS2\+V_{41}$, with $\beta _{1}=45$.

Recently B.~Chevallier (unpublished) found some additional examples. His
approach looks different (it is based on Viro's patching method), but
there is an evidence that both the approaches are close or even
equivalent: the approach based on the Horikawa theorem also relies on
Viro's method to construct auxiliary curves, while Chevallier's approach
also uses double quadrics. The new
types
with the total Betti number or Euler characteristic close to extremal
ones are\mnote{toZv: the sentence and the list changed;\Dg: one type
removed after Zv}
$$
\gather
18S\+S_1\+V_9,\quad 17S\+2S_1\+V_3,\quad
17S\+S_2\+V_5,\,\dots,\,
4S\+S_8\+V_{23},\\
3S\+S_1\+V_{39},\quad 3S\+S_7\+V_{27},\quad
2S\+S_1\+V_{41},\quad S\+S_{18}\+V_3,
\quad S\+S_{10}\+V_{21}.
\endgather
$$

It is worth mentioning that if an $M$-quintic with 25 components does
exist, then, according to the known inequalities and congruences, it must
be of one of the following topological types:
$$
24S\+V_5,\quad 23S\+\SS1\+V_3,\quad 23S\+\SS2\+V_1,\quad
22S\+\SS1\+\SS1\+V_1.
$$
The two last types are prohibited (the forth one, O.~Viro, and the third
one, by V.~Kharlamov and I.~Itenberg, cf\.\cite{KhOber}). The problem of
existence of $M$-surfaces of degree 5 of the two other topological types
is still open.

\subsection{Rational Surfaces}\label{C2.4}
In accordance with the general setting of this survey, by a \emph{real
rational surface} we mean a rational (over~$\C$) surface supplied with a
real structure. It is worth emphasizing that such a surface does {\bf
not} need to be rational over~$\R$; Theorem~\ref{C2.4.A} below gives
plenty of examples. Real rational surfaces were classified by
Comessatti~\cite{Comessatti:1}--\cite{Comessatti:3}. Following
Comessatti's approach, let us confine ourselves to \emph{minimal}
(over~$\R$) real rational surfaces. (Recall that a nonsingular
surface~$X$ is called (relatively) minimal if any holomorphic degree~$1$
map $X\to X'$ is a biholomorphism. A real surface $(X,\conj)$ is minimal
over~$\R$ if the above condition holds for $\conj$-equivariant maps. A
surface is minimal over~$\C$ if and only if it does not contain a
$(-1)$-curve, and it is minimal over~$\R$ if and only if it does not
contain a real $(-1)$-curve or a pair of disjoint conjugate
$(-1)$-curves.)

\claim[Minimal real rational surfaces]\label{C2.4.A}
The following is a complete list of rigid isotopy types of minimal
over~$\R$ rational surfaces with anti-holomorphic involution\rom:
\roster
\item\local1
real projective plane $\Cp2$\rom:\quad$\RX=V_1$\rom;
\item\local2
real quadric $\Cp1\times\Cp1$\rom: there are four types\rom: $\RX=S_0$,
$\RX=S_1$, and two nonequivalent
surfaces with $\RX=\varnothing$\rom;
\item\local3
ruled rational surfaces $\Sigma_m$, $m\ge2$\rom:
\smallskip
\halign{\hbox to.5\hsize{\qquad\qquad#\hfil}&\qquad#\hfil\cr
$m$ even\rom:\quad$\RX=\varnothing$ or $S_1$,& $m$
odd\rom:\quad$\RX=V_2$\rom;\cr}
\smallskip
\item\local4
real conic bundles over $\Cp1$ whose reducible fibers are all real and
consist of pairs of complex conjugate exceptional curves\rom:\quad
$\RX=mS$, where $2m\ge4$ is the number of reducible fibers\rom;
\item\local5
Del Pezzo surfaces of degree $d=K^2=1$ or~$2$\rom:
\smallskip
\halign{\hbox to.5\hsize{\qquad\qquad#\hfil}&\qquad#\hfil\cr
$d=1$\rom:\quad$\RX=V_1\dplus4S$,& $d=2$\rom:\quad$\RX=3S$ or $4S$.\cr}
\endroster
\endclaim

(From the complex point of view a Del Pezzo surface of degree~$d\le9$ is
obtained from~$\Cp2$ by blowing up $(9-d)$ points in general position,
i.e., such that no three of them lie in a line and no six of them lie in a
quadric.) Note that the Del Pezzo surface of degree~2 with $\RX=3S$ can
also be represented as a conic bundle over~$\Cp1$ with six reducible
fibers. Note also that in the above list only surfaces~\loccit1--\loccit3
are minimal over~$\C$.

Theorem~\Ref{C2.4.A} depends on a classification of anti-birational
transformations of~$\Cp2$, which, in turn, is based on considering the
fundamental points of such a transformation. A representation of
Comessatti's proof in the modern language can be found
in~\cite{CilibertoPedrini}.

\subsection{Abelian, $K3$-, and Enriques Surfaces}\label{C2.5}
The next group in the classification of complex algebraic surfaces are
surfaces of Kodaira dimension~$0$, i.e., those whose minimal models are
either abelian, or hyperelliptic, or~$K3$-, or Enriques
surfaces, each type, except hyperelliptic surfaces, forming one
deformation family. Since the minimal model of a complex surface of
Kodaira dimension $\ge0$ is unique (see, e.g.,
Hartshorne~\cite{Hartshorne}), it must be real if the original surface is
real. Hence, it suffices to consider surfaces minimal over~$\C$; any other
surface can be obtained from a minimal (over~$\C$) one by a sequence of
blow-ups at real points or pairs of complex conjugate points.

The classification of real hyperelliptic surfaces is still in progress.
Below we treat the other three types.\mnote{\Dg: par added}

An abelian surface~$X$ (which topologically is a $4$-torus $S^1\times
S^1\times S^1\times S^1$) can be characterized by $b_1(X)=4$ and $K_X=0$.
Anti-holomorphic involutions on such surfaces were classified by
Comessatti~\cite{Comessatti:3}, whose result is the following:

\claim[Real Abelian Surfaces]\label{C2.5.A}
Up to equivariant deformation an anti-holo\-mor\-phic involution $\conj$
on a complex torus~$X$ is determined by the real point set~$\RX$ and, if
$\RX=\varnothing$, by the action of $\conj_*$ on $H_1(X;\Z)$. There are
five classes\rom: $\RX=k\SS1$ with $k=1,2,4$ and two classes with
$\RX=\varnothing$.
\endclaim

A description of the moduli space of algebraic abelian surfaces with
non\-empty real part can be found in Silhol~\cite{Silhol}.

\proof[Sketch of the proof]
Consider the uniformization $X=\C^2/L$, where $L\subset\C^2$ is a lattice
of rank~4, and lift $\conj$ to an anti-automorphism~$c$ of~$\C^2$. If
$\RX\ne\varnothing$, the lift can be chosen to be an involution
preserving~$L$. Since $c$ has eigenvalues~$\pm1$ of multiplicity~2 each,
there are three isomorphism classes of the restriction $c\big|_L$, and
for each class there is a unique up to homotopy equivariant embedding
to~$\C^2$ (considered with the standard conjugation).

If $\RX=\varnothing$, the above arguments apply to the adjoint action
$\Ad\conj$ on~$X$ considered as group, and $\conj$ itself is the
composition of $\Ad\conj$ and the affine shift by an element
$a\in\Ker(1+\Ad\conj)$, which is to be considered up to
$\Im(1-\Ad\conj)$. Up to isomorphism, there are two solutions, one with
$\Ad\conj=\id$, and one with $\Ad\conj=\id\oplus\sigma$, where $\sigma$
is a real structure on an elliptic curve $C'=\C/L'$ such that $\R C=S^1$.
\endproof

A complex $K3$-surface~$X$ is characterized by $b_1(X)=0$, $K_X=0$.
(Topologically $X$ is a simply connected $4$-manifold with
$H_2(X;\Z)=2(-E_8)\oplus3U$). In order to state the classification result
for real $K3$-surfaces, let us introduce the following notion: a
\emph{topological type} (i.e., a class of surfaces with homeomorphic real
parts) is called \emph{extremal} (within a fixed deformation family of
complex surfaces) if it cannot be obtained from another topological type
by a \emph{Morse simplification}, i.e., Morse modofication decreasing the
total Betti number. (A Morse simplification is either removing a
spherical component or contracting a handle.) The following result is due
to V.~Kharlamov~\cite{Kh:K3} and N.~Nikulin~\cite{Nikulin:forms}:

\claim[Real $K3$-surfaces]\label{C2.5.B}
A $K3$-surface~$X$ with anti-holomorphic involution is determined up to
equivariant deformation by the topological type of its real part~$\RX$
and its type \rom(i.e., whether $[\RX]$ is or is not homologous to zero
in $H_2(X)$, see~\ref{C1.4}\rom).  The topological types of the real
parts are those and only those which can be obtained by a sequence of
Morse simplifications from the following extremal types\rom:
\roster
\item\local1
$M$-surfaces\rom:\quad$\SS{10}\dplus S$, $\SS6\dplus5S$,
$\SS2\dplus9S$\rom;
\item \local2
$(M-2)$-surfaces\rom:\quad$\SS7\dplus2 S$, $\SS3\dplus6S$\rom;
\item\local3
pair of tori\rom:\quad $2\SS1$.
\endroster
The topological types of the real part of a real $K3$-surface are listed
in Table~\ref{K3.nc}.

\topinsert

\def\mbox#1{\hbox to6pt{\hss$#1$\hss}}
\def\.{\mbox{\bullet}}
\let\,\.
\def\*{\mbox{*}}
\def\d{\mbox{\cdot}\kern-6pt\mbox{\circ}}
\hbox to\hsize{\hss\vbox{\offinterlineskip
\ialign{\hss\eightrm#\quad\vrule&&\strut\hbox to15pt{\eightrm\hss#\hss}\cr
$a$&\cr
&\cr
0&&\*&&&&&&&&\*&&&&&&&&\*\cr
1&\.&&\,&&&&&&\,&&\,&&&&&&\.&&\,\cr
2&&\.\*&&\.&&\*&&\,&&\,\*&&\,&&\*&&\.&&\.\*&&\.\cr
3&&&\.&&\,&&\,&&\,&&\,&&\,&&\.&&\.&&\.\cr
4&&&&\.&&\.\*&&\,&&\,\*&&\,&&\.\*&&\.&&\.\*\cr
5&&&&&\.&&\,&&\.&&\,&&\.&&\.&&\.\cr
6&&&&&&\.&&\.&&\.\*&&\.&&\.\*&&\.\cr
7&&&&&&&\.&&\.&&\.&&\.&&\.\cr
8&&&&&&&&\.&&\.\d&&\.&&\.\cr
9&&&&&&&&&\.&&\.&&\.\cr 
10&&&&&&&&&&\.\d&&\.\cr
11&&&&&&&&&&&\.&&&&&&\cr
\noalign{\hrule}
&1&2&3&4&5&6&7&8&9&10&11&12&13&14&15&16&17&18&19&20&&r\cr}}\hss}
\figure[Table]\label{K3.nc}
Topological types of real $K3$-surfaces
\endfigure

\bigskip\eightpoint\noindent
\* and \d\ stand for a surface of type~$\I{}$ and \. stands for a surface
of type~$\II$. The surfaces of type~$\I{}$ denoted by \d\ have $\RX=2S_1$
($(a,r)=(8,10)$) or $\RX=\varnothing$ ($(a,r)=(10,10)$); the surfaces
with $a+r=22$ (the right edge) have $\RX=(r-10)S$; the other surfaces
have $\RX=S_g\+kS$, where $g=\12(22-a-r)$ and $k=\12(r-a)$.

\endinsert
\endclaim

Proof of Theorem~\Ref{C2.5.B} is based on the global Torelli theorem for
$K3$-surfaces, see~\cite{BPV}, which describes their moduli space. This
reduces the problem of classification of anti-holomorphic involutions on
$K3$-surfaces up to rigid isotopy to the arithmetical problem of
classification of involutive isometries on $2(-E_8)\oplus3U$ up to
automorphism of the lattice. We omit the details; see~\cite{Nikulin:forms}
or~\cite{DIK}.

The remaining deformation type, Enriques surfaces, can be characterized by
$2K_X=0$, $K_X\ne0$. (Another definition is the following: $X$ is an
Enriques surface if $\pi_1(X)=\Z_2$ and the universal covering of~$X$ is a
$K3$-surface.) Like the previous cases, a real Enriques surface~$X$ is
determined up to deformation by the topology of the complex conjugation
involution $\conj\:X\to X$. However, the topology of the real part~$\RX$
and its type are not sufficient anymore to determine the topology of the
involution; one needs more sophisticated invariants. Below we cite the
principal results on the deformation classification of real Enriques
surfaces; details and a complete list of the deformation classes are
found in~\cite{DIK}.

\claim[Topology of the real part]\label{C2.5.C}
There are\/ \rom{87} topological types of Enriques surfaces with
anti-holomorphic involution. Each of them can be obtained by a sequence
of Morse simplifications from one of the\/ \rom{22} extremal types listed
below. Conversely, with the exception of the two types $6\SS0$ and
$\SS1\dplus5\SS0$, any topological type obtained in this way is realized
by a real Enriques surface.

The\/ \rom{22} extremal types are\/\rom: \bgroup
\edef\1{\kern\the\parindent}
\def\fourcol{\par\bgroup\interlinepenalty200
 \openup1pt\halign\bgroup\1\1##&&\hbox to.23\hsize{\1\1$##$\hss}\cr}
\def\endfourcol{\crcr\egroup\egroup\par\smallbreak}
\def\lab(#1)#2{\llap{$(#1)$\space}\hbox to0pt{$#2$\hss}}
\roster
\item
$M$-surfaces\/\rom: \fourcol
&\lab(a){\chi(\RX)=8\:}     &&\lab(b){\chi(\RX)=-8\:}\cr
&4\VV1\dplus2\SS0,      &2\VV2\dplus4\SS0,
&\VV{11}\dplus\VV1,         &\VV8\dplus\VV4,\cr
&\VV2\dplus2\VV1\dplus3\SS0,&\VV4\dplus5\SS0,
&\VV{10}\dplus\VV2,         &\VV7\dplus\VV5,\cr
&\VV3\dplus\VV1\dplus4\SS0, &\VV2\dplus\SS1\dplus4\SS0,
&\VV9\dplus\VV3,        &2\VV6,\cr
&                           &&&\VV{10}\dplus\SS1;\cr
\endfourcol
\item
$(M-2)$-surfaces with $\chi(\RX)=0$\rom: \fourcol
&\VV4\dplus2\VV1,           &\VV5\dplus\VV1\dplus\SS0,
&2\VV3\dplus\SS0,       &\VV6\dplus2\SS0,\cr
&\VV3\dplus\VV2\dplus\VV1,  &\VV4\dplus\VV2\dplus\SS0,
&\VV4\dplus\SS1\dplus\SS0,  &2\VV2\dplus\SS1;\cr
\endfourcol
\item
pair of tori\/\rom:\quad $2\SS1$.
\qed
\endroster
\egroup
\endclaim

Let $X$ be a real Enriques surface and $\tX\to X$ the covering
$K3$-surface. The real structure $\conj$ on~$X$ lifts to two commuting
real structures $t\ix1,t\ix2\tX\to\tX$, whose real parts $\R\tX\ix1$,
$\R\tX\ix2$ are disjoint. The images $\RX\ix{i}$ of $\R\tX\ix{i}$,
$i=1,2$, are also disjoint. Thus, the real part $\RX$ of~$X$ splits into
two subsets $\RX\ix1$, $\RX\ix2$ (which are unions of whole components),
called \emph{halves}. The \emph{half decomposition}
$\RX=\RX\ix1\cup\RX\ix2$ is an important topological (and deformation)
invariant of the surface; it is designated by an expression of the form
$$
\RX=\{\text{components of $\RX\ix1$}\}\+\{\text{components of $\RX\ix2$}\}.
$$
Clearly, one can speak about the type ($\I0$, $\Iu$, or~$\II$) of each
half $\RX\ix{i}$ separately, as well as about the types ($\I{}$ or $\II$)
of $\R\tX\ix{i}$. Furthermore, $\RX\ix1$ is naturally identified with the
fixed point set of the descent of~$t\ix1$ to $\tX/t\ix2$; thus, one can
also consider the type of $\RX\ix1$ in $\tX/t\ix2$. The latter quotient
is diffeomorphic to a rational surface, if $\RX\ix2\ne\varnothing$, or to
an Enriques surface otherwise. In the former case the possible types are
$\I{}$ and~$\II$, in the latter case, $\I0$, $\Iu$, and~$\II$.

\proposition\label{K3toE}
The covering $\R\tX\ix1\cup\R\tX\ix2\to\RX$ is the orientation double
covering.
\qed
\endproposition

\theorem[Half decomposition]\label{C2.5.D}
Each half of a real Enriques surface is either
\roster
\item\local1
$\Ga V_g\+a V_1\+b S$ with $g>1$, $a\ge0$, $b\ge0$, $\Ga=0,1$, or
\nopagebreak
\item\local2
$2V_2$, or\nopagebreak
\item\local3
$S_1$.
\endroster
With the exception of the cases $\RX=kS$ and $\RX=V_{2r}\+kS$ a half
decomposition of a real part~$\RX$ as in~\ref{C2.5.C} is realizable if and
only if it satisfies~\loccit1--\loccit3 above. The realizable half
decompositions of the exceptional real parts are listed in
Figure~\ref{halves.table}.
\qed

\topinsert
\begingroup
\dimendef\step=3  \step=10pt \dimendef\thin=4  \thin=.2pt
\dimendef\thick=5 \thick=.4pt \dimendef\gap=6   \gap=1pt
\countdef\loopcount=2

\def\grid#1#2{\vbox{\offinterlineskip \parindent 0pt \hsize#1\step
 \advance\hsize by\thick \advance\hsize by\gap
 \dimen0 \hsize \kern\gap \advance\hsize-\thin
 \def\vert##1{\vbox to0pt{\hsize##1 %
  \vss\vrule height\dimen0 width##1}\kern-\thin}%
 \def\zbox##1{\vbox to0pt{\hsize0pt \vss\hbox{\raise##1 \box0 }}}
 \def\labelbox##1{\setbox0 \hbox to0pt{$\scriptstyle##1$}%
  \dimen1=-\ht0 \ht0=0pt \dp0=0pt}
 \def\xlabel##1{\line{\kern\hsize\labelbox{\,\,##1\hss}\zbox{1pt}\hss}}
 \def\ylabel##1{\line{\labelbox{\hss##1\,}\advance\dimen1 by\dimen0 %
  \zbox{\dimen1}\hss}}
 \def\pointbox##1(##2,##3){\dimen1\step \multiply\dimen1 by##2 %
  \dimen2\step \multiply\dimen2 by##3
  \setbox0\vbox to0pt{\hsize 0pt\vss\kern-2\thick\kern2\thin\line
   {\hss\kern2\thick\kern2\thin$##1$\hss}\vss}%
  \moveright\dimen1\hbox{\zbox{\dimen2}}}
 \def\black{\pointbox\bullet}
 \def\white{\pointbox\circ}
 \loopcount#1 \loop\hrule height\thin width\hsize \kern\step
  \kern-\thin \advance\loopcount-1 \ifnum\loopcount>0 \repeat
 \hrule height\thick \line{\vert\thick
 \loopcount#1 \loop \kern\step \vert\thin
  \advance\loopcount-1 \ifnum\loopcount>0 \repeat \hfil}
 \line{\vbox{#2}}}}

\def\graph#1#2{\grid#1{\xlabel{a}\ylabel{b}
 \def\(##1,##2){\black(##1,##2)}#2}}
\def\capt#1{\hidewidth$\scriptstyle#1$\hidewidth}
\def\line#1{\hbox to\hsize{#1}}
\def\qqquad{\kern4ex}

\line{\hfil\openup0pt\vbox{\halign
{#\qqquad&\qqquad#\qqquad&\qqquad#\qqquad&\qqquad#\cr
\graph5{\(0,0)\(0,1)\(0,2)\(0,3)\(0,4)\(1,0)\(1,1)\(1,2)\(1,3)\(1,4)
 \(2,0)\(2,1)\(3,0)\(3,1)\(4,0)\(4,1)\(2,2)}&
\graph5{\(0,0)\(0,1)\(0,2)\(0,3)\(0,4)\(1,0)\(1,1)\(1,2)\(1,3)\(1,4)}&
\hbox to5\step{\hfil\graph2{\(0,0)\(0,1)\(0,2)}\hfil}& \hbox
to5\step{\hfil\graph1{\(0,0)\(0,1)}\hfil}\cr \noalign{\vglue-2pt}\cr
\capt{\{a\SS0\}\dplus\{b\SS0\},}& \capt{\{\VV4\dplus
a\SS0\}\dplus\{b\SS0\}}& \capt{\{\VV6\dplus a\SS0\}\dplus\{b\SS0\}}&
\capt{\{\VV8\dplus a\SS0\}\dplus\{b\SS0\},}\cr \capt{\{\VV2\dplus
a\SS0\}\dplus\{b\SS0\}}&&& \capt{\{\VV{10}\dplus
a\SS0\}\dplus\{b\SS0\}}\cr}}\hfil}
\endgroup
\figure\label{halves.table}
Half decompositions of exceptional real Enriques
surfaces
\endfigure

\endinsert
\endtheorem

Proposition~\ref{K3toE} follows essentially from the fact that $\tX$ is
$\Spin$ and $X$~is not; thus, the Enriques involution reverses the
$\Spin$-orientation of~$\R\tX$ (see~\ref{Spin.orient}). Together with the
restrictions to the topology of~$\R\tX$ (see~\ref{C2.5.B}) this gives the
possible forms~\loccit1--\loccit3 for a half of~$\RX$. The exceptional
half decompositions are treated using certain congruences similar to the
extremal congruences~\ref{C1.5.A}.

In order to prohibit almost all nonexisting topological types it suffices
to combine~\ref{K3toE} (and~\ref{C2.5.B}) with the Smith
inequality~\ref{C2.1.A}, generalized Comessatti-Petrovsky
inequality~\ref{C2.1.B}, and extremal congruences~\ref{C1.5.A}. After
that, there are only two types, $\SS1\dplus5\SS0$ and $3\VV2$, and one
series, $\SS1\dplus\VV1\dplus\dots$, still left to be prohibited. This is
done by considering Kalinin's spectral sequence (see~\cite{DIK} for
details).

For a technical reason real Enriques surfaces are subdivided into three
classes, those of hyperbolic, parabolic, or elliptic type, depending on
whether the minimal Euler characteristic of the components of~$\RX$ is
negative, zero, or positive.

\theorem[Deformation classes]\label{S.def-cl}
With few exceptions listed below the deformation class of a real Enriques
surface~$X$ with a distinguished half~$\RX\ix1$ is determined by the
topology of its half decomposition. The exceptions are\rom:
\roster
\item\local1
$M$-surfaces of parabolic and elliptic type, i.e., those with
$\RX=2V_2\+4S$, $V_2\+2V_1\+3S$, or $4V_1\+2S$\rom; the invariant is the
Pontrjagin-Viro form\rom;
\item\local2
surfaces with $\RX=2V_1\+4S$\rom;
\item\local3
other surfaces with a half $\RX\ix1=4S$\rom; the additional invariants are
the types, $\Iu$, $\I0$, or~$\II$, of~$\RX\ix1$ in~$X$ and
$\tX/t\ix2$\rom;
\item\local4
surfaces with $\RX=\{V_{10}\}\+\{\varnothing\}$,
$\{V_4\+S\}\+\{\varnothing\}$, $\{V_2\+4S\}\+\{\varnothing\}$, and
$\{2S\}\+\{2S\}$\rom; the additional invariant is the type, $\Iu$
or~$\I0$, of~$\RX$ in~$X$\rom;
\item\local5
surfaces with $\RX=2V_1\+3S$\rom; the additional invariant is the type,
$\Iu$ or~$\II$, of~$\RX$ in~$X$\rom;
\item\local6
surfaces with $\RX=\{S_1\}\+\{S_1\}$.
\qed
\endroster
\endtheorem

The additional invariant needed to distinguish the deformation classes in
cases \iref{S.def-cl}2 and~\ditto6 are more sophisticated and depend on
the special arithmetical properties of the induced
$(\ZZ\times\ZZ)$-action on $H_2(\tX;\Z)$ (in case~\ditto2) or the
topology of mutual position of the two halves $\R\tX\ix{i}=S_1\+S_1$
in~$\tX$ (in case~\ditto6).

The deformation classification of real Enriques surfaces is based upon
two different approaches: Donaldson's trick and arithmetical calculations
using the global Torelli theorem for $K3$-surfaces. (Although each of the
two approaches could do the job, the complexity of the resulting
auxiliary problem depends essentially on the surface in question. The
former, more geometric approach, works well for surfaces of hyperbolic
and parabolic types; the latter, for those of elliptic type.) Donaldson's
trick employs the hyper-K\"ahler structure on~$\tX$: one can change the
complex structure so that the Enriques involution becomes
anti-holomorphic and $t\ix2$ becomes holomorphic; then the quotient
$\tX/t\ix2$ is a real rational surface (unless $\RX\ix2=\varnothing$) and
the problem reduces to the study of real rational surfaces with a
nonsingular anti-bicanonical curve (the image of $\RX\ix2$). In the
second approach one uses the global Torelli theorem to establish a
one-to-one correspondence between the deformation classes of real
Enriques surfaces and the isomorphism classes of $(\ZZ\times\ZZ)$-actions
(induced by~$t\ix1$ and~$t\ix2$) on the intersection lattice
$H_2(\tX;\Z)=2(-E_8)\oplus3U$; the action must satisfy certain arithmetic
conditions to ensure the realizability by anti-holomorphic involutions and
the freeness of the composition $t\ix1\circ t\ix2$. The resulting
arithmetical problem is solved using advanced theory of discriminant
forms, and the relation between the properties of the action and the
topology of $(X,\conj)$ is established using mainly Kalinin's spectral
sequence.

\end



\ifx\rokhloaded\undefined\input rokh.def \fi

\input epsf

{\catcode`\@=11
\gdef\proclaimfont@{\sl}
\global\let\proclaimfont\proclaimfont@}

\def\smalldelim#1#2{#1{\vcenter{\hbox{$\scriptscriptstyle#2$}}}}
\def\smallleft{\smalldelim\mathopen(}
\def\smallright{\smalldelim\mathclose)}

\let\ix\index

\def\scap{\mathbin{\scriptstyle\cap}}
\def\plus{\mathbin{\scriptstyle\sqcup}}
\let\+\plus

\document

\section{Curves on Surfaces}\label{C3}

\subsection{Congruences}\label{C3.1}
A topological curve (i.e., disjoint union of embedded circles) $L$ in a
topological surface~$S$ is said to be \emph{totally homologous to zero}
if the inclusion homomorphism $H_1(L)\to H_1(S)$ is trivial. Further,
define the \emph{genus} of a topological surface~$S$ to be
$$
g(S)=b_0(S)-\frac12\bigl(\chi(S)+b_0(\partial S)\bigr).
$$
(As there seems to be no commonly accepted definition of genus of a
nonorientable or disconnected surface, we choose the one that suits our
needs. Certainly, for a connected orientable surface our definition gives
the usual genus. In general, the genus is the total number of handles, a
cross-cup
counted as half a handle.)
Note that a surface of genus~$0$ is a union of spheres with holes, and a
surface of genus~$1/2$ is a union of spheres with holes and one real
projective plane with holes. It is also easy to check that if a closed
surface~$S$ is represented as the union $S_1\cup S_2$ of two surfaces so
that their common boundary $\partial S_1=\partial S_2$ is totally
homologous to zero in~$S$, then $\dim\Im[\inc_*\: H_1(S_i)\to
H_1(S)]=2g(S_i)$, $i=1,2$.

\theorem\label{C3.1.A}
Let $X$ be a complex surface with real structure, $A$ an ample even real
divisor on~$X$ with $\RA\ne\varnothing$, and $Y$ the double covering of
$X$ branched over $A$. Assume that $H_1(X)=0$.  The following are
necessary conditions for $\RY_-$ to be an $M$- or $(M-1)$-surface\rom:
$X$ is an $M$-surface, $\RA$~lies entirely in a single component
of~$\RX$, all the other components of~$\RX$ belong to~$\RX_-$ and are
covered trivially by~$\RY_-$, and~$\RA$ is totally homologous to zero
in~$\RX$. If these conditions are satisfied, then
\roster
\item\local1
$\RY_-$ is an $M$-surface if and only if $A$ is an $M$-curve and $\RX_+$
is of genus~$0$. In this case $\chi(\RX_+)=\frac14[A]^2\bmod8$.
\item\local2
$\RY_-$ is an $(M-1)$-surface if and only if either $A$ is an
$(M-1)$-curve and $\RX_+$ is of genus~$0$, or $A$ is an $M$-curve and
$\RX_+$ is of genus~$1/2$. In this case
$\chi(\RX_+)=\frac14[A]^2\pm1\bmod8$.
\endroster
\endtheorem

\proof
Let $X$, $Y$, and~$A$ be, respectively, an $(M-d_X)$- and
$(M-d_Y)$-surfaces and an $(M-d_A)$-curve.  Using the exact sequences of
pairs $(\RX,\RA)$ and $(\RX,\RX_+)$, the obvious fact that
$\Gb_*(\RX,\RX_+)=\Gb_*(\RX_-)$, and the expressions for $\Gb_*(\RY_-)$
and $\Gb_*(Y)$ found, respectively, in the left and right hand sides of
the last inequality in~\ref{C2.1.A} (cf.~\ref{C1.2.C} and its proof), one
obtains
$$
d_Y=2d_X+d_A+2(\Gb_0(\RX)-c_--1)+(\Gd+\Ge_+)-2\Gb_1(X),\eqtag\label{eqC3.1}
$$
where $c_-$ is the number of closed components of~$\RX_-$ covered
trivially, and
$$
\gather
\Gd=\dim\Im[\inc_*\:H_1(\RA)\to H_1(\RX)],\\
\Ge_+=\dim\Im[\inc_*\:H_1(\RX_+)\to H_1(\RX)].
\endgather
$$
This gives the conditions for $\RY_-$ to be an $M$- or $(M-1)$-surface.
The congruences follow from the\ basic extremal congruences~\ref{C1.5.A}
applied to~$\RY_-$ and~$\RX$.
\endproof

Note that, unlike the case of plane curves, in general $\RY_-$ and $\RY_+$
cannot both be close to $M$-surfaces. More precisely, arguing as in the
proof of \ref{C3.1.A} one obtains (in the obvious notation)
$d_{Y_-}+d_{Y_+}-\Gb_2(X)=\Gb_2(X)+\Gb_*(A)-(\Gb_*(\RX_-)+
\Gb_*(\RX_+)-2c_0(\RX_-) -2c_0(\RX_+))$, and from~\ref{C2.1.A} it follows
that
$$
d_{Y_-}+d_{Y_+}\ge\Gb_2(X).
$$

Given $X$, $A$ and $Y$ as above, let us discuss whether the covering
$\RY_-\to\RX_-$ is trivial and whether $\RY_-$ is orientable. Since
$A$~is even and ample, there exists a real curve $E\subset X$ such that
$2[E]$ is an odd multiple of~$[A]$. (Indeed, one can divide~$[A]$ by~2 and
replace the result with its large odd multiple in order to make it
effective.) Denote by~$e$ the class of~$[E]$ in $\iH_2(X)$,
see~\ref{C1.7}. Since $H_2(X;\Z)$\mnote{\Dg: edited after Zv} is torsion
free and $\vr_2[\RE]=e$, both~$e$ and the class of~$[\RE]$ in $H_1(\RX)$
are uniquely determined by~$[A]$. Next statement follows immediately from
the definition of the characteristic class of a double covering and the
multiplicative properties of~$\bv_2$, see~\ref{C1.7.B}.

\claim\label{C3.1.B}
Let $X$, $A$, and~$E$ be as above. Then the characteristic class of the
covering $\RY_-\to\RX_-$, considered as a homomorphism
$H_1(\RX_-)\to\Z_2$, is the restriction to $H_1(\RX_-)$ of the map
$x\mapsto\vr_2x\circ e=x\circ[\RE]$.
\qed
\endclaim

\claim[Corollary]\label{C3.1.C}
Let $X$, $A$, and~$E$ be as above. If $\RX_+$ is of genus~$0$ \rom(or of
genus~$1/2$\rom), the projection~$\RY_-\to\RX_-$ is the orientation double
covering if and only if $[\RE]=w_1(\RX)$ \rom(respectively,
$[\RE]=w_1(\RX)+\Gk c$, where $c$ is the class of a one-sided circle
in~$\RX_+$ and $\Gk\in\Z_2$\rom).

If $X$ is an $M$-surface, $[\RE]=w_1(\RX)$ if and only if $\frac12[A]$ is
a characteristic element of~$H_2(X)$.
\qed
\endclaim

Note that $\frac12[A]$ is a characteristic element of $H_2(X)$ if and only
if $Y$ is $\Spin$. This follows from the projection formula for $w_2$ and
the fact that the covering projection induces an epimorphism $H_2(Y)\to
H_2(X)$, see~\ref{C1.1.B}.
Combining~\ref{C3.1.A} and~\ref{C3.1.C} with Nikulin's
congruence~\ref{C2.2.F}
gives the following result:

\claim[Nikulin's congruence for curves]\label{C3.1.D}
Let $X$ be a complex surface with
real structure,
$H_1(X)=0$, and $A$ an ample even real divisor on~$X$. Assume that
\roster
\item"--"
$X$ is an $M$-surface and $A$ is an $M$-curve,
\item"--"
$\RA$~is not empty and lies entirely in a single component of~$\RX$,
\item"--"
the other components of~$\RX$ belong to~$\RX_-$ and are orientable,
\item"--"
$\RA$ is totally homologous to zero in~$\RX$,
\item"--"
$\frac12[A]$ is a characteristic element in $H_2(X)$.
\endroster
Assume, further, that for some integer $r\ge0$ the Euler characteristics
of all the closed components of~$\RX_-$ are divisible by~$2^{r+1}$ and
those of the other components are divisible by~$2^r$ for some
integer~$r\ge0$. Then\rom:
\roster
\item"--"
if
$\RX_+$ is of genus~$0$,
then $\chi(\RX_-)=0\bmod{2^{r+3}}$\rom;
\item"--"
if
$\RX_+$ is of genus~$1/2$,
then $r=0$ and $\chi(\RX_-)=\pm1\bmod{8}$.
\qed
\endroster
\endclaim

\remark{Remark}
Complementary Fiedler's congruence, see~\cite{Fie:cong}, obtained by means
of Rokhlin-Guillou-Marin formula, can
also be generalized from plane curves to
curves on surfaces.
\endremark

\subsection{Generalizations of Rokhlin's formula of complex orientations}\label{C3.2}
Consider a nonsingular complex surface~$X$ with real structure and a
dividing real curve $A\subset X$ without real singular points. Denote by
$A^k$ the $\R$-components of $A$ and assume that each of $A^k$ is
nonsingular over~$\C$. Let $C_1,\dots,C_n$ be the connected components
of~$\RX$, and for each $i=1,\dots,n$ let $C_{ij}$,
$j=\nobreak1,\dots,m_i$, be the orientable components of $C_i\sminus\RA$.
If~$C_i$ is orientable, fix an orientation on it and endow~$C_{ij}$ with
the induced orientation; otherwise pick an orientation of~$C_{ij}$
arbitrarily. Denote by~$[C_i]$ (for orientable~$C_i$) and~$[C_{ij}]$ the
fundamental classes in $H_2(\RX,\RA;\Z)\subset H_2(\RX,\RA;\Q)$
corresponding to the chosen orientations. Obviously, $[C_{ij}]$ form a
basis of $H_2(\RX,\RA;\Z)$ and $[C_i]=\sum_{j=1}^{m_i}[C_{ij}]$ whenever
$C_i$ is orientable.

Since $A$ is dividing, so are its $\R$-components and, for each $k$,
$\RA^k$ separates~$A^k$ into two halves~$A_\pm^k$. Let~$[A_\pm^k]\in
H_2(X,\RA;\Z)$ be the class realized by~$A_\pm^k$ with its complex
orientation. Obviously, $[A_+^k]+[A_-^k]=[A^k]$. (We identify $[A^k]$,
$[A_\pm^k]$, etc\. with their images under various inclusion and
coefficient homomorphisms.) The halves $A_\pm^k$ induce a pair of opposite
orientations on~$\RA^k$, called the \emph{complex orientations}. The
corresponding fundamental classes are denoted by~$[\RA_\pm^k]$.
Let~$[\RA_\pm]=\sum[\RA_\pm^k]$.

Let $\inc\:A\to X$ be the inclusion. The most straightforward
generalization of Rokhlin's formula~\ref{Rokh.form}
(underlined by V.~Zvonilov \cite{Zvonilov:RF} under slightly
stronger assumption) is obtained under
the assumption that
$[\RA_+]$ vanishes in $H_1(\RX;\Q)$, which is equivalent to the existence
of certain rational numbers~$x_{ij}$, $i=1,\dots,n$, $j=1,\dots,m_i$,
such that
$$
[\RA_+]=\partial\sum x_{ij}[C_{ij}].\eqtag\label{weights}
$$
If $C_i$ is nonorientable, $x_{ij}$ are uniquely determined by
$\RA\subset\RX$ and its complex orientation. If $C_i$ is orientable, the
class $\sum_{j=1}^{m_i}x_{ij}[C_{ij}]$ is well defined up to a multiple
of~$[C_i]$. In this case, unless $C_i$ is a torus, there is a unique
choice of~$x_{ij}$ so that $\sum x_{ij}\chi(C_{ij})=0$. (If $C_i$ is a
torus, $\sum x_{ij}\chi(C_{ij})$ is determined by~$\RA$ and its complex
orientation.)

Denote by~$x$ the image of $\sum x_{ij}[C_{ij}]$ in $H_2(X,\RA;\Q)$. There
are unique classes $\xi_+,\xi_-\in H_2(X;\Q)$ taken, respectively, into
$[A_+]-x$, $[A_-]+x$ by the relativization homomorphism. It is clear that
$\xi_++\xi_-=[A]$. Note that, unlike the case of plane projective curves,
these classes may not be determined by the complex orientation of~$\RA$;
thus, they provide a new rigid isotopy invariant of a dividing curve
in~$X$.
According to V.~Zvonilov,\mnote{\Dg: added after Zv} there are examples
when this invariant distinguishes curves with the same complex
orientations.

If $A$ is nonsingular, $A_+\cap A_-=\RA$. In general, $A_+\cap A_-$
differs from $\RA$ by a finite number of points. Denote by $d^{+-}$ the
number of the points, counted with their multiplicities. Clearly, $0\le
d^{+-}\le\frac12\sum_{i<j}[A_i]\circ[A_j]$.

\theorem\label{C3.2.A}
Let $A\subset X$ be as above and $\inc_*[\RA_+]=0$ in $H_1(\RX;\Q)$. Then
$$
 \xi_+\circ\xi_-=d^{+-}+\sum_{i,j=1}^{n,m_i}x_{ij}^2\chi(C_{ij}).
$$
Furthermore, for each orientable component~$C_i$ of~$\RX$
$$
 \xi_+\circ[C_i]=-\xi_-\circ[C_i]=\sum_{j=1}^{m_i}x_{ij}\chi(C_{ij}).
$$
\endtheorem

\proof
The statement is proved similar to Rokhlin's formula, using the
Lagrangian property of $\RX$, (i.e., the fact that the multiplication by
$\sqrt{-1}$ transforms the tangent bundle of $\RX$ into its normal bundle
with the opposite orientation).
\endproof

\claim[Complement]\label{moreint}
One has $\xi_+\circ [A]=\xi_-\circ[A]=\frac12 [A]^2$.
\qed
\endclaim

In the case of projective plane $\xi_+=\xi_-=\frac12[A]$, and replacing
the left hand side of the first equation in~\ref{C3.2.A} by $\frac14[A]^2$
gives true Rokhlin's formula \ref{Rokh.form}
(as well as the Rokhlin-Mishachev formula
for curves of odd degree and Rokhlin-Zvonilov-Viro formula for reducible
curves). In general one cannot assert that $\xi_+=\xi_-$ and give an
implicit form to the above substitute of Rokhlin's formula. From
comparing the exact sequences of pairs $(A,\RA)$ and $(X,\RX)$ it follows
that $\inc_*[\RA_+]=0$ and $\xi_+=\xi_-$ (under an appropriate choice
of~$x_{ij}$) if and only if the kernel of the inclusion homomorphism
$H_2(A,\RA;\Q)\to H_2(X,\RX;\Q)$ is nontrivial. There is another
sufficient condition for $\xi_+=\xi_-$, which follows from the second
equation in~\ref{C3.2.A} and the obvious fact that $\xi_+-\xi_-$ is
$\conj_*$-invariant:

\claim[Proposition]\label{C3.2.B}
If $\inc_*[\RA_+]=0$ in $H_1(\RX;\Q)$ and the $(+1)$-eigenspace~$H^+$
of~$\conj_*$ on $H_2(X;\Q)$ is spanned by the fundamental classes of the
orientable components of~$\RX$, then $\xi_+=\xi_-$, provided that the
coefficients~$x_{ij}$ are chosen so that
$\sum_{j=1}^{m_j}x_{ij}\chi(C_{ij})=0$ for each orientable component~$C_i$
of~$\RX$. \rom(Under the hypotheses the last equation holds automatically
for all toroidal components.\rom)
\qed
\endclaim

If $X$ is K\"ahler, since $\dim H^+=\frac12(b_2(X)+\chi(\RX))-1$, the last
hypothesis of~\ref{C3.2.B} is satisfied, for example, in the following
cases:
\roster
\item
$b_2(X)+\chi(\RX)=2$;
\item
$\RX$ is an orientable $M$-surface, $b_0(\RX)\le2$, and $\RX$ is not
a pair of tori;
\item
$\RX$ is an orientable $(M-1)$-surface, $b_0(\RX)=1$, and $\RX$ is
not a torus.
\endroster

\remark{Remark\rom: Versions of Arnol$\,'\!$d congruence}
Under
a stronger assumption $\inc_*[\RA]=0$ in $H_1(\RX;\Z)$,
the coefficients~$x_{ij}$ in~\eqref{weights} can be chosen integral,
and \ref{C3.2.A} implies that $\chi(B)=\xi_+\circ\xi_-\bmod 4$, where
$B$ with $\partial B=\RA$ is the union of $C_{ij}$ with $x_{ij}$ odd.
On the other hand, if $\inc_*[\RA]=0$ in $H_1(\RX)$ and $B$ is any
piece of $\RX$ such that $\partial B=\RA$, then $\chi(B)=\Cal
P(\xi)-\chi(A)-\frac12[A]^2\bmod4$, where $\Cal P$ is the Pontrjagin
square and $\xi\in H_2(X)$ is the class of $A_++B$. A piece $B$ with
$\partial B=\RA$ exists, e.g., if $\inc_*[A]=0\in H_2(X)$, and in this
case $\chi(B)=\Cal P(\xi)\bmod 4$. If, in
addition, $H_1(X)=0$, then $X$ is $\ZZ$-Galois maximal and
$\xi=\frac12A+\bv_2w_1(B)$, where $B$ is a half of $\RX$ corresponding to
the double covering; thus, $\Cal P(\xi)$ can be evaluated using the
Pontrjagin-Viro form.
\endremark

Another generalization of Rokhlin's formula, suggested by J.-Y.
Welschinger \cite{W1}, uses the assumption $\inc_*[\RA]\in
2l\inc_*[H_1(\RA;\Z)]\subset H_1(\RX;\Z)$, $l\in\Z$. (In fact, the
assumption in~\cite{W1} is a bit more general.) It unifies original
Rokhlin's formula, Arnol$'$d congruence, and some of Mikhalkin's
congruences~\cite{Mikhalkin}. The proof uses the canonical quadratic map
${\Cal P}_{\conj}\:H_2(X;\Z/2l\Z)\to\Z/4l\Z$, which extends the reduction
of the twisted quadratic form $H_2(X;\Z)\to\Z$, $x\mapsto x\circ\conj x$.

Under the above assumption there are integers~$x_{ij}$, $i=1,\dots,n$,
$j=1,\dots,m_i$, such that $[\RA_+]=\partial\sum x_{ij}[C_{ij}]\bmod 2l$.
Denote by~$x$ the image of $\sum x_{ij}[C_{ij}]$ in $H_2(X,\RA;\Z)$.
There are unique classes~$\xi_+,\xi_-\in H_2(X;\Z/2l\Z)$ taken,
respectively, to $[A_+]-x$, $[A_-]+x$ by the relativization homomorphism.
As before, $\xi_++\xi_-=[A]$.

\claim[Theorem]\label{Welsch1}
Given a dividing curve $A$
with $\inc_*[\RA]\in2l\inc_*[H_1(\RA;\Z)]$,
$l\in\Z$, one has
$$
{\Cal P}_{\conj}(\xi_+)=-\sum_{i,j=1}^{n,m_i}x_{ij}^2\chi(C_{ij})\bmod4l.
\rlap\qed
$$
\endclaim

In some cases, e.g., if $H_2(X;\Z)$ is generated by real algebraic curves
with known position in respect to $\RA$, the element $\xi_+$ can be
calculated explicitly in terms of the complex orientation of~$\RA$. Then
\ref{Welsch1} gives some new prohibitions on the complex orientation.
Some precise statements can be found in \cite{W1}.

\subsection{Adjunction formulas}\label{adjunction}
Rokhlin's complex orientation formula is closely related to the Bishop-Lai
generalized adjunction formula stated below: they both count the
algebraic number of virtual double points and, when they both apply, they
lead to the same result; moreover, their geometric proofs go along the
same lines.

The traditional \emph{adjunction formula} for a nonsingular curve $D$ in
a nonsingular surface $X$ states that
$$
2g(D)-2=[D]^2+[K_X]\circ[D],\eqtag\label{eqC3.4}
$$
or,
equivalently,
$$
\chi(D)+D^2=c_1(X)\scap D.
$$
If $D$ is singular, one can interpret the defect
$\12(\chi(D)+D^2-c_1(X)\scap D)$ as the virtual number of double points
of~$D$.

The formula can easily be generalized from algebraic curves to immersed
real $2$-manifolds (see~\cite{Lai}, \cite{EKh}). Let $F$ be a real
oriented compact $2$-manifold (possibly disconnected and with boundary)
and $f\:F\to X$ an immersion into a complex surface~$X$. If $\partial
F\ne\emptyset$, we assume that either $f$ has no complex points on
$\partial F$ (a point $x\in F$ is complex if $df(T_x F)\subset
T_{f(x)}X=\C^2$ is a complex line) or all the points of $\partial F$ are
complex. Pick a nondegenerate vector field $n\:\partial F\to f^*(TX)$
normal to $\partial F$ as follows: if $f$ has no complex points on
$\partial F$, let $n\:\partial F\to TF\subset f^*(TX)$ be normal to
$\partial F$ in~$F$; otherwise, $TF\subset f^*(TX)$ is a complex line
subbundle of $f^*(TX)$ over $\partial F$ and $n\:\partial F\to f^*(TX)$
is taken normal to~$TF$.

Denote by $\nu(f)\in\Z$ the obstruction to extending $\sqrt{-1}\,n$ to a
vector field normal to~$F$, and by $c_1(f)\in\Z$, the obstruction to
extending $(\tau,n)$, where $\tau$ is a nondegenerate tangent field on
$\partial F$, to a complex $2$-frame field in $f^*(TX)$ over $F$. The
latter is the Maslov index of the framing with respect to to the
membrane~$F$. If $\partial F=\varnothing$, then $\nu(f)$ is the normal
Euler number of~$F$ in $f^*(TX)$ and $c_1(f)=c_1(X)\cap f_*[F]$. In this
notation, the Bishop-Lai formula states
$$
\chi(F)+\nu(f)-c_1(f)=2d_-(f),\qquad
\chi(F)+\nu(f)+c_1(f)=2d_+(f),\eqtag\label{general.adjunction}
$$
where $d_\pm(f)$ are determined by the complex points of~$f$. More
precisely, if $f$ is generic, then $d_\pm(f)=e_\pm(f)-h_\pm(f)$, where
$e_\pm$ (respectively, $h_\pm$) are the numbers of positive and negative
elliptic (respectively, hyperbolic) points of~$f$. In general, $d_+(f)$
and $d_-(f)$ can be defined as the intersection numbers of the image
of~$F$ in the Gauss fibration of oriented real $2$-planes in $f^*(TX)$
with the cycle of complex lines equipped by their complex or,
respectively, anti-complex orientation. If a component of $\partial F$
consists of complex points, it is first shifted off the cycle in the
tangent directions $\tau^*\otimes n+\sqrt{-1}\,\tau^*\otimes\sqrt{-1}\,n$
determined by the real $2$-planes~$(\tau,n)$. If $f$ is generic, $e_+$
and~$h_+$ (respectively, $e_-$ and~$h_-$) are the numbers of the positive
and negative intersection points of~$F$ with the complex (respectively,
anti-complex) orientation cycle. A definition of elliptic and hyperbolic
points in local coordinates can be found, e.g., in~\cite{EKh}.

\subsubsection{}\label{piecewise.smooth}
The definition of $\nu(f)$ and $c_1(f)$ and the Bishop-Lai formula extend
to a wider class of surfaces: one can assume that $F$ and $f\:F\to X$ are
sewed (along boundary components) of surfaces $F_i$ and immersions
$f_i\:F_i\to X$ satisfying the above conditions, so that any two adjacent
pieces are normal in $X$. (We do not impose any requirements on the sewn
boundary fields,
as
we do not use additivity relations.) Note also that, if $F$ is closed, the
algebraic number $i(f)$ of virtual double points of~$f$ is given by
$$
i(f)=\frac12\bigl((f_*[F])^2-\nu(f)\bigr)=
\frac12\bigl(\chi(F)+(f_*[F])^2\pm c_1(f)\bigr)-d_\pm(f).
\eqtag\label{double.pts}
$$

The following statement is a straightforward consequence
of~\eqref{general.adjunction}.

\claim\label{mixedcycleformula}
Let $F$ be a closed oriented real $2$-manifold, $X$ a complex surface,
and $f\:F\to X$ a piecewise smooth map as in~\ref{piecewise.smooth} whose
image has no negative complex points. Then
$$
\chi(F)+\nu(f)=c_1(f).\rlap\qed
$$
\endclaim

\claim[Corollary]\label{mixedcor}
In the notation of~\ref{C3.2.A} one has
$$
\xi_\pm^2=\frac12 [A]^2-d^{+-}-\sum_{i,j=1}^{n,m_i}x_{ij}^2\chi(C_{ij}).
$$
\endclaim

\proof
Apply \ref{mixedcycleformula} twice: first
to $\xi$ (or
to each of its smooth pieces), and then to~$A$. (Certainly, the same
result
follows from~\ref{C3.2.A} and~\ref{moreint}, cf\. the beginning
of this
section.)
\endproof

In the case of plane curves
\ref{mixedcor} turns into Rokhlin's
formula.

\subsection{An application of the adjunction inequality}\label{adj<}
Let $A\subset X$ be as in~\ref{C3.2}. Assume that $\inc_*[\RA_+]=0\in
H_1(\RX;\Z)$ and denote by $x_{ij}$ the integers resolving
\eqref{weights}. The cycle~$\xi_+$, see \ref{C3.2}, can be represented as
$f_*[F]$ for some map $f\:F\to X$ and some surface~$F$ sewed of~$A_+$ and
$\ls|x_{ij}|$ copies of~$C_{ij}$. Such a cycle has no negative complex
points and, in view of~\eqref{double.pts} and~\ref{mixedcycleformula},
the number $i_+(f')$ of positive double points of a generic
bifurcation~$f'$ of~$f$ satisfies the \emph{adjunction inequality}
$$
i_+(f')\ge\frac12\bigl(\chi(F)+(f_*[F])^2-c_1(f)\bigr).
$$
The results obtained using the inequality are weaker than Rokhlin's
formula (which, in this case, is equivalent to the Bishop-Lai formula).
However, as it was observed by G.~Mikhalkin~\cite{Mikh2},
new prohibitions on
\emph{reducible} real curves can be obtained applying the adjunction
inequality to similar mixed cycles with the orientation on some of
holomorphic pieces reversed. (For the sake of simplicity
we
ignore the somewhat more cumbersome case where the irreducible components
have real intersection points; on the other hand, we treat curves on
arbitrary K\"ahler surfaces.)

Thus, fix an integer $s\ge1$ and let $f'\:F\to X$ be identical to $f$
over all $C_{ij}$ and all $A^k_+$ with $k>s$ and equal to $\conj\circ f$
over all $A^k_+$ with $1\le k\le s$. Denote by $d_1,d_2\in H_2(X;\Z)$ the
elements realized by $\bigcup_{k>s}A_k$ and $\bigcup_{k\le s}A_k$,
respectively, and let $x^{++}$ be the algebraic intersection number of
$\bigcup_{k>s}A_+^k$ and $\bigcup_{k\le s}A_+^k$ (recall that, due to the
assumption made, see~\ref{C3.2}, the boundaries of the two
chains\mnote{\Dg: edited after Zv} are disjoint).

\claim[Theorem]\label{adjprohi} Let $A\subset X$ be as in \ref{C3.2}
and $\inc_*[\RA_+]=0$ in $H_1(\RX;\Z)$.
If $X$ is K\"ahler and $\omega\cap(d_1-d_2)>0$, where $\omega$ is the
K\"ahler class, then
$$
j_-\ge x^{++}+\frac12(c_1(X)\scap d_2-d_1\circ d_2),
$$
where $j_-=\sum\12\ls|x_{ij}|(\ls|x_{ij}|-1)$, the sum taken over all
$\ls|x_{ij}|\ge 2$ corresponding to the components $C_{ij}$ homeomorphic
to a disc.
\endclaim

\proof
The inequality is the difference between \ref{mixedcycleformula} applied
to $f$ and the adjunction inequality applied to the usual Lagrangian
perturbation of~$f'$, which provides an immersion with $j_-+\12(d_1\circ
d_2)-x^{++}$ negative double points.
\endproof

\claim[Corollary]\label{adjcor}
Let $X$, $d_1$, $d_2$ be as above and $\omega\cap(d_1-d_2)>0$
\rom(see~\ref{adjprohi}\rom). Assume that the $(+1)$-eigenspace of
$\conj_*$ on $H_2(X;\Q)$ is spanned by the fundamental classes of the
orientable components of $\RX$ and that for every~$k$ each connected
component of $\RA^k$ bounds a disc. Then
$$
j_-\ge z^+-z^-+\frac12 c_1(X)\scap d_2-\frac14(d_1\circ d_2),
$$
where $z^\pm$ are the numbers of positive and negative injective pairs of
ovals such that one of the ovals is in $A_k$ with $k\le s$ and the other
one, in $A_k$ with $k>s$. \rom(In order to define injective pairs for each
connected component of~$\RA$ one picks and fixes a disk bounded by this
component.\rom)
\endclaim

\proof
Construct the Rokhlin cycles by
filling\mnote{toZv: changement} $\bigcup_{k>s}A_+^k$ and $\bigcup_{k\le
s}A_-^k$ with disks. Then $x^{++}$ is found by calculating the
intersection of the cycles in two ways: homologically
(applying~\ref{C3.2.B}) and geometrically.
\endproof

The following example of application of~\ref{adjcor} (combined with
Rokhlin's formula) is also due to Mikhalkin. Let $A$ be a real plane
curve of degree $2k$ with the arrangement
$\langle\smash{\frac32}k(k-1)\rangle\+1\langle\12(k-1)(k-2)\rangle$
(Harnack arrangement), and $C$ a real conic with $\R C$ disjoint from
$\RA$. Then {\proclaimfont at most $(k-2)$ of the interior ovals of~$\RA$
are inside the oval of~$\R C$.}

\subsection{Orevkov's formula}\label{stepa}
S.~Orevkov~\cite{Or} observed that, if a plane $M$-curve of degree $d$
has a nest of depth $k-1$, $k=[d/2]$, Rokhlin's formula splits into two
separate identities. Under the assumption the nonempty ovals are pairwise
nested and, thus, are divided into \emph{positive} and \emph{negative}
(e.g., with respect to the outermost oval of the nest). Denote by
$\pi^+_+$ and $\pi^-_+$ (respectively, $\pi^+_-$ and $\pi^-_-$) the
number of positive and negative injective pairs whose inner oval is empty
and outer oval is positive (respectively, negative).

\claim[Orevkov's formula]\label{orev:COF}
For any plane $M$-curve of degree $2k$ with a nest of depth $k-1$ one has
$$
\pi^-_+-\pi^+_+=k^2_+,\quad \pi^+_- - \pi^-_-=k^2_-,
$$
where $k_+$ and $k_-$ are the number of positive and negative nonempty
ovals.
\endclaim

The formula above and a similar formula for curves of odd degree give a
new information on the complex orientation in degrees~$\ge7$. Both the
formulas extend to nonmaximal dividing curves.

The original proof (see~\cite{Or}) is based on counting the self-linking
numbers of the braids cut by $\C A$ on the boundary of the two $4$-balls
formed by the imaginary points of the imaginary lines through a point in
an innermost oval of the nest. \hbox{J.-Y.~Welschinger}~\cite{W2} observed
that Orevkov's formula can also be obtained from the adjunction formula
applied to appropriate mixed singular cycles and generalized the formula
to curves on conic bundles. For simplicity, we state his result for
geometrically ruled surfaces,
i.e., conic bundles without singular fibers.

Let $\pi\:X\to B$ be a holomorphic equivariant map of a real surface to a
real curve. Assume that $B$ is a dividing curve and each fiber~$F$
of~$\pi$ is a smooth irreducible rational curve. Assume, in addition,
that $\pi$ has a real section~$E$. The fundamental classes of~$F$ and~$E$
form a basis in $H_2(X;\Z)$.

Denote by $B_1,\dots,B_p$ the connected components of $\R B$ and by
$X_1,\dots,X_p$, their pull-backs in $\RX$. (Each $X_j$ is a torus or a
Klein bottle.) Pick a half~$B_+$ of~$B$ and equip $\R B$ with the
corresponding complex orientation.

Let $A\subset X$ be a nonsingular irreducible dividing real $(M-2r)$-curve
with a contractible oval on each real component of $X$ (the latter
assumption is not essential for \ref{Welsch2}). Pick a half~$A_+$ of~$A$
and equip $\RA$ with the corresponding complex orientation. For each $i$,
$1\le i\le p$, consider the components of $\RA$ which project to~$B_i$
with a positive (respectively, negative) degree and denote by $n^+_i$
(respectively, $n^-_i$) their total degree. For each~$i$ there is a
positive integer~$k_i$ such that $n^+_i+n^-_i+2k_i=[A]\circ [F]$. If
$k_i=1$ for each~$i$, $1\le i\le p$, we say that $A$ is a \emph{curve
with a deep nest}. Clearly, {\proclaimfont for a curve with a deep nest
the weights $n^+=n^+_i$ and $n^-=n^-_i$ do not depend on $i$}.

From now on assume that $A$ is a curve with a deep nest. Then each real
fiber of~$\pi$ meets~$A$ in at most two imaginary points. If a real
fiber~$F$ does meet~$A$ at two imaginary points, pick a point $x\in \R F$
not in~$E$. Such a choice can be made continuously. Moreover, if all the
\emph{turn points} (i.e., the points of tangency of~$\RA$ and a fiber)
lie outside $E$, there is a canonical choice: trace the circle~$\gamma$
in~$F$ through the two imaginary points and the point of~$E$ and take
for~$x$ the intersection point of~$\gamma$ and ~$\R F$ other than the
point of $E$. The chosen points form a finite collection of arcs~$S_j$,
each connecting two turn points. The union of $S_j$ is completed to a
closed, not necessarily connected or simple, curve $\Gamma^+$
(respectively, $\Gamma^-$) by arcs of~$\RA$ projecting to~$\R B$ with
positive (respectively, negative) local degree. The complex orientation
of $\RA$ defines a coorientation of~$S_j$ at its boundary points; it
extends to a continuous coorientation of $S_j$. Denote the extended
coorientation by~$\omega$.  Without lost of generality one can assume
that the intersection points of the interiors of~$S_j$ and~$\RA$ are
transversal, and attribute to each point sign~$+1$ if $\omega$ coincides
with the orientation of~$\RA$ and $-1$ otherwise. Let $\Lambda^{\pm}$ be
the sum of the signs over all the double points of~$\Gamma^\pm$.

To make $\Lambda^\pm$ more invariant, remove the zig-zags on~$\RA$ (each
zig-zag changes one of $\Lambda^\pm$ by $n^\pm-1$ or~$1$, depending on
whether it crosses~$\RE$ or not). Adjusted so, $\Lambda^\pm$ are
preserved by the isotopies of nonorientable components of~$\RX$ identical
on~$\RE$ and on each component of~$\RA$ intersecting~$\RE$. More
generally, $\Lambda^\pm$ are preserved by an isotopy of nonorientable
components of~$\RX$ identical on~$\RE$ if during the isotopy the turn
points do not cross~$\RE$ (each crossing changing
both $\Lambda^\pm$ by~$n^\pm$). An
isotopy of an orientable component of~$\RX$ preserves~$\Lambda^\pm$,
e.g., if during the isotopy (as well as at the initial position) $\RA$
does not intersect~$\RE$.

\theorem\label{Welsch2}
Under the above assumption there are integers $\mu^\pm$, $g^\pm$, $a^\pm$
such that
$$
\Lambda^{\pm}+n^\pm a^\pm+\mu^\pm-g^\pm=
 \frac12(\chi(B)- n^\pm E^2)(n^\pm+1),
$$
and
$$
\gather
g^\pm\ge c^\pm,
\quad g^++g^-+2-\mu^+-\mu^+=c^++c^-+2r,\\
\mu^\pm\ge 1, \mu^++\mu^-\le p+1,\\
a^++a^-=a,\quad a^\pm+E^2(n^\pm+1)=\RE\circ\Gamma^\pm\bmod2,\\
\ls|\RE\cap\Gamma^\pm|\le a+E^2(n^\pm+1)\le
 a^\pm+E^2b-\ls|\RE\cap\Gamma^\mp|,
\endgather
$$
where $a[F]+b[E]=[A]$ and $c^+$, $c^-$ are the numbers of the components
of~$\RA$ projecting to a component of $\R B$ with a positive
\rom(respectively, negative\rom) degree.
\endtheorem

\proof
Consider the surface $W^{\pm}=(A^\pm\cap\pi^{-1}(B^+))\cup
(A^\mp\cap\pi^{-1}(B^-))$ and fill its holes with embedded discs
contained in the pull-back of $\R B$, with the segments $S_j$ as
diameters (e.g., one can use arcs of the circles described in the
construction of $S_j$). The statement follows from the adjunction formula
\ref{mixedcycleformula} applied to the resulting singular cycle
$\Sigma^\pm$, using $\nu=(\Sigma^\pm)^2+2\Lambda^\pm$, $\Sigma^\pm\circ
F=n^\pm+1$,
introducing $a^\pm$ by $[\Sigma^\pm]=a^\pm[F]+(n^\pm+1)[E]$ and taking
for $g^\pm$ and $\mu^\pm$ respectively the genus and the number of
connected components of $W^{\pm}$.
\endproof

\claim[Corollary]\label{WelschC}
If $A$ as in~\ref{Welsch2} does not intersect~$E$, one has
$$
\Lambda^{\pm}+\mu^\pm-g^\pm=\frac12(\chi(B)+E^2n^\pm)(n^\pm+1).\rlap\qed
$$
\endclaim

\remark{Remark}
To obtain~\ref{orev:COF} from~\ref{WelschC}, blow up a point in the
innermost oval of a nest of depth $k-1$; then one has $\chi(B)=2$,
$E^2=-1$, $\mu^\pm=1$, $g^\pm=K_\pm$, and
$\Lambda^\pm=\Pi^-_\pm-\Pi^+_\pm$, where $K$ and $\Pi$ differ from $k$
and $\pi$ in~\ref{orev:COF} by the fact that the innermost oval chosen
contributes to~$K$ and does not contribute to~$\Pi$, whereas it does not
contribute to~$k$ and contributes to~$\pi$.
\endremark

\subsection{The Bezout Theorem}\label{C3.3}
The Bezout theorem plays an important r\^ole in theory of real plane
projective algebraic curves. For an arbitrary surface~$X$, in its complex
form this famous theorem states that {\proclaimfont given two
curves~$A,B\subset X$ without common component, the number of
intersection points of~$A$ and~$B$, counted with their multiplicities,
equals $[A]\circ[B]$}. An important fact is that the multiplicity of an
intersection point of two analytic curves without common components is a
positive integer. This gives the following real version of the Bezout
theorem (cf. a remark in~\ref{C1.7}): {\proclaimfont if $X$ is real and
$A,B\subset X$ are real curves without common components, then the number
of intersection points of~$\RA$ and~$\R B$, counted with their
multiplicities, does not exceed $[A]\circ[B]$ and is congruent to
$[A]\circ\nobreak[B]\bmod\nobreak2$.}

One of the typical applications of the Bezout theorem is tracing an
auxiliary curve~$D$ through a certain number of appropriately chosen
points and estimating the number of intersection points of~$\R D$ and the
given curve~$\RA$ (it is in this way that Harnack obtained his bound for
the number of real components of plane curve). Sometimes one also needs
to know that
$\R D$\mnote{\Dg: changed after Zv} is connected or, more generally, that
all the intersection points belong to the same component of $\R D$. The
simplest way to assure this is to consider rational curves.

Thus, in order to apply the Bezout theorem on an arbitrary surface~$X$ one
needs to know, first, whether a given class $x\in H_2(X;\Z)$ can be
realized by a (formal linear combination of) analytic curves, second, the
dimension of the space of such curves, and third, the genus of such
curves. The answer to the first question in the case when $X$ is
projective is given by the Lefschetz-Hodge theorem
(see~\cite{Kodaira-Spencer}): {\proclaimfont a homology class $x\in
H_2(X;\Z)$ can be realized by a formal linear combination of algebraic
curves on~$X$ if and only if the class Poincar\'e dual to~$x$ belongs
to~$H^{1,1}(X)$} (see~\ref{C1.3}). Though the dimension $h^{1,1}(X)=\dim
H^{1,1}(X)$ for surfaces is a topological invariant,
$h^{1,1}(X)=\Gs^-(X)+1$, the position of~$H^{1,1}$ in $H^2(X;\Z)$ depends
on the complex structure on~$X$ and may change when $X$ varies within its
deformation type. Thus, in general one can only assert that the rational
multiples of the hyperplane section class (those which are still integral
classes) can be realized by combinations of algebraic curves. However, if
$h^{1,1}(X)=b_2(X)$, i.e., $\Gs^+(X)=1$, all classes in~$H_2(X;\Z)$ are
algebraic. Examples are rational surfaces and Enriques surfaces.

The dimension of a linear system~$\ls|D|$ can, in certain cases, be
estimated using the \emph{Riemann-Roch theorem}:
$$
\dim\ls|D|+\dim\ls|K_X-D|\ge
\frac12\bigl([D]^2-[K_X]\circ[D]\bigr)+p_a(X)-1,\eqtag\label{eqC3.2}
$$
where $p_a(X)=\frac12(\Gs^+(X)-b_1(X)+1)$ is the \emph{arithmetic genus}
of~$X$. If $D$ is ample (recall that for an effective divisor on a
surface this is equivalent to $[D]^2>0$ and $[D]\circ[D']\ge0$ for all
irreducible components of~$D$), then $\ls|K_X-D|=\varnothing$ and
in~\ref{eqC3.2} equality is attained:  $$
\dim\ls|D|=\frac12\bigl([D]^2-[K_X]\circ[D]\bigr)+p_a(X).\eqtag\label{eqC3.3}
$$

The following statement is an easy consequence of \eqref{eqC3.3}.

\claim[Proposition]\label{interpltn}
Let $D$ be an ample real divisor and
$s_D=\frac12\bigl([D]^2-[K_X]\circ[D]\bigr)+p_a(X)$. Then through any
collection of $s_D$ points one can trace a real curve in $\ls|D|$.
Moreover, for any $(s_D-s)$-point set~$S$, $s\ge0$, the real curves in
$\ls|D|$ containing~$S$ form a real linear system of dimension~$\ge s$.
\qed
\endclaim

Unfortunately, this result is not sufficient even to generalize the well
known Hilbert bounds\mnote{\Dg: edited after Zv} to the higher numbers of
nests. The difficulty is that the number of connected components of a
curve given by~\ref{interpltn} is not easy to control. In general one can
only claim that this number is between~$1$ and
$2+\frac12([D]^2+[K_X]\circ[D])$ i.e., the Harnack bound. Thus, it would
be useful to know, in addition, through what number of points one can
trace a connected component of a curve in $\ls|D|$. (This and other
related questions, more or less explicitly raised in~\cite{Ro2}, seem to
be of an independent interest, regardless of bounds on the nests).

Denote by $c(D)$ the maximal number $c$ such that through any $c$ points
of $\R X$ in general position one can trace a real curve in $\ls|D|$ with
connected real point set. Denote, further, by $c'(D)$ the maximal number
$c$ such that through any $c$ points of $\R X$ (not necessarily in general
position) one can trace a connected component of the real point set of a
real curve in $\ls|D|$.

\claim[Proposition]\label{c'-c}
$c'(D)\ge c(D)$ for any divisor $D$.
\endclaim

\claim[Lemma]\label{hsdorff}
Let $\{F_n\}$ be a sequence of continuous functions on a Hausdorff compact
space $Y$ converging uniformly to a function $F$. If the zero set $Z_n$
of $F_n$ is connected for each $n$, then the zero set $Z$ of $F$ has a
connected component which contains $\bigcap_n Z_n$.
\endclaim

\proof
Assume that $\bigcap_n Z_n\ne\varnothing$. It is contained in~$Z$. Pick a
point $z\in\bigcap_n Z_n$ and an open-closed subset~$C$ of~$Z$
containing~$z$. Since $Z$ is closed in~$Y$, so are~$C$ and $Z\sminus C$
and, hence, they are separated by a neighborhood~$U$ (i.e, $U\supset C$
and $\Cl U\cap(Z\sminus C)=\varnothing$).

Since the convergence $F_n\to F$ is uniform, there is an integer~$N$ such
that $Z_n\subset U\cup(Y\sminus\Cl U)$ for all $n>N$. Then, since all
$Z_n$ are connected and intersect $C\subset U$, each~$Z_n$ with $n>N$ is
contained in~$U$ and, hence, $\bigcap_n Z_n$ is contained in $C=U\cap Z$.
It remains to notice that in a compact Hausdorff space the intersection
of all open-closed subsets containing a given point is the component of
the point.
\endproof

\proof[Proof of \ref{c'-c}]
Let $c$ be a number such that through any generic set of $c$ points of
$\R X$ one can trace a real algebraic curve in $\ls|D|$ with connected
real point set. We need to prove that an arbitrary $c$-point set in $\R
X$ is contained in a connected component of a real curve in~$\ls|D|$.

Consider a $c$-point set $\Sigma\subset\R X$. By the assumption, there is
a sequence~$\{\Sigma_n\}$ of $c$-point sets in~$\R X$ converging to
$\Sigma$ and such that each~$\Sigma_n$ is contained in a real connected
curve~$\R A_n$ in~$\ls|D|$. Let $\xi$ be the $\R$-bundle over $\R X$
corresponding to $\ls|D|$. Its tensor square, as the square of any
$\R$-bundle, is trivial. Let $F_n$, $n\ge1$, be a real section of the
$\xi$ representing~$A_n$. The space of sections
is\mnote{\Dg: edited after Zv} finite dimensional. Hence, after a proper
scaling $\{F_n\}$ contains a convergent subsequence. As usual, one can
assume that $\{F_n\}$ converges itself. It remains to apply
Lemma~\ref{hsdorff} to~$F_n^2$, considered as functions on~$\R X$.
\endproof

The following additivity property is an
immediate consequence of the definitions.

\claim[Proposition]\label{additivity}
For any divisors $D_1, D_2$ one has
$$
 c(D_1+D_2)\ge c(D_1)+c(D_2)-1\quad\text{and}\quad
 c'(D_1+D_2)\ge c'(D_1)+c'(D_2)-1.
\qed
$$
\endclaim

Some other related problems are also worth mentioning. First, it is
natural to consider curves with a bounded number of components.  Then
$c(D)$ becomes the first member $c_1(D)$ of a larger family of
characteristics $c_{\le k}(D)$, defined similar to $c(D)$ with connected
curves replaced by curves with at most $k$~connected components. Second,
one can consider the minimax $\nu(r,D)$ and maximin $\mu(r,D)$ of the
number of connected components in a linear subsystem of given
dimension~$r$ in~$\ls|D|$.  Another possibility is to consider other
topological characteristics of a linear system, say, typical Morse
numbers, etc.

\claim[Marin's Bound]\label{marin.gen}
If $D$ is very ample and $\chi(X)+D^2$ is odd, any generic pencil of real
curves in
$\ls|D|$\mnote{\Dg: changed after Zv} contains at least one $(M-j)$-curve
with $j\ge 3$.
\endclaim

This result can be restated as $\mu(1,D)\le M(D)-3$, where $M(D)$ is the
Harnack number $2+\frac12([D]^2+[K_X]\circ[D])$.

\claim[Klein's statement]\label{nodevelop}
If a real curve of type~\rom{I} undergoes a Morse surgery through a
nondegenerate double point, its number of connected components can not
increase.
\endclaim

Note that the number of connected components of the curve does not need
to change, cf\.~\ref{typechange} below. In other words, Klein's statement
means that for any~$q$ the number of connected components, considered as
a piecewise constant function on the space $\R C_q$ of curves of degree
$q$, attains a `local
maximum'\mnote{\Dg: edited after Zv} at each component of $\R C_q\sminus\R
D_q$ (where $\R D_q$ is the hypersurface of singular curves) formed by
real curves of type~I. Using complex orientations one can check that, in
addition to~\ref{marin.gen} and~\ref{typechange}, there is no wall (i.e.,
principal stratum of~$\R D_q$) adjacent to two components formed by
curves of type~I. (There is an exception: $q=2$ and the wall
corresponding to turning the oval inside out; however, this wall is not
adjacent to two different components!)

\remark{Remark}
There are extensions of~\ref{nodevelop} to surfaces
(see~\cite{Viro.congres}) and to hypersurfaces in a projective space
(see~\cite{Kalinin}).
\endremark

\proof[Proof of \ref{nodevelop}]
Let $A_t$, $t\in[-\Ge,\Ge]$, be a family of real curves such that: $A_0$
is a curve with a double nondegenerate point, $A_{<0}$ are of type~I, and
the number of connected components decreases from $\R A_{>0}$ to $\R
A_{<0}$. The last assumption implies that the singular connected
component of $\R A_0$ lifted to the normalization of~$A_0$ is a
topological circle. In particular, the number of real components of the
normalization of~$A_0$ is equal to the number of components of $\R
A_{<0}$. On the other hand, the genera of~$A_0$ and~$A_{\ne 0}$ differ
by~1: $g(A_0)=g(A_{\ne0})-1$. Thus, due to, e.g., Klein's congruence
(which states that the number of real components of a dividing curve has
the same parity as the genus, cf.~Introduction), $A_0$ is not a curve of
type~I. Therefore, two conjugate points in $A_0\sminus\R A_0$ can be
joined by a path in $A_0\sminus\R A_0$. Hence, such a path exists also in
each close curve $A_t$, $t\ne0$, which contradicts to the assumption that
the curves $A_t$ with $t<0$ are of type~I.
\endproof

For $(M-2)$-curves Klein's statement can be strengthened.

\claim[Proposition]\label{typechange}
If two $(M-2)$-curves are obtained from each other by a Morse surgery
through a double nondegenerate point, one of them is of type~\rom{I} and
the other one is of type~\rom{II}.
\endclaim

\proof
Denote by $A_t$, $t\in[-\Ge,\Ge]$, a family of real curves such that:
$A_0$ is a curve with a double nondegenerate point and $A_{\ne 0}$ are
$(M-2)$-curves.  Since the number of connected components does not change
from $\R A_{>0}$ to $\R A_{<0}$, the singular connected component of $\R
A_0$ lifted to the normalization of~$A_0$ turns into two topological
circles. In particular, the number of real components of the
normalization of~$A_0$ is greater by~1 than the number of components
of~$\R A_t$ with~$t\ne0$. On the other hand, $g(A_0)=g(A_{\ne 0})-1$.
Thus, the normalization is an $M$-curve and, hence, it is of type~I, and
so is its perturbation coherent with the complex orientation (while the
other one is of type~II).
\endproof

\proof[Proof of \ref{marin.gen}]
Consider graph~$\Gamma$ with 4 vertices representing $M$-, $(M-1)$-, and
$(M-2)$-curves of type~I and $(M-2)$-curves of type~II. Two vertices are
connected by an edge if there are adjacent curves representing them. Due
to~\ref{nodevelop} and~\ref{typechange}, $\Gamma$ has only three edges:
they connect ~$M$ and $M-1$, \ $M-1$ and $M-2$ of type~II, and $M-2$ of
type~II and $M-2$ of type~I. In particular, $\Gamma$ supports no cycle of
odd length.

A generic pencil of real curves in~$\ls|D|$, considered as a line in the
corresponding projective space, intersects the discriminant hypersurface
of singular curves transversally, the intersection points representing
surgeries of the real parts of the curves.  The degree of the
discriminant has the same parity as $\chi(X)+D^2$.
Hence, if the pencil does not contain an $(M-j)$-curve with $j\ge3$, it
defines an odd length cycle in the graph constructed above, which is a
contradiction.
\endproof

To conclude this section we apply the above results to derive some bounds
on nests. The first theorem is stated so that it depends only on a lower
bound for~$c'(D)$.

\claim[Proposition]\label{disjoint nests}
If a real curve ~$C$ has $c\le c'(D)$ disjoint nests of depths
$h_1,\dots,h_c$, so that no oval of the nests of depth $>1$ envelopes all
the other ovals of the nests, then
$$
h_1+\dots+h_c+[\12(c'(D)-c)]\le\12 C\circ D.
$$
In particular, if~$C$ has $c\le c'(D)-1$ disjoint nests of depths $h_1,
\dots, h_c$, so that no oval of the nests of depth $>1$ envelopes all the
other ovals of the nests, and
$$
h_1+\dots+h_c+[\12(c'(D)-c)]=[\12C\circ D],
$$ then the union of the nests exhausts the whole curve.
\endclaim

\proof
Take a point on (or inside, if $h>1$) the innermost oval of each nest  and
$c'(D)-c$ more points on one of these ovals.  Trace a connected component
of a real curve of~$\ls|D|$ trough the $c'(D)=c+(c'(D)-c)$ points obtained.
It intersects the original curve in at least $2(h_1+\dots+h_c)+r$ points,
where $r=2\bigl[\frac12(c'(D)-c)\bigr]$. Now the bound follows from the
Bezout theorem, and the extremal property is its direct consequence.
\endproof

If the union of the nests does not exhaust the real part of the curve,
the estimate can be sharpened slightly, by choosing the additional
$c'(D)-c$ points on an oval not in the nests.

\claim[Proposition]\label{nests outside}
Let a real curve~$C$ have $s=\frac12D\circ(D-K)-l$ nests outside each
other. If $s>\mu(l,D)$, then there are $s-\mu(l,D)+1$ nests \rom(among
the given ones\rom) whose total weight is at most $\frac12C\circ
D-(\mu(l,D)-1)$.
\endclaim

\proof
Take a point on the innermost oval of each nest and consider the linear
system of curves in~$\ls|D|$ through these points. Its dimension is at
least~$l$, and by the definition of $\mu(l,D)$ the system contains a
curve~$D_0$ with $\le\mu(l,D)$ connected components. Hence, at least
$s-\mu(l,D)+1$ nests contribute at least twice their depth to the
intersection with~$D_0$. Each other nest contributes at least two, by its
innermost oval. It remains to apply the Bezout inequality.
\endproof

\claim[Corollary]\label{total weight}
If
a curve~$C$ has $\12D\circ(D-K)$
nests outside each other, then among these nests there are $r=-D\circ
K-1$ nests whose total weight
is at most $\frac12 C\circ D-\frac12D\circ(D+K)-1$.
\endclaim

\proof
This is a special case of~\ref{nests outside} corresponding to $l=0$.
Indeed, $\mu(0,q)=\frac12D\circ(D+K)+2$ and
$r=\frac12D\circ(D-K)-\mu(0,D)+1$.
\endproof

\subsection{Plane curves through given points}\label{s16.1}
In the case of plane curves a linear system $\ls|D|$ is determined by the
degree $q$ of $D$. Thus, one deals with sequences $\{c(q)\}$,
$\{c'(q)\}$. Their first four terms are known.

\claim[Proposition]\label{16.1.BB}
One has $c(1)=c'(1)=2$, $c(2)=c'(2)=5$, $c(3)=c'(3)=8$, and
$c(4)=c'(4)=13$.
\endclaim

Note that, in fact,
in any real pencil of
cubics there is a real cubic with the connected real part.

\proof
The real point set of any real curve of degree 1 or 2 is connected.
Due to~\ref{c'-c} it remains to prove that
$$
8\le c(3),\, c'(3) < 9,\quad 13\le c(4),\, c'(4) < 14.
$$

The bound $c(3)\ge 8$ follows from \ref{rat.cubic}. To show that
$c'(3)<9$, take an irreducible real curve~$A$ of degree~3 with two real
components and pick 5 points on its two-sided component and 4 points on
the one-sided component. By Poincar\'e duality, any real curve~$A'$ of
degree~3 through the 9~chosen points has a fifth intersection point with
the one-sided component (or, as a special case, one of the four chosen
points is multiple).  The Bezout theorem implies then that $A'$ coincides
with $A$, but no component of $A$ contains all the 9~points.

The bound $c(4)\ge 13$ is a straightforward consequence of
\ref{marin.gen}.

To show that $c'(4)<14$, take a real curve~$A$ of degree~4 with 4~real
components and pick 11 points on one of the ovals and one point on each
of the other ones. By Poincar\'e duality, any real curve~$A'$ of degree~4
through the 14 chosen points intersects each oval of~$A$ at least one
more point (or has a multiple intersection point). Thus, the total number
of intersection points is $>16$, and the Bezout theorem implies that $A'$
coincides with $A$ (as a curve of degree $\ge 4$ with the maximal number
of components, $A$ is nonsingular and irreducible). It remains to notice
that none of the component of $A'=A$ contains all the 14 points.
\endproof

\remark{Remark}
In the above proof of $c'_3<9$ the 9~points are chosen so that they do
not lie on a real cubic other than~$A$. From the proof it follows also
that a connected real cubic can be traced through 9~points of an
$M$-cubic~$C$ if and only if they belong to a cubic~$C'$ other than~$C$.
The existence of~$C'$ can be expressed in terms of the group structure
on~$C$: it exists if and only if the sum of the 9~points equals~$0$ (or,
more precisely, 9~times an inflection point). Certainly, the set of
9~points should be invariant with respect to the complex conjugation;
then $C'$ can be chosen real. As a consequence, the 9-tuples {\bf not} on
a connected cubic form an open nonempty subset in the space of all
9-tuples of points of~$\Rp2$.
\endremark

The precise values of the other terms of~$\{c(q)\}$, $\{c'(q)\}$ are not
known. Below we give a few simple general lower and upper bounds.

\claim[Proposition]\label{16.1.C}
For any $q$ one has $c(q)\ge3q-1$.
\endclaim

\subsubsection**{Digression: real rational curves}
As far as we know, the following problem is still open: {\proclaimfont is
it possible to trace an irreducible real rational curve \rom(or, more
precisely, its connected component\rom) of degree $q$ through any set of
$3q-1$ real points in general position.} In~\cite{Ro2} the question
is answered in the affirmative; however, the proof has been never
published; probably, it contained a gap. Alternatively, one can lift the
general position assumption and pose the question of existence of some
(not necessarily irreducible) rational curve.

The first nontrivial case (and the only one where the complete answer is
known) is $q=3$: through $8$ generic points one can trace $12$ cubics;
depending on the position of the points, $8$, $10$, or $12$ of them are
real. (All the three values occur; all the $12$ rational cubics in a
pencil are real if and only if the pencil contains two cubics with a
solitary double point.)

There also are a few simple results on collections of points with
prescribed multiplicities. For example, if $q-1$ of the $3q-1$ points
coincide, the answer is in the affirmative. In fact, in an affine
coordinate system such that the multiple point, $M$, and one of the other
$2q$\mnote{\Dg: changed after Zv} points, $N$, lie at infinity, the
problem becomes linear: if the $x$-axis contains~$N$ and the $y$-axis
contains~$M$, it reduces to interpolating $2q-1$ values $y_j$ by a
rational function $y=P(x)/Q(x)$ of multidegree $(\deg P,\deg
Q)=(q-1,q-1)$. The explicit solution is called \emph{Cauchy interpolation
formula}.
As this construction shows, through {\bf any} $2q$~real points one can
trace a real rational curve of degree~$q$. The result can also be
obtained by a proper small perturbation of a family of $q$~lines through
the given points, considering in the space $\R C_q$ of real curves of
degree $q$ the subvariety of curves having $\frac12(q-1)(q-2)$ double
points and passing through the given points. Replacing lines by conics or
rational cubics yields slightly better results: a real rational curve of
degree~$q$ can be traced through any $\frac52q$ or, respectively,
$\frac83q$ real points.

\claim[Proposition]\label{rat.cubic}
Given 8 generic points in $\Rp2$, denote by $c_+$ and $c_-$ the numbers
of real cubics through them with, respectively, a solitary double point
and a real crossing point. Then $c_--c_+=8$ and $c_-+c_+\le12$. In
particular, there are at least $8$ real rational cubics through the given
points.
\endclaim

\proof
The pencil of cubics through the 8 given points defines a (singular)
fibration $X\to\Rp1$, where $X$ is $\Rp2$ blown up at the 8~points, and
the statement follows from counting the Euler characteristic of~$X$.
\endproof

Note that in a similar way one can prove that in any pencil of cubics
there is one with connected real part.

\proof[Proof of \ref{16.1.C}]
In view of~\ref{16.1.BB} it suffices to consider $q>4$. The proof depends
on the residue $q\bmod4$.

If $q=4k+1$, $k\in\Z$,  then $3(4k+1)-1=12k+2$. Trace a real line through
any two of the $12k+2$ points. The other $12k$~points decompose into
$k$~groups of 12~points.  Each group and one of the first 2~points can be
joined by a real curve of degree~4 with connected real point set. The
product of the line and the constructed curves of degree~4 is of degree
$4k+1$; its real point set is connected and contains all the $12k+2$
points.

If $q=4k+2$, $k\in\Z$,  then $3(4k+2)-1=12k+5$. Trace a real conic through
any five of the $12k+5$ points. The other $12k$ points decompose into
$k$~groups of 12~points. Each group and one of the first 5~points can be
joined by a real curve of degree~4 with connected real point set.  The
product of the conic and the constructed curves of degree~4 is of degree
$4k+2$; its real point set is connected and contains all the $12k+5$
points.

If $q=4k+3$, $k\in\Z$,  then $3(4k+3)-1=12k+8$, and a curve is
constructed as above, starting with a cubic through any eight of the
$12k+8$ points.

If $q=4k$, $k\in\Z$, then the bound $12k-1$ can be replaced with $12k+1$.
Pick one of the points. The other $12k$ points decompose into $k$~groups
of 12~points. Each group and the fixed point can be joined by a real
curve of degree~4 with connected real point set.  The product of these
curves has degree $4k$, its real point set is connected and contains all
the $12k+1$ points.
\endproof

\claim[Proposition]\label{16.1.D}
If $q$ is odd, then $c'(q)\le\12q(q+3)-[\18(q^2-1)]$. If $q$ is even, then
$c'(q)\le\12q(q+3)-[\18q(q-2)]$.
\endclaim

\proof
First, consider the case of $q$~odd, $q\ge3$. (For $q=1$ the inequality
follows from $c'(1)=2$). Put $q=2k-1$ and $m=k+1$. Then the bound takes
the form $c'(q)\le\frac32k(k+1)-1$. Take an $M$-curve~$A$ of degree~$m$.
Its real part has $\frac12{k(k-1)}+1$ connected components. Pick $3k$
points on one of the components and 3~points on each of the others,
totally $\frac32 k(k+1)$ points. If $k$ is even, the distinguished
component should be one-sided. Suppose that $c'(q)\ge\frac32k(k+1)$. Then
there is a real curve $B$ of degree $q$ with a connected component which
contains all the $\frac32k(k+1)$ points. By the Poincar\'e duality, each of
the components of~$\R A$ intersects~$\R B$ in at least one additional
point (or one of the intersection points is multiple). Thus, the number
of intersection points of~$A$ and~$B$, counted with their multiplicity,
is at least $2k^2+k+1>mq$, and, by the Bezout theorem, $B$ contains~$A$
as a component (as an $M$-curve, $A$~is irreducible). Thus, $B$ decomposes
into two curves, $A$ and~$B'$, of degree $m=k+1$ and $q-m=k-2$,
respectively. Furthermore, $\RA$ is contained in a connected component
of~$\R B'$; hence, $A$ and~$B'$ intersect in at least $k(k-1)>(k+1)(k-2)$
points. Since $A$ is irreducible, this contradicts to the Bezout theorem.

If $q$ is even put $q=2k$ and $m=q+4=2(k+2)$ and consider a Wiman
$M$-curve~$A$ of degree $m$ (i.e., a curve with the real scheme
$\langle(k+1)(k+3)+1\rangle\+\12k(k+1)\langle1\langle1\rangle\rangle$; it
can be constructed by a small perturbation of the square of an $M$-curve
of degree $k+2$). Suppose that
$c'(q)\ge\frac12{q(q+3)}-\frac18{q(q-2)}+1=\frac32k(k+1)+2k+1$. Pick one
point in the inner oval of each nest of~$A$ and one point in $k^2+3k+1$
of $k^2+4k+4$ empty outer ovals, totally $\frac32k(k+1)+2k+1$ points, and
trace through them a connected component of a real curve~$B$ of
degree~$q$. Then $A$ and~$B$ intersect in at least $4k^2+8k+2>4k^2+8k$
points, which contradicts to the Bezout theorem.
\endproof

\remark{Remark}
The case of $q$~odd can be treated in an alternative way, similar to that
of $q$~even. Surprisingly, this does not improve the bound.
\endremark

\claim[Corollary]\label{16.1.E}
One has
$14\le c(5)\le c'(5)\le 17$,\enspace
$17\le c(6)\le c'(6)\le 24$, and
$20\le c(7)\le c'(7)\le 29$.
\qed
\endclaim

The upper bound for $c'(6)$ can easily be improved:

\claim[Proposition]\label{16.1.F}
$c'(6)\le 23.$
\endclaim

\proof
Suppose that $c'(6)\ge 24$, pick a (necessarily irreducible) $M$-curve of
degree~9 with the scheme $J\+\langle1\langle
a\rangle\rangle\+\langle1\langle b\rangle\rangle\+\langle1\langle
c\rangle\rangle\+ \langle1\langle d\rangle\rangle$, $a,b,c,d>0$,
$a+b+c+d=24$, and trace a connected component of a real sextic through
all the 24 empty ovals. The two curves intersect in at least $48+8=56>54$
points.
\endproof

\subsection{Obstructions to rigid isotopy }\label{s16.3}
The topological characteristics of a real curve, such as its real scheme,
type, and complex scheme
(i.e., real scheme enhanced with the complex orientations) for curves of
type~I are preserved under rigid isotopies. Another powerful source of
obstructions is the Bezout theorem applied to appropriate auxiliary
curves.

To our knowledge, obstructions derived from the Bezout theorem have never
been studied systematically.  Thus, we confine ourselves to an
explanation of the principal idea and its simplest realizations taken
from \cite{Fie:pencil} and \cite{GM.zamecha}.

Let us start from auxiliary curves of the lowest possible degree: consider
lines. In this case we look for pairs of points, in $\Rp2\sminus\R A$ for
simplicity, selected topologically in such a way that the line~$L$
through them would intersect the curve only at real points. If such a
selection is possible, the position of~$\R A$ with respect to~$\R L$ is
preserved under rigid isotopies of~$A$. In particular, the position
of~$\R A$ is the same with respect to all real lines close to~$L$.

To make this idea more constructive, let us introduce the following
definition. A real line~$L$ with 2 marked points $p_1,p_2\in\R L$ outside
the real part of a real nonsingular curve~$A$, is called a \emph{lock}
of~$A$, if (1) the number of points in $S=\R L\cap\R A$ is equal to the
degree of~$A$, (2) for any pair of points (possibly
coinciding)\mnote{\Dg: added after Zv} in~$S$ bounding a segment in $\R
L\sminus(\{p_1,p_2\}\cup S)$ and an arc in $\R A\sminus S$, the cycle
formed by the segment and the arc either is not homologous to~0 in $\Rp2$
or envelops at least one oval of~$A$. It is clear that a lock extends
continuously to any rigid isotopy of~$A$ and, moreover, regardless of the
extension the topological type of the quadruple $(\Rp2,\R A\cup\R L,\R
L,\{p_1, p_2\})$ does not change. Even this simple consideration gives us
two invariants: existence of a lock and the topological types of the
quadruples corresponding to all locks.

As an example, consider two curves of degree 6 with the scheme
$\langle2\rangle\+\langle1\langle2\rangle\rangle$ constructed from
2~products of ellipses, as shown in Figure~\ref{f16.1}. The first curve
has no lock. For the second one, any pair of points, one inside each of
the two
outer ovals, define a lock.  Hence, the curves are not
rigidly isotopic.

\topinsert
 \line{\hss\vtop{\topglue0pt \hbox{\epsffile{2prod2.pdx}}}\qquad
  \vtop{\topglue0pt \hbox{\epsffile{2prod1.pdx}}}\hss}
\figure\label{f16.1}
 Sextic without and with a lock
\endfigure
\endinsert

There are situations where simultaneous consideration of several locks is
useful. Example: the two curves of degree~5 with 4~ovals
shown in Figure~\ref{f16.2}.
In both the curves any pair of ovals defines a
lock; hence, any triple defines a triple lock. The first curve has a
triple of ovals whose triple lock has the following properties: one of
the connected components of the complement of the three lines contains
the fourth oval and does not intersect the one-sided component. The other
curve has no such triple. Hence, the two curves are not rigidly isotopic.

\topinsert
 \line{\hss\vtop{\topglue0pt \hbox{\epsffile{quintic2.pdx}}}\qquad
  \vtop{\topglue0pt \hbox{\epsffile{quintic1.pdx}}}\hss}
\figure\label{f16.2}
 Two quintic with different triple locks
\endfigure
\endinsert

In fact, in both the examples the curves differ by their complex scheme:
as one can see from the construction, in each example one curve is of
type~$\I{}$ and the other, of type~$\II$. Pairs of not rigidly isotopic
curves with the same complex scheme appear first in degree~$7$. An
example, due to A.~Marin~\cite{GM.zamecha}, is shown in
Figure~\ref{f16.3}. The two curves with the complex scheme
$\<J\+1^+\<1^+\+2^->\+6^+\+5^->$ are constructed by the Hilbert method,
with a difference at the last step. The three inner ovals define a triple
lock, and the two curves differ by the distribution of the outer ovals
among the connected components of the complement of the three lines.

Another example of not rigidly isotopic curves of degree~$7$ with the
same complex scheme was constructed by T.~Fiedler~\cite{Fie:pencil}. (In
his example the curves are not dividing.)

\topinsert
 \line{\hss\vtop{\topglue0pt \hbox{\epsffile{marin1.pdx}}}\qquad
  \vtop{\topglue0pt \hbox{\epsffile{marin2.pdx}}}\hss}
\figure\label{f16.3}
 Two curves of degree seven with different triple locks
\endfigure
\endinsert

There are examples where quadrics are used for locks instead of lines. In
one of them, due to I.~Itenberg, the inner ovals of the
two negative injective pairs of the sextic $\<1\<3^+\+2^->>\+\<5>$ are
separated from each other by a quadric trough the three other inner ovals
and two outer ovals. Certainly, locks of higher degrees can also be
considered. We interrupt this theme because it is not yet sufficiently
understood.

Back to Marin's example, note that it is not known whether the two curves
are equivariantly isotopic. More generally, the following two fundamental
problems are still open: {\proclaimfont Are there Bezout like
obstructions which are not flexible\rom? Do there exist curves which are
equivariantly isotopic but not rigidly isotopic\rom?} (When this survey
was ready, S.~Fiedler-Le Touz\'e and S.~Orevkov announced that they found a
flexible affine sextic not realizable algebraically because of a Bezout
like obstructions in respect to cubics.)

\subsection{Curves on Quadrics}\label{C3.4}
We conclude this section with a brief account of known results on curves
on nonsingular quadrics, which seem to be the only case (except~$\Cp2$)
that has been seriously studied. Let $X$ be such a quadric.  From the
complex point of view $X=\Cp1\times\Cp1$, and up to biholomorphism $X$
admits two anti-holomorphic involutions with nonempty real part; the
resulting real surfaces are the hyperboloid, which is an $M$-surface with
real part homeomorphic to torus, and the ellipsoid, an $(M-1)$-surface
with $\RX\cong S^2$.

Fix a pair $P_1$, $P_2$ of generatrices of~$X$. The fundamental classes
$[P_1]$, $[P_2]$ form a basis of $H_2(X)\cong\Z\oplus\Z$; the intersection
form is given by $[P_1]^2=[P_2]^2=0$, $[P_1]\circ[P_2]=1$.  For any other
curve $A\subset X$ one has $[A]=m_1[P_1]+m_2[P_2]$ for some nonnegative
integers $m_1$, $m_2$. The pair $(m_1,m_2)$ is called the \emph{bidegree}
of~$A$.  If $(x_0:x_1)$, $(y_0:y_1)$ are homogeneous coordinates
in~$P_1$, $P_2$ respectively, then $A$ is given by a
polynomial\mnote{\Dg: edited after Zv}
$$
F(x_0,x_1;y_0,y_1)=\sum_{i,j=1}
^{m_1,m_2}a_{ij}x_1^ix_0^{m_1-i}y_1^jy_0^{m_2-j}
$$
homogeneous of degree~$m_1$ in~$(x_0,x_1)$ and of degree~$m_2$
in~$(y_0,y_1)$. Such curves form a projective space of dimension
$m_1m_2+m_1+m_2$; if $A$ is nonsingular, the genus of~$A$ equals
$(m_1-1)(m_2-1)$ (see~\eqref{C3.4}; $K_X$ has bidegree~$(-2,-2)$\,).

In the case of hyperboloid $\conj$ acts by conjugating all the four
coordinates; hence, $A$ is real if and only if all $a_{ij}$ are real.
In the case of ellipsoid $\conj$ acts via
$[(x_0:x_1),(y_0:y_1)]\mapsto
[(\bar y_0:\bar y_1),(\bar x_0:\bar x_1)]$,
and $A$ is real if and only if $a_{ij}=\bar a_{ji}$ (in particular,
$m_1=m_2$).

From now on we assume that $A$ is a real curve. First, let $(X,\conj)$ be
a hyperboloid. The fundamental classes $[\R P_1]$, $[\R P_2]$, endowed
with some orientations (which are to be fixed) form a basis of
$H_1(\RX)\cong\Z\oplus\Z$. The real part~$\RA$ may have components of two
types: contractible in~$\RX$ and noncontractible. The contractible
components are called \emph{ovals}, their number is denoted by~$l$. Each
oval bounds a disk in~$\RX$, which is called the \emph{interior} of this
oval. All the noncontractible components realize the same nontrivial
class $(c_1,c_2)\in H_1(\RX)$, where $c_1$, $c_2$ must be relatively
prime. The number of such components is denoted by~$h$. In order to
encode the isotopy type of~$\RA\subset\RX$ we will use the notation
$$
\<(c_1,c_2),\,\operator{scheme}_1,(c_1,c_2),\,\operator{scheme}_2,
\dots,(c_1,c_2),\,\operator{scheme}_h>,
$$
where $\,\operator{scheme}_1,\dots,\,\operator{scheme}_h$ are the schemes
of ovals of~$\RA$ in the annulus bounded by two consequent
noncontractible components.

If both $m_1$ and~$m_2$ are even, $\RA$ separates~$\RX$ into two parts
$\RX_\pm$. Obviously, $\chi(\RX_-)=-\chi(\RX_+)$. If all the components
of~$\RA$ are ovals, we denote by $\RX_+$ the part of genus~$1$.

\claim[Prohibitions for curves on a hyperboloid]\label{C3.4.A}
Let $A$ be a real curve of bidegree~$(m_1,m_2)$ on a hyperboloid~$X$, and
$l$, $h$, $c_1$, $c_2$ as above. Then $A$
satisfies the following restrictions\rom:

\medskip
\noindent{\bf The Harnack inequality\rom:}
\roster
\item\local1
$l+h\le(m_1-1)(m_2-1)+1$.
\endroster
\medskip
\noindent{\bf Congruences\rom:} Assume that $m_1$,
$m_2$ are even and $\frac12(m_1c_2+m_2c_1)=c_1c_2\bmod2$.
Then \rom(with the above convention on the choice
of~$\RX_+$\rom)
\roster
\item[2]\local2
If $A$ is an $M$-curve, then $\chi(\RX_+)=\frac12m_1m_2\bmod8$\rom;
\item\local3
If $A$ is an $(M-1)$-curve, then
$\chi(\RX_+)=\frac12m_1m_2\pm1\bmod8$\rom;
\item\local4
If $A$ is an $(M-2)$-curve and $\chi(\RX_+)=\frac12m_1m_2+4\bmod8$, then
$A$ is of type~$\I{}$\rom;
\item\local5
If $A$ is of type~$\I{}$, then $\chi(\RX_+)=0\bmod4$.
\endroster
\medskip
\noindent{\bf Consequences of the Bezout Theorem\rom:}
\roster
\item[6]\local6
$hc_i=m_1\bmod2$ and $h\ls|c_i|\le m_1$ for $i=1,2$\rom;
\item\local7
if $\RA$ has $r$ pairwise disjoint nests of depths $s_1,\dots,s_r$,
then
$$
\sum_{j=1}^rs_j\le\frac12\bigl(m_1-h\ls|c_1|+\intpart{r/2}(m_2+h\ls|c_2|)\bigr)
$$
\rom(and a similar inequality with $m_1,c_1$ interchanged with
$m_2,c_2$\rom).
\endroster
\endclaim

\proof
\loccit1 is the classical Harnack-Klein bound.
\loccit2--\loccit5 follow from~\ref{C4.Mikh.B}
and~\ref{C4.Mikh.D}.
Statement~\loccit6 is the real version of the Bezout theorem
(see~\ref{C3.3}) applied to~$A$ and curves of bidegree~$(1,0)$
and~$(0,1)$, and~\loccit7 is proved by applying the Bezout theorem to~$A$
and the curve of bidegree~$(\intpart{r/2},1)$ traced through $r$~points
inside the innermost ovals of the nests (see~\cite{Zvonilov}).
\endproof

\remark{Remark}
In a weaker form the above congruences appeared in \cite{Mat:91;Mat:92}.
This inspired
\cite{Mikhalkin}, where
the congruences were improved and generalized
to take the form of~\ref{C4.Mikh.B}.
\endremark

If $(X,\conj)$ is an ellipsoid,
all the components of~$\RA$ are contractible and are called ovals; their
number is denoted by~$l$.
Fix an \emph{exterior point} $\infty\in\RX\sminus\RA$
and define the \emph{interior} of an oval $C\subset\RX$
as the component of $\RX\sminus C$ which does not contain~$\infty$. Then a
natural partial order arises, which can be encoded as usual.

Since $\RX\cong S^2$, any curve~$\RA$ divides~$\RX$ into two parts,
denoted by~$\RX_\pm$.

\claim[Prohibitions for curves on an ellipsoid]\label{C3.4.B}
Let $A$ be a real curve of bidegree~$(m,m)$ on an ellipsoid~$X$,
and $l$ the number of components of~$\RA$. Then $A$
satisfies the following restrictions\rom:

\medskip
\noindent{\bf The Harnack inequality\rom:}
\roster
\item\local1
$l\le(m_1-1)(m_2-1)+1$.
\endroster
\medskip
\noindent{\bf Congruneces\rom:}
Let $m$ be odd.
\roster
\item[2]
\local2 If $A$ is an $M$-curve, then
$\chi(\RX_\pm)=\frac12(m^2+1)\bmod8$\rom;
\item\local3
If $A$ is an $(M-1)$-curve, then
$\chi(\RX_\pm)=\frac12(m^2+1)\pm1\bmod8$\rom;
\item\local4
If $A$ is an $(M-2)$-curve and $\chi(\RX_\pm)=\frac12(m^2+1)+4\bmod8$,
then $A$ is of type~$\I{}$\rom;
\item\local5
If $A$ is of type~$\I{}$, then $\chi(\RX_\pm)=1\bmod4$.
\endroster
\medskip
\noindent{\bf Consequences of the Bezout Theorem\rom:}
\roster
\item[6]\local6
The total number of ovals in any three pairwise disjoint nests does
not exceed~$m$.
\endroster
\endclaim

\proof\nofrills\enspace
is similar to that of~\ref{C3.4.A}.
\endproof



Next two well known statements are found, e.g., in~\cite{Zv4}. Their
standard proofs, reproduced below, serve as a model for a number of
similar cases (such as rational and hyperelliptic curves on ruled
surfaces, intersections of quadrics, or homotopy classification of
nondegenerate homogeneous polynomial vector fields).\mnote{\Dg: changed
after Zv}

\claim[Curves of Bidegree~$(m,1)$]\label{C3.4.C}
A curve of bidegree~$(m,1)$ on a hyperboloid is determined up to rigid
isotopy by its real scheme, which is $\<(w,1)>$ with $-m\le w\le m$ and
$w=m\bmod2$. All such curves are of type~$\I{}$.

There are two rigid isotopy classes of curves of bidegree~$(1,1)$ on an
ellipsoid\rom: $\<0>$ and $\<1>$, the latter being of type~$\I{}$.
\endclaim

\proof
A nonsingular curve of bidegree~$(m,1)$ on a hyperboloid has equation
of the form
$$
y_0a(x_0,x_1)+y_1b(x_0,x_1)=0,
$$
where $a$ and~$b$ are nonzero homogeneous polynomials of degree~$m$
without common roots. The roots of~$a$ and~$b$ being fixed, the
polynomials are recovered uniquely up to constant factors. Any pair of
consecutive real root of one of these polynomials not separated by real
roots of the other one can be removed to the complex domain. After all
such pairs are removed, $\ls|w|$ is the number of real roots of~$a$
(or~$b$); for each nonzero\mnote{\Dg: edited after W} value of~$\ls|w|$
there are two classes of curves which differ by the signs of~$a$ and~$b$
at a fixed generic point.

Any curve of bidegree~$(1,1)$ on an ellipsoid is cut by a plane, and
it is obvious that there are only two possibilities.
\endproof

\claim[Curves of bidegree~$(m,2)$]\label{C3.4.D}
A curve of bidegree~$(m,2)$ on a hyperboloid is determined up to
rigid isotopy by its real scheme, which is either
$$
\<(0,1),\Ga_1,(0,1),\Ga_2,\dots,(0,1),\Ga_h>
$$
with $0\le h\le m$, $h=m\bmod2$, and $h+\sum\Ga_i\le m$, or
$$
\<2(\12{w},1)> \text{ for $m$ even or } \<(w,2)> \text{ for $m$ odd}
$$
with $-m\le w\le m$ and $w=m\bmod2$. The latter scheme and the $M$-schemes
are of type~$\I{}$\rm; the others are of type~$\II$.

There are three rigid isotopy classes of curves of bidegree~$(2,2)$ on an
ellipsoid\rom: $\varnothing$, $\<1>$, and $\<2>$. The only scheme of
type~$\I{}$ is~$\<2>$.
\endclaim

The proof given below works, more or less without changes, whenever the
curves under consideration are hyperelliptic (so that the hyperelliptic
structure is induced from a ruling of the surface). For this reason we
reproduce all the essential details.

\proof
A curve of bidegree~$(m,2)$ on a
hyperboloid has equation of the form
$$
y_0^2a(x_0,x_1)+2y_0y_1b(x_0,x_1)+y_1^2c(x_0,x_1)=0,
$$
where $a$, $b$, and~$c$ are homogeneous polynomials of degree~$m$. Such a
curve is nonsingular if and only if the discriminant $D=b^2-ac$ has no
multiple roots.

Let $V=\{(x_0:x_1)\in\RP^1\,|\,D(x_0,x_1)\ge0\}$. Consider an affine part
$x_1\ne0$. The connected components of the set consisting of the
nonsingular curves with fixed $V$ and fixed roots $u_1,\dots,u_m$ of~$a$
(which may be multiple) are determined by the choice of the values
$b(u_i)=\pm\sqrt{D(u_i)}$, $i=1,\dots,m$, i.e., by the choice of the sign
of~$b$ at the real roots of~$a$ (which belong to~$V$) and one of the two
branches of the square root at each pair of complex conjugate roots. The
components can be joined using the following moves:
\roster
\item
a pair of consecutive real roots of~$a$ at which $b$ has the same
sign can be removed to the complex domain; a pair of complex
conjugate roots can be moved to any component of~$V$;
\item
the sign of~$b$ at a real root of~$a$ not separated from~$\partial V$
by other roots can be changed (by moving this root to~$\partial V$);
\item
the branch of $\sqrt{D}$ at any pair of complex conjugate roots
of~$a$ can be changed (by moving this pair to a pair of roots of~$D$).
\endroster
Thus, if $V=\varnothing$ or $\partial V\ne\varnothing$, the sign of~$b$
at each (real or imaginary) root of~$a$ can be reversed and the curve is
determined up to rigid isotopy by the components of~$V$ and the parity of
the number of real roots of~$a$ in each component:  the corresponding
component of the curve is an oval if this number is even or of
class~$(1,0)$ if it is odd.  If $V=\RP^1$, then it contains $\ls|w|$~real
roots of~$a$ so that $b$~has opposite signs at each pair of consecutive
roots. This gives the second group of real schemes in~\ref{C3.4.D}.

Consider now the case of ellipsoid. All empty curves are rigidly isotopic
(in fact, in any bidegree), as they form a convex set. If the curve is
nonempty, we blow up a real point of this curve and blow down the proper
inverse images of the two imaginary generatrices through this point. The
result is a nonsingular cubic in~$\RP^2$, and the inverse transformation
is given by a pair of complex conjugate point of this cubic. The two
rigid isotopy classes of cubics correspond to the two classes of curves
on the ellipsoid.
\endproof

\claim[Curves of bidegree~$(3,3)$]\label{C3.4.E}
A curve of bidegree~$(3,3)$ on a hyperboloid or an ellipsoid is
determined up to rigid isotopy by its real scheme and, in the case
of hyperboloid, its type.  On a hyperboloid the real schemes are
$$
\<(1,\pm1),\Ga>\text{ with $\Ga\le4$,}\quad
\<3(1,\pm1)>,\quad\<(3,\pm1)>,\quad\<(1,\pm3)>.
$$
The schemes $\<(1,\pm1),4>$, $\<3(1,\pm1)>$, $\<(3,\pm1)>$, and
$\<(1,\pm3)>$ are of type~$\I{}$\rom; the scheme $\<(1,\pm1),2>$ is
of indefinite type\rom; the other schemes are of type~$\II$.  On an
ellipsoid the real schemes are
$$
\<\Ga>\text{ with $\Ga\le5$ and }\<1{\<1{\<1>}>}>.
$$
The schemes $\<5>$ and $\<1{\<1{\<1>}>}>$ are of type~$\I{}$\rom; the
others are of type~$\II$.
\qed
\endclaim

Similar to~\ref{C2.3.A}, the proof of~\ref{C3.4.E} is based on enumerating
the walls in the space of curves; the latter problem is mainly reduced to
the classification of plane quartics with a fixed tangent. Details can be
found in~\cite{DZ}.

The remaining results of this section
(Propositions~\ref{C3.4.G}--\ref{C3.4.I}) only list the possible schemes
for certain bidegrees. All the schemes not on the lists are prohibited
by~\ref{C3.4.A} and~\ref{C3.4.B}. The listed schemes have been realized
by D.~Gudkov, A.~Usachev, and E.~Shustin (see~\cite{Gudkov},
\cite{Gudkov-Usachev}, and~\cite{Gudkov-Shustin}). The assertions on the
type of the schemes follow mainly from Rokhlin's formula of complex
orientations~\ref{C3.2.A}. Note also that the complex orientations are
those and only those which do not contradict to Rokhlin's and Orevkov's
formulas. Details and realization of the schemes of indefinite type can
be found in Zvonilov~\cite{Zvonilov}.

\claim[Schemes of bidegree~$(m,3)$]\label{C3.4.G}
A curve of bidegree~$(m,3)$ with $m>3$ on a hyperboloid can only have one
of the following real schemes\rom:
\roster
\item\local1
$\<3(w,1)>$ with $3\ls|w|\le m$ and $w=m\bmod2$\rom;
\item\local2
$\<(w,3)>$ with $\ls|w|\le m$, $w=m\bmod2$, and $w\ne0\bmod3$\rom;
\item\local3
$\<(m,\pm1)>$\rom;
\item\local4
$\<(w,1),\Ga>$ with $\ls|w|\le m-2$, $w=m\bmod2$, and $0\le\Ga\le2m-2$.
\endroster
Schemes~\loccit1---\loccit3 are of type~$\I{}$\rom; schemes~\loccit4 are
of type~$\I{}$ when $\Ga=2m-2$ and are of type~$\II$ when $\Ga$ is odd or
zero.\mnote{\Dg: edited after Zv}
\qed
\endclaim

All the schemes above (including complex orientations) were realized by
A.~Korchagin~\cite{Korchagin}, who used Viro's patching method. According
to V.~Zvonilov, they can as well be realized by $T$-curves.\mnote{\Dg:
added after Zv}
A new promising approach
to the study of trigonal curves on ruled surfaces by means of
Grothendieck's \emph{dessins d'enfant}
applied to the discriminant of the curve
was recently suggested by S.~Orevkov~\cite{Or:RET}.

\claim[Schemes of bidegree~$(4,4)$ on a hyperboloid]\label{C3.4.H}
The following is a complete list of real schemes of bidegree~$(4,4)$ on a
hyperboloid\rom:
\roster
\item\local1
$\<4(1,0)>$, $\<4(0,1)>$, and $\<4(1,\pm1)>$\rom;
\item\local2
$\<2(2,\pm1)>$ and $\<2(1,\pm2)>$\rom;
\item\local3
$\<(1,\pm1),\Ga,(1,\pm1),\Gb>$ with $0\le\Ga+\Gb\le8$ and
$\Ga\ge\Gb$\rom;
\item\local4
$\<(1,0),\Ga,(1,0),\Gb>$ and $\<(0,1),\Ga,(0,1),\Gb>$ with either
$0\le\Ga+\Gb\le6$ and $\Ga\ge\Gb$, or $(\Ga,\Gb)=(8,0)$, $(4,4)$,
$(7,0)$, or~$(4,3)$\rom;
\item\local5
$\<\Ga\dplus1{\<\Gb>}>$ with either $0\le\Ga+\Gb\le7$, or
$(\Ga,\Gb)=(0,9)$, $(4,5)$, $(8,1)$, $(0,8)$, $(3,5)$, $(4,4)$,
$(7,1)$, or~$(8,0)$\rom;
\item\local6
$\<2{\<1>}>$\rom;
\item\local7
$\<0>$.
\endroster
Schemes~\loccit1, \loccit2, and~\loccit6 are of type~$\I{}$.
Scheme~\loccit7 is of type~$\II$. Schemes~\loccit3 are of type~$\I{}$ if
$\Ga+\Gb=8$, are of type~$\II$ if $\Ga+\Gb$ is odd or $(\Ga,\Gb)=(0,0)$,
$(1,1)$, or~$(4,0)$, and are of indefinite type otherwise.
Schemes~\loccit4 are of type~$\I{}$ if $(\Ga,\Gb)=(8,0)$, $(4,4)$,
or~$(5,1)$, are of indefinite type if $(\Ga,\Gb)=(3,3)$, $(4,0)$,
$(2,2)$, or~$(1,1)$, and are of type~$\II$ otherwise. Schemes~\loccit5
are of type~$\I{}$ if $(\Ga,\Gb)=(0,9)$, $(4,5)$, $(8,1)$, $(1,6)$,
or~$(5,2)$, are of indefinite type if $(\Ga,\Gb)=(3,4)$, $(7,0)$,
$(2,3)$, or~$(4,1)$, and are of type~$\II$ otherwise.
\qed
\endclaim

\claim[Schemes of bidegree~$(4,4)$ on an ellipsoid]\label{C3.4.I}
The following is a complete list of real schemes of bidegree~$(4,4)$ on
an ellipsoid\rom:
\roster
\item\local1
$\<\Ga\dplus1{\<\Gb>}>$ with $0\le\Ga+\Gb\le9$ and $\Ga\ge\Gb$\rom;
\item\local2
$\<1{\<1{\<1{\<1>}>}>}>$\rom;
\item\local3
$\<0>$.
\endroster
Schemes~\loccit2 and~\loccit1 with $\Ga+\Gb=9$ are of type~$\I{}$\rom;
schemes~\loccit3 and~\loccit1 with $\Ga+\Gb\le4$ or even or
$(\Ga,\Gb)=(5,0)$ are of type~$\II$\rom; the other schemes are of
indefinite type.
\qed
\endclaim

The double covering of a quadric ramified in a curve of bidegree $(4,4)$
is a $K3$-surface whose Picard group contains $U(2)$ (i.e., a pair $a$,
$b$ with $a^2=b^2=0$ and $ab=2$), and, vice versa, any $K3$-surface~$X$
with a distinguished subgroup $U(2)\subset\operatorname{Pic}X$ is
obtained in this way (see, e.g.,~\cite{Nikulin:forms}). Thus, the curves
of bi-degree~$(4,4)$ can be studied similar to plane sextics. In
particular, this techniques can provide their rigid isotopy
classification.
First steps in this direction are made
in~\cite{Matsuoka}.
In the case of hyperboloid it will not be a
surprise if the list of~\ref{C3.4.H} turns out to coincide with the
classification up to rigid isotopies and automorphisms of the
surface. (In the case of ellipsoid there probably are additional
invariants.) Note that some curves on the hyperboloid, e.g.,
$\<1\<9>>$, $\<4\dplus 1\<5>>$, $\<1\<8>>$, $\<3\dplus 1\<5>>$, and
$\<4\dplus 1\<4>>$, are not chiral, see~\cite{Kh:rigid}.

As to curves on surfaces of higher degree
little is known.
We can only cite~\cite{Segre} and~\cite{Mikh:cubic}, which treat,
respectively, real lines and real quadric sections on real cubic surfaces.

\end



\ifx\rokhloaded\undefined\input rokh.def \fi

\let\dd\partial
\def\aug{\operator{aug}}
\def\Sm{\operator{Sm}}

\let\GD\Delta

\def\barX{X'}
\def\barA{A'}
\def\scup{\mathbin{\scriptstyle\cup}}
\def\scap{\mathbin{\scriptstyle\cap}}


\begin

\appendix{Topology of involutions}\label{apA}

\subsection{Smith theory}\label{Smith.theory}\label{s7}
Most results cited in this section are due to P.~A.~Smith; proofs can be
found, e.g., in~\cite{Bredon, Chapter~3}. Throughout the section we
consider a topological space~$X$ with involution $c\:X\to X$ and denote
$F=\Fix c$ and $\barX=X/c$.

\theorem\label{1.1.1}\label{Smith.seq}
There are two natural, in respect to equivariant maps, exact sequences,
called \rom(homology and cohomology\rom) \emph{Smith sequences}
of~$(X,c)$\rom:
$$
\gather
@>>>H_{p+1}(\barX,F)@>\,\,\GD>>H_p(\barX,F)\oplus H_p(F)
  @>\,\,{\tr^*}+{\inc_*}>>H_p(X)
  @>\,\,\pr_*>>H_p(\barX,F)@>>>\rlap{\,,}\\
@>>>H^p(\barX,F)@>\,\,\pr^*>>H^p(X)@>\,\,{\tr_*}\oplus{\inc^*}>>
  H^p(\barX,F)\oplus H^p(F)@>\,\,\GD>>H^{p+1}(\barX,F)@>>>\rlap{\,.}
\endgather
$$
The homology and cohomology connecting homomorphisms~$\GD$ are given by
$$
x\mapsto x\scap\Go\oplus\dd x\quad\text{and}\quad x\oplus f\mapsto
x\scup\Go+\Gd f,
$$
respectively, where $\Go\in H^1(\barX\sminus F)$ is the characteristic
class of the double covering $X\sminus F\to\barX\sminus F$. The images
of~${\tr^*}+\inc_*$ and~$\pr^*$ consist of invariant classes\rom:
$\Im\tr^*\subset\Ker(1+c_*)$ and $\Im\pr^*\subset\Ker(1+c^*)$.
\pni
\endtheorem

\subsubsection{}\label{1.1.1'}
The \emph{transfer} homomorphisms $\tr^*$ and~$\tr_*$ in~\ref{1.1.1} can
be described as follows. Represent~$X$ as the $S^0$-bundle associated
with a relative $D^1$-bundle $T\to\barX$ whose characteristic class
is~$\Go$. Namely, take the cylinder of the projection $X\to\barX$ and for
each $f\in F$ contract the fiber over $f$.  Then $\tr^*$ and $\tr_*$ are
the compositions
$$
\spreadlines{-\jot}
\align
\tr^*\:\quad&
  H_p(\barX,F)@>\,\,\th^{-1}>>H_{p+1}(T,X)@>\,\,\dd>>H_p(X),\\
\tr_*\:\quad&
  H^p(X)@>\,\,\Gd>>H^{p+1}(T,X)@>\,\,\th^{-1}>>H^p(\barX,F),
\endalign
$$
where $\th$ stands for the Thom isomorphisms.

The following immediate consequences of~\ref{1.1.1}, which we state in
the homology setting, has an obvious counterpart for cohomology.

\corollary[Corollary \rm(of~\ref{1.1.1})]\label{Sinqty}\label{1.1.2}
\roster
\item\local1
\strut$\dim H_*(F)+2\sum_{p}\dim\Coker({\tr^p}+\inc_p)=\dim H_*(X)$
\rom(Smith identity\rom)\rom;
\item\local2
\strut$\dim H_*(F)\le\dim H_*(X)$ \rom(Smith inequality\rom)\rom;
\item\local3
\strut$\dim H_*(F)\le\dim H^1(\ZZ;H_*(X))$ \rom(Borel-Swan
inequality\rom)\rom;
\item\local4
\strut$\dim H_*(F)=\dim H_*(X)\bmod2$\rom;
\item\local5
\strut$\chi(X)=\chi(F)+2\chi(\barX,F)$\rom;
\item\local6
\strut$\chi(X)=2\chi(\barX)-\chi(F)$
\endroster
\rom(Recall that $H^1(\ZZ;H_*(X))=\Ker(1+c_*)/\Im(1+c_*)$.\rom)
\pni
\endcorollary

\subsubsection{}\label{1.1.8}
If \iref{1.1.2}2 turns into an equality, i.e., $\dim H_*(F)=\dim H_*(X)$,
one says that $c$ is an \emph{$M$-involution}. If \iref{1.1.2}3 turns
into an equality (which is equivalent to
$\Im({\tr^*}+\inc_*)\supset\Ker(1+c_*)$), $c$ is called
\emph{\rom($\ZZ$-\rom)Galois maximal}. (This terminology is introduced by
V.~A.~Krasnov \cite{Kr:GM}. R.~Thom~\cite{ThomStras} calls a dimension
$p\in\N$ \emph{regular} for $(X,c)$ if
$\Im({\tr^p}+\inc_p)\supset\Ker(1+c_p)$.)

In general, due to~\iref{1.1.2}2 and~\ditto4, one has $\dim H_*(F)=\dim
H_*(X)-2d$ for some integer $d\ge0$; in this case $c$ is called an
\emph{$(M-d)$-involution}.

\subsubsection{Geometrical construction of the Smith sequences}\label{1.1.9}
Introduce the \emph{Smith chain complexes}
$$
\align
\Sm_*(X)&=\Ker[(1+c_*)\:S_*(X)\to S_*(X)],\\
\Sm_*(X,F)&=\Ker[(1+c_*)\:S_*(X,F)\to S_*(X,F)].
\endalign
$$
and \emph{Smith homology} $H_r(\Sm_*(X))$ and $H_r(\Sm_*(X,F))$. There is
a canonical isomorphism $\Sm_*(X,F)=\Im[(1+c_*)\:S_*(X)\to S_*(X)]$ and a
canonical splitting $\Sm_*(X)=S_*(F)\oplus\Im(1+c_*)$. Furthermore,
$\tr^*\:S_*(\barX,F)\to\Sm_*(X,F)$ is an isomorphism, and in view of the
above identifications the Smith sequences are the long homology and
cohomology exact sequences associated with the short exact sequence of
complexes
$$
0@>>>\Sm_*(X)@>\,\,\operator{inclusion}>>S_*(X)
  @>\,\,1+c_*>>\Sm_*(X,F)@>>>0.
$$

\theorem[Steenrod squares in the Smith sequences]\label{1.1.10'}
For any $x\in H_p(\barX,F)$ and $y\in H^p(X)$ one has
$\Sq\tr^*x=\tr^*\Sq(x\scap(1+\Go))$ and $\tr_*\Sq
y=(\Sq\tr_*y)\scup(1+\Go)$.
\endtheorem

\proof
The statement follows from~\ref{1.1.1'} and the Thom-Cartan formula for
sphere bundles (see~\cite{ThomAEN}):  since the total characteristic
class of the relative $D^1$-bundle $T\to\barX$ is $1+\Go$, one has
$\th^{-1}\Sq(x\scap(1+\Go))=\Sq\th^{-1}x$ and $\th^{-1}\Sq
y'=(\Sq\th^{-1}y')\scup(1+\Go)$, where $y'=\Gd y$.
\endproof

\corollary\label{1.1.10}
For any $x\in H_p(\barX,F)$ and $y\in H^p(X)$ one has
$\Sq_1\tr^*x=\tr^*\Sq_1x+\inc_*\dd x$ and
$\tr_*\Sq^1y=\Sq_1\tr_*y+\Gd\inc^*y$.
\endcorollary

\proof
Since the Smith sequences are exact, $\tr^*(x\scap\Go)=\inc_*\dd x$ and
$(\tr_*y)\scup\Go=\Gd\inc^*y$. Hence, the statement follows
from~\ref{1.1.10'}.
\endproof

\subsubsection{Relative version of the Smith sequences}\label{1.1.11}
If $A\subset X$ is a $G$-invariant subspace, there is a natural exact
sequence
$$
\align
@>>>H_{p+1}(\barX,F\cup\barA)&@>\,\,\GD>>
  H_p(\barX,F\cup\barA)\oplus H_p(F,F\cap\barA)
  @>\,\,{\tr^*}+{\inc_*}>>\\
&\qquad\qquad@>\,\,{\tr^*}+{\inc_*}>>H_p(X,A)
  @>\,\,\pr_*>>H_p(\barX,F\cup\barA)@>>>
\endalign
$$
(and a similar sequence in cohomology), where $\barA=A/c$. Its properties
and consequences are similar to those of the absolute Smith sequence.  In
particular, the Smith inequality~\iref{1.1.2}2 turns into
$$
\dim H_*(F,F\cap A)\le\dim H_*(X,A).\eqtag\label{1.1.12}
$$

\subsection{Kalinin's spectral sequence}\label{C1.7}
In this section we describe a relatively new powerful tool, so called
\emph{Kalinin's spectral sequence}. Originally it was introduced by
I.~Kalinin~\cite{Kalinin} as a stabilized version of the Borel-Serre
spectral sequence of the fibration $X\times_{\ZZ}S^\infty\to\Rp\infty$,
where $X$~is a topological space with involution~$c$ and $S^\infty$ is
considered with the standard antipodal involution. The sequence starts at
the homology $H_*(X)$ and converges to the total homology $H_*(\Fix c)$.
The resulting filtration $\CF^*$ on $H_*(\Fix c)$ and the isomorphisms
$\vr_p$ between the limit term of the spectral sequence and
$\Grad_{\CF}H_*(\Fix c)$ were discovered by O.~Viro geometrically before
Kalinin's work and were primarily related to the Smith exact sequence.
The relation between the two construction is established
in~\cite{Degt:Stiefel}. As is shown in~\cite{DegtKh:2} (preprint version)
or~\cite{DIK}, Kalinin's spectral sequence can be derived from the Smith
exact sequence as well (see~\ref{C1.7.D}).

A construction similar to Kalinin's spectral sequence was studied by
V.~A.~Krasnov in a series of papers; see, e.g.,~\cite{Kr:94}.

Below we give an alternative, geometrical, description of the homology
version of Kalinin's spectral sequence and Viro homomorphisms~$\vr_*$ and
state their main properties. (As usual, the cohomology version is obtained
by raising all indices and reversing all arrows; we mention it when
speaking about multiplications and Poincar\'e duality.) Proofs of these
results and/or further references can be found in~\cite{DegtKh:2} and
\cite{DIK}. The differentials~$\rd_p$ of the spectral sequence and Viro
homomorphisms~$\vr_p$ are often regarded as additive relations (i.e.,
partial multivalued homomorphisms) $H_p(X)\dasharrow H_{p+r-1}(X)$ and
$\bv_p\: H_*(\Fix c)\dasharrow H_p(X)$, respectively. As our approach is
geometrical, we have to use chains; depending on the nature of~$X$ one
may regard them as singular, simplicial, smooth, etc. To assure
convergence of the spectral sequence, $X$ must satisfy certain
conditions, which depend on the homology theory chosen (e.g., sheaf
theories and locally compact finite dimensional spaces). Certainly, these
conditions are always fulfilled for compact smooth manifolds and finite
$CW$-complexes.

Thus, let us fix a good (see above) topological space~$X$ with
involution~$c$. Consider the partial homomorphisms $\bv_p\:H_*(\Fix
c)\dasharrow H_p(X)$ and the $\Z$-graded spectral sequence
$(\rH_*(X),\rd_*)$ defined as follows:
\roster
\item
$\bv_0$ is zero on $H_{\ge1}(\Fix c)$ and its restriction to $H_0(\Fix
c)$ coincides with the inclusion homomorphism\rom;
\item
$\vr_p$ is defined on a\/ \rom(nonhomogeneous\/\rom) element~$x\in
H_*(\Fix c)$ represented by a cycle $\sum x_i$ (where $x_i$ is the
$i$-dimensional component of~$x$) if and only if there exist some chains
$y_i$ in~$X$, $1\le i\le p$, so that $\partial y_1=x_0$ and $\partial
y_{i+1}=x_i+(1+c_*)y_i$ for $i\ge1$.  In this case $\vr_p x$ is
represented by the class of $x_p+(1+c_*)y_p$ in $H_p(X)$\rom;
\item
$\RH1_*=H_*(X)$ and $\Rd1_*=1+c_*$\rom;
\item
$\rd_p$ is defined on a cycle $x_p$ in~$X$ if and only if there are chains
$y_p=x_p$, $y_{p+1}$,\dots, $y_{p+r-1}$ so that $\partial
y_{i+1}=(1+c_*)y_i$. In this case $\rd_px_p=(1+c_*)y_{p+r-1}$.
\endroster

\claim[Spectral sequence of involution]\label{C1.7.A}
The  homomorphisms~$\bv_*$ and spectral sequence $(\rH_*,\rd_*)$ are
natural with respect to equivariant maps. Furthermore, $\rH_*$ and $\rd_*$
do form a spectral sequence \rom(i.e., $\rd_p$ are well defined
homomorphisms $\rH_p\to\rH_{p+r-1}$ and\/
$\RH{r+1}_p=\Ker\rd_p/\Im\rd_{p-r+1}$\rom), and this sequence converges to
$H_*(\Fix c)$ via~$\bv_*$, i.e., $\bv_p$ induces an \rom(honest\rom)
isomorphism $\CF^p/\CF^{p+1}\to\iH_p$, where
$\CF^p=\operatorname{Domain}\vr_p=\Ker\vr_{p-1}$.
\qed
\endclaim

\remark{Remark}
A detailed analysis of the construction above shows that $\bv_p$ can as
well be regarded as a true homomorphism
$$
\bv_p\:\Ker[\bv_{p-1}\:H_*(\Fix c)\dasharrow\RH{r-1}_{p-1}]\to\rH_p,
 \quad r>1
$$
(where $\bv_p\:H_*(\Fix c)\dasharrow\RH1_p=H_p(X)$ is defined as the
inclusion homomorphism on $H_p(\Fix c)$ and zero on $H_{>p}(\Fix c)$\,).
These homomorphisms are, in fact, related to the original, nonstabilized,
Borel-Serre spectral sequence (see~\cite{Degt:Stiefel} for details).
\endremark

The cohomology version $\rH^*\Rightarrow H^*(\Fix c)$ is constructed
similarly and is dual to the homology one.  The cup- and cap-products
in~$X$ naturally descend to, respectively, a $\ZZ$-algebra structure
in~$\rH^*$ and $\rH^*$-module structure in~$\rH_*$. If $X$ is a connected
closed $N$-manifold and $\Fix c\ne\varnothing$, the fundamental class
$[X]$ survives to~$\iH_N$ and the multiplication
$\capprod[X]\:\rH^p\to\rH_{N-p}$ is an isomorphism (Poincar\'e duality).
This, in the usual way, defines a pairing
$\circ\:\rH_p\otimes\rH_q\to\rH_{p+q-N}$, called \emph{Kalinin's
intersection pairing}.

\claim[Kalinin's intersection pairing]\label{C1.7.B}
Let $X$ be a smooth closed connected $N$-man\-i\-fold with a smooth
involution $c\:X\to X$.  Denote by~$w(\nu)$ the total Stiefel-Whitney
class of the normal bundle~$\nu$ of~$\Fix c$ in~$X$.  Then for
$a\in\CF^p$ and $b\in\CF^q$ one has $w(\nu)\capprod(a\circ b)\in
\CF^{p+q-N}$ and
$$
\bv_p a\circ\bv_q b=\bv_{p+q-N}[w(\nu)\capprod(a\circ b)].\rlap\qed
$$
\endclaim

\remark{Remark}
The above identity contains, as a special case, the
Borel-Haefliger~\cite{BH} version of the Bezout theorem over the reals
(with the $\ZZ$-valued intersection): given two real algebraic cycles
$A$, $B$ of complementary dimension in a real algebraic manifold~$X$,
their intersection number in~$X$ is congruent $\bmod\,2$ to the
intersection number of~$\RA$ and~$\R B$ in~$\RX$. To see this, it is
sufficient to notice that $[\RA]$ is the highest degree term of
$\bv^{-1}[A]$.
\endremark

The homology Steenrod operations $\Sq_t\: H_*(X)\to H_{*-t}(X)$ also
descend to $\rH_*$. In order to describe their relation to the ordinary
Steenrod operations in $H_*(\Fix c)$, introduce the \emph{weighed Steenrod
operations}
$$
\hSq_t\: H_p(X)\to H_{\le p}(X),\qquad
 x\mapsto\sum_{0\le j\le t}{\tsize\binom{P-p}{t-j}}\Sq_jx,
$$
where $P>p+t$ is a power of~2. (The binomial coefficients do not depend
on~$P$, see, e.g., \cite{SteenrodEpstein, Lemma~I.2.6}.) Then one has:

\claim[Steenrod squares in Kalinin's spectral sequence]\label{C1.7.C}
If $x\in\CF^p$ and $t\ge0$, then \smash{$\hSq_tx\in\CF^{p-t}$} and
$$
\Sq_t\bv_px=\bv_{p-t}\hSq_tx.\rlap\qed
$$
\endclaim

Kalinin's spectral sequence can alternatively be derived from the Smith
exact sequence~\ref{Smith.seq}.
Below we give the corresponding description of the differentials~$\rd_*$
and Viro homomorphisms~$\vr_*$; the construction of the spectral sequence
and filtration~$\CF^*$ can be found in~\cite{DIK}.)

\claim[Kalinin's spectral sequence and Smith exact sequence]\label{C1.7.D}
The differentials~$\rd_p$ and Viro homomorphisms~$\vr_p$ are given by
$$
\gather
 \rd_p=\tr^{p+r-1}\circ\GD_{p+r-1}^{-1}\circ\ldots\circ
  \GD_{p+1}^{-1}\circ\pr_p\:H_p(X)\dasharrow H_{p+r-1}(X),\\
\vspace{1\jot}
 \bv_p({\textstyle\sum_{r\le p}x_i})=\inc_p x_p+\tr^p y'_p
  \:H_*(\Fix c)\dasharrow H_p(X),
\endgather
$$
where $\inc_*\: H_*(\Fix c)\to H_*(X)$ is the inclusion
homomorphism, $\tr^*$, $\pr_*$, and $\GD_*$ are as in
the Smith exact sequence~\ref{Smith.seq}, and, in the latter equation,
$x_r\in H_r(\Fix c)$ and $y'_p\in H_p(X',\Fix c)$ is defined recursively
via $y'_0=0$, \ $y'_{r+1}=\GD_{r+1}^{-1}(x_r\oplus y'_r)\in
H_{r+1}(X',\Fix c)$.
\qed
\endclaim

\subsubsection{Lifting involutions}\label{apdx.lifts}
Let $(X,c)$ be a topological space with involution, $Y\to X$ a double
covering, and $\Go\in H^1(X)$ its characteristic class. Then $c$ lifts to
an involution on~$Y$ if and only if $\Go$ survives to $\iH^1(X)$. If this
is the case, there obviously are two lifts $c_1,c_2\:Y\to Y$, which
commute with each other and with the deck translation~$\tau$ of the
covering and satisfy the identity $c_1\circ c_2=\tau$.

An obvious necessary condition for $\Go$ to survive to $\iH^1$ is that it
must be $c^*$-invariant (i.e., survive to $\RH2^1$). Note that if $X$ is
connected and $\Fix c\ne\varnothing$, this condition is also sufficient
as all the differentials landing in $\rH^0$ are trivial.

Everything above applies to branched coverings as well, with~$X$ replaced
with $X\sminus A$, where $A$ is the branch set.

\subsubsection{Complex surfaces}\label{Kalinin.surface}
Kalinin's spectral sequence, as well as the related intersection pairing
and Steenrod operations have a transparent geometrical meaning in the
case when $X$ is a complex surface and $c=\conj$ a real structure. Let
$C_i$, $i=1,\dots,k$, be the components of~$\RX$, $\<C_i>\in H_0(\RX)$
the corresponding classes, and $C_I=\sum_{i\in I}\<C_i>$ for
$I\subset\{1,\dots,k\}$. The values $\bv_0\<C_i>$, $\bv_1\Ga$ (with
$\Ga\in H_1(\RX)$), and $\bv_2[C_i]$ are always well defined and coincide
with the images of the classes under the inclusion homomorphisms. The
value $\bv_1\<C_I>$ is well defined whenever the number $\ls|I|$ of
summands is even; it is represented by the collection $(1+\conj_*)y_1$ of
invariant circles, where $y_1$ is a collection of arcs in~$X$ with one
boundary point in each of the summands~$C_i$, $i\in I$ (and no other
boundary points). If, under an appropriate choice of~$y_1$, the curve
$(1+\conj_*)y_1$ is null homologous, $(1+\conj_*)y_1=\dd y_2$, then
$\bv_1\<C_I>=0$ and $(1+\conj_*)y_2$ represents~$\bv_2\<C_I>$. Similarly,
if $\bv_1\Ga=0$, i.e., $\Ga$ is homologous to $(1+\conj_*)y_1$ for a
chain~$y_1$\mnote{\Dg: edited after Zv} in~$X$, there is a chain~$y_2$
in~$X$ such that $\dd y_2=\Ga+(1+\conj_*)y_1$, and $(1+\conj_*)y_2$
represents~$\bv_2\Ga$. Finally, if $\bv_1\Ga=\bv_1\<C_I>$ for some
$0$-class $C_I$ and $1$-class $\Ga\in H_1(\RX)$, then $\bv_2(\Ga+\<C_I>)$
is defined and is represented by $(1+\conj_*)y_2$, where $\dd
y_2=\Ga+(1+\conj_*)y_1$ and $y_1$ is an appropriately chosen collection
of arcs with one boundary point in each of~$C_i$, $i\in I$, and no other
boundary points.

Elements of the form $\bv_2[C_i]$, $\bv_2\<C_I>$, $\bv_2\Ga$, and
$\bv_2(\Ga+\<C_I>)$ span~$\iH_2(X)$; with an abuse of the language we let
$\sum\bv_2x_j=\bv_2\sum x_j$, provided that the latter is well defined,
even if the summands are not. According to~\ref{C1.7.B}, Kalinin's
intersection form on~$\iH_2(X)$ is the one given by Table~\ref{tabKif},
where the intersection $\Ga\ii\Gb$ is regarded as an element of
$H_0(\RX)$, and $\aug(\Ga\ii\Gb)$ and $\aug_i(\Ga\ii\Gb)$ are,
respectively, the total intersection number and its portion falling
into~$C_i$. The notation $\Gd_{ij}$~stands for the Kronecker symbol:
$\Gd_{ii}=1$ and $\Gd_{ij}=0$ if $i\ne j$.  The intersection form extends
linearly to the classes of the form $\bv_2(\Ga+\<C_I>)$, as if $\bv_2\Ga$
and $\bv_2\<C_I>$ were well defined.

\topinsert
\table\label{tabKif}
Kalinin's intersection form
\endtable
\offinterlineskip
\def\sp{\quad}
\let\a\Ga
\let\b\Gb
\def\d[#1]{\Gd_{#1}}
\def\\{\cr\noalign{\hrule}\omit\cr}
\let\c\ii
\line{\hfill\vbox{
\ialign{\vrule height11pt depth5pt\hfil\sp$#$\sp\hfil\vrule&&
 \hfil\sp$#$\sp\hfil\vrule\cr
\noalign{\hrule}\omit\cr
             &\bv_2\<C_I>        &\bv_2\a        &\bv_2[C_i]\\
 \bv_2\<C_K> &0                  &0              &\sum_{k\in K}\d[ik]\\
 \bv_2\b     &0                  &\aug(\a\c\b)   &\aug_i(\b\c\b)\\
 \bv_2[C_k]  &\sum_{i\in I}\d[ik]&\aug_k(\a\c\a) &\d[ik]\chi(C_i)\\
}}\hfill}
\endinsert

The Bockstein homomorphism $\Sq_1\:\iH_2(X)\to\iH_1(X)$ is given by
\ref{C1.7.C}: one has $\Sq_1\bv_2[C_i]=\bv_1w_1(C_i)$ and
$\Sq_1\bv_2(\Ga+\<C_I>)=\bv_1\Ga$.

We conclude this section with a few examples of
calculation of the filtration~$\CF^*$ and Viro homomorphisms~$\vr_*$. The
first example may serve as a model for regular complete intersections
in~$\Cp{n}$ (see \cite{Kalinin}), the second one, as a local model for
involutions on smooth manifolds.

\claim[Filtrations in $\Rp{n}$ \ROM(see~\cite{Kalinin})]\label{C1.7.E}
Let $X=\Cp{n}$ and $c=\conj$. Then the filtrations~$\CF^*$ on
$H_*(\Rp{n})$ and $\CF_*$ on $H^*(\Rp{n})$ are given by
$$
\CF^r=\Sq H_{\ge{r/2}}(\Rp{n}),\qquad \CF_r=\Sq H^{\le{r/2}}(\Rp{n}),
$$
where $\Sq$ stands for both homology and cohomology total Steenrod
squares.
\qed
\endclaim

\claim[Filtration in Thom spaces\/ \ROM(see~\cite{DegtKh:2})]\label{C1.7.F}
Let $\nu$ be an $m$-dimensional vector bundle over a finite cell
complex~$F$, and let~$T$ and~$\partial T$ be the associated disk and
sphere bundles, respectively, supplied with the antipodal involution.
Then the filtrations~$\CF^*$ on $H_*(F)$ and~$\CF_*$ on $H^*(F)$
associated with Kalinin's spectral sequence of pair~$(T,\partial T)$ are
given by
$$
\CF^{p}=w(\nu)^{-1}\capprod H_{\ge p-m}(F),\qquad
 \CF^{p}=w(\nu)\cupprod H^{\le p-m}(F),
$$
where $w(\nu)=1+w_1(\nu)+w_2(\nu)+\dots$ is the total Stiefel-Whitney
class of~$\nu$. The Viro homomorphisms are
$$
\gather
\vr_p\:\CF^p\to H_p(T,\partial T),\quad
 x\mapsto\th^{-1}[w(\nu)\capprod x]_{p-m},\\
\vr^p\: H^p(T,\partial T)\to H^*(F)/\CF_{p-1},\quad
 x\mapsto w(\nu)\cupprod\th^{-1}x,
\endgather
$$
where $\th$ is the Thom isomorphism and $[x]_r$ denotes the
$r$-dimensional component of a nonhomogeneous homology class~$x$.
\qed
\endclaim

Next statement is a generalization of~\ref{C1.7.E};
it can be deduced from~\ref{C1.7.E}
using the
splitting
principle.

\claim[Filtrations in Grassmann manifolds]\label{gras}
Let $X=G_k(\C^n)$ be the Grassmann manifold of $k$-planes in~$\C^n$
and $c=\conj$ the standard conjugation. Then the
filtrations~$\CF^*$ on $H_*(G_k(\R^n))$ and $\CF_*$ on $H^*(G_k(\R^n))$
are given by
$$
\CF^r=\Sq H_{\ge{r/2}}(G_k(\R^n)),\qquad \CF_r=\Sq H^{\le{r/2}}(G_k(\R^n)),
$$
where $\Sq$ stands for both homology and cohomology total Steenrod
squares.
\qed
\endclaim

The following consequence of~\ref{gras} is, in fact, equivalent to it.
A related identity is found in~\cite{Kr:94}.

\claim[Corollary]\label{ChtoWh}
Let $X$ be a topological space with involution $c\:X\to X$ and $E$ a real
vector bundle over~$X$ in the sense of Atiyah \cite{At:Kreal} \rom(i.e.,
$E$ is a $\C$-vector bundle supplied with an anti-linear involutive
lift~$c_E$ of~$c$\rom). Then for all~$i$ one has $\bv^{2i}c_i(E)=\Sq
w^i(\R E)$, where $\R E=\Fix c_E$ is an $\R$-vector bundle over~$\Fix c$.
\qed
\endclaim

\subsection{Involutions on manifolds}\label{involutions}
Let $X$ be a smooth closed oriented $2k$-manifold and $c$ a smooth
involution, preserving the orientation if $k$ is even and reversing it
otherwise. Define the \emph{twisted intersection form}
$$
B^c\:(H_{k}(X;\Z)/\Tors)\otimes(H_{k}(X;\Z)/\Tors)\to\Z,\quad
 B^c(x,y)=x\ii c_*y.
$$
This form is symmetric; its signature~$\Gs(c)$ is called the
\emph{signature of the involution}. If $k$ is odd, the signature is $0$,
otherwise $\Gs(c)=\Gs_+-\Gs_-$, where $\Gs_\pm$ are the signatures of the
restrictions of the intersection form to the eigenspaces $H^{\pm1}_{k}$
of~$c_*$.

\theorem\label{2.1.1}
If all the components of~$F=\Fix c$ are of dimension~$k$, one has
$\Gs(c)=F\ii F$, where $F\ii F$ is the normal Euler number of~$F$ in~$X$.
\pni
\endtheorem

Since the real part $\RX$ of a complex manifold $X$ with real structure
is totally real, all the components of~$\RX$ are of dimension
$n=\dim_{\C}X$ and $\RX\ii\RX=(-1)^{n/2}\chi(\RX)$. If $X$ is supplied
with an invariant K\"ahler metric, the real part is Lagrangian.

\claim[Corollary]\label{Hst}
If $X$ is a $2k$-dimensional compact complex manifold with real structure
$\conj$, one has $\Gs(\conj)=(-1)^k\chi(\RX)$.
\qed
\endclaim

If $X$ is a closed compact (not necessarily orientable) manifold and
$\dim X=2k$ is even, one can define a $\ZZ$-version $B^c\:H_k(X)\otimes
H_k(X)\to\ZZ$ of the twisted intersection form. It is also symmetric.

\claim[Arnol$'$d lemma]\label{Arnold.lemma}
Let $X$ be a closed smooth manifold with involution~$c$ and $\dim X=2k$
even. If $\dim\Fix c\le k$, then the fundamental class~$[\Fix c]_k$ of
the union of the $k$-dimensional components of~$\Fix c$ realizes
in~$H_k(X)$ the characteristic class of $B^c\:H_k(X)\otimes H_k(X)\to\ZZ$.
\pni
\endclaim

\remark{Remark}
In general the characteristic class of~$B^c$ is the $k$-dimensional
component of $\inj_!(v(\tau)\,v^{-1}(\nu))$, where $\inj\:\Fix c\to X$ is
the inclusion, $\tau$ and~$\nu$ are, respectively, the tangent and normal
bundles of~$\Fix c$, and $v$ stands for the total Wu class. See~\cite{CM}
for details.
\endremark

\end



\ifx\rokhloaded\undefined\input rokh.def \fi

\begin

\def\Bil{\operatorname{Bil}}
\def\Qu{\operatorname{Qu}}
\def\Pin{\operatorname{Pin}}
\def\ind{\operatorname{ind}}
\def\Card{\operator{Card}}
\def\sint{\operator{int}}
\def\gl{\frak l}
\def\gm{\operator{\frak{gm}}}
\def\gM{\frak M}
\def\PD{\operator{D\null}}

\let\dd\partial

\appendix{Integral lattices and quadratic forms}\label{apB}

\noindent
In this appendix we\mnote{\Dg: edited after Zv} give a brief introduction
into arithmetic of integral lattices and finite quadratic forms, which is
a commonly used tool in topology of real algebraic varieties. Proofs,
mostly straightforward, can be found in any advanced textbook,
e.g.,~\cite{MilnorHusemoller} and~\cite{Scharlau}.

\subsection{Finite quadratic forms}\label{forms}
A \emph{finite symmetric bilinear form} is a pair $(\CL,b)$, where $\CL$
is a finite abelian group and $b\:\CL\otimes\CL\to\Q/\Z$ is a symmetric
bilinear form.  A \emph{finite quadratic form} on~$\CL$ is a map
$q\:\CL\to\Q/2\Z$ satisfying the relations $q(nx)=n^2q(x)$, $x\in\CL$,
$n\in\Z$, and $q(x+y)-q(x)-q(y)=2b(x,y)$, $x,y\in\CL$, where $b$ is a
symmetric bilinear form; $b$ is determined by~$q$ and is called the
bilinear form \emph{associated with}~$q$, whereas $q$ is called a
\emph{quadratic extension} of~$b$.

If a quadratic form~$q$ (bilinear form~$b$) is understood, it is often
abbreviated as $q(x)=x^2$ (respectively, $b(x,y)=xy$).

If $q_1$, $q_2$ are two quadratic extensions of the same bilinear
form~$b$, the difference $q_1-q_2$ is a linear form.  Furthermore, for any
$x\in\CL$ one has $q_i(x)=b(x,x)\bmod\Z$, $i=1,2$; hence, $q_1-q_2$ takes
values in $\ZZ=\Z/2\Z\subset\Q/2\Z$. Thus, the set $\Qu(\CL,b)$ of
quadratic extensions of a given bilinear form~$b$ on~$\CL$ is an affine
space over $\Hom(\CL,\ZZ)$. In particular, unless $\CL$ has $2$-torsion,
any bilinear form~$b$ on~$\CL$ has a unique quadratic extension.

A symmetric bilinear form~$b$ is called \emph{nondegenerate} if the
\emph{correlation homomorphism} $\CL\to\Hom(\CL,$ $\Q/\Z)$, $x\mapsto
b(x,\,\cdot\,)$, is bijective. A quadratic form is called
\emph{nondegenerate} if so is its associated bilinear form. In general,
the kernel of the correlation homomorphism is denoted by $\Ker\CL$ or
$\CL\ort$.

The nondegenerate finite quadratic and bilinear forms constitute monoids
$\Bil$ and $\Qu$, respectively, the addition operation being the
orthogonal sum. A set of generators of $\Bil$ is given, e.g.,
in~\cite{Wall}; for each generator $\CL$ is either cyclic or (not
orthogonal) sum of two cyclic groups. Generators of~$\Qu$ can be obtained
by lifting those of~$\Bil$. In this survey we deal explicitly only with
forms on groups of period~$2$ (however, proofs of the results on real
Enriques surfaces make an extensive use of forms on groups of period~$4$,
whose theory can be found in~\cite{DIK}); they constitute direct summands
$\Bil_2$ and $\Qu_2$ of $\Bil$ and $\Qu$, respectively. The Wall
generators of $\Bil_2$ are the following:
$$
\CL=\Z_{2}\:\quad
 a_{2}=\Bigl[\frac1{2}\Bigr];\qquad
\CL=\Z_{2}\times\Z_{2}\:\quad
 u_{2}=\biggl[\matrix0&\frac1{2}\\\frac1{2}&0\endmatrix\biggr].
$$
Up to isomorphism $a_2$ and~$u_2$ have two quadratic extensions each:

\vglue\abovedisplayskip \bgroup\openup\jot
\centerline{%
for $a_2$:\quad$\Bigl[\pm\dfrac12\Bigr]$,\qquad for $u_2$:\quad
$\CU_2=\biggl[\matrix0&\frac12\\\frac12&0\endmatrix\biggr]$ and
$\CV_2=\biggl[\matrix1&\frac12\\\frac12&1\endmatrix\biggr]$.%
}\vglue\belowdisplayskip \egroup

The \emph{Brown invariant} $\Br q$ of a nondegenerate finite quadratic
form~$q$ on~$\CL$ is the $(\bmod\,8)$-residue defined via
$$
\tsize \exp\bigl(\frac14i\pi\Br q \bigr)= \ls|\CL|^{-\frac12}\sum_{x\in
\CL}\exp\bigl(i\pi q(x)\bigr).
$$

\claim\label{smallBr}
One has $\Br\Bigl[\pm\dfrac12\Bigr]=\pm 1$, $\Br\CU_2=0$, and
$\Br\CV_2=4$.
\qed
\endclaim

A finite group~$\CL$ of period~$2$ is a $\ZZ$-vector space. If $\CL$ is
equipped with a finite quadratic form, it is called a \emph{quadratic
space}. A quadratic (bilinear) form on a $\ZZ$-space takes values in
$(\12\Z)/2\Z\subset\Q/2\Z$ (respectively, in $(\12\Z)/\Z\subset\Q/\Z$).
By means of the canonical identifications $(\12\Z)/2\Z=\Z_4$ and
$(\12\Z)/\Z=\ZZ$, $\12\mapsto1$, one can consider quadratic (bilinear)
forms on $\ZZ$-spaces as maps $\CL\to\Z_4$ (respectively,
$\CL\otimes\CL\to\ZZ$).

A vector $w\in\CL$ is called a \emph{characteristic element} of~$\CL$ if
$xw=x^2\bmod2$ for any $x\in\CL$.
Since on groups of period~$2$ the map $x\mapsto x^2\bmod2$ is linear,
characteristic elements exist and form a coset $\bmod(\Ker\CL)$; it is
called the \emph{characteristic coset} of~$\CL$. If $\CL$ is
nondegenerate, a characteristic element is unique.  The subgroup of~$\CL$
spanned by $\Ker\CL$ and a characteristic element will be denoted by
$\GO=\GO(\CL)$.

A quadratic space is called \emph{even} if $0$ is its characteristic
element, i.e., if for all $x\in\CL$ one has $x^2=0\bmod2$. Otherwise, it
is called \emph{odd}.

\claim\label{2forms}
The monoid of nondegenerate finite quadratic forms on groups of period~$2$
is generated by~$\Cal A^\pm=[\pm\12]$, $\CU_2$, and~$\CV_2$, which are
subject to the relations $2\CU_2=2\CV_2$, $3\Cal A^\pm=\Cal
A^\mp\oplus\CV_2$, and $\Cal A^\pm\oplus\Cal A^+\oplus\Cal A^-=\Cal
A^\pm\oplus\CU_2$.

Two nondegenerate quadratic spaces are isomorphic if and only if they have
the same rank, parity, and Brown invariant.
\pni
\endclaim

A quadratic form $q$ on~$\CL$ is called \emph{informative} if
$q|_{\CL\ort}=0$. The notion of Brown invariant extends to informative
forms: since $q$ vanishes on~$\CL\ort$, it descends to a quadratic form
$q'\:\CL/\CL\ort\to\Q/2\Z$, and one lets $\Br q=\Br q'$. The Brown
invariant is additive\rom: given a pair~$(\CL_i,q_i)$, $i=1,2$, of
informative quadratic forms, one has $\Br(\CL_1\oplus \CL_2,q_1\oplus
q_2)=\Br(\CL_1,q_1)+\Br(\CL_2,q_2)$.

A subform~$\CN$ of an informative quadratic form~$(\CL,q)$ is called
\emph{informative} if $\CN\ort\subset\CN$ and $q|_{\CN\ort}=0$.  (Clearly,
an informative subform is an informative form; hence, its Brown invariant
is well defined.) A subform $\CS\subset\CL$ is called \emph{isotropic}
(more precisely, $q$-isotropic) if $q|_{\CS}=0$. Clearly, a subform $\CN$
of a nondegenerate form $(\CL,q)$ is informative if and only if $\CN\ort$
is isotropic.

The following is a direct consequence of a Gauss summation formula.

\claim\label{Sort/S}
If $\CS$ is a $q$-isotropic subform of a nondegenerate quadratic form
$(\CL,q)$, then $q$ descends to a nondegenerate quadratic form on
$\CS\ort/\CS$ and $\Br(S\ort/S)=\Br\CL$.
\pni
\endclaim

\corollary\label{inf-space}
If $\CN$ is an informative subform of an informative quadratic
form~$(\CL,q)$, then $\Br(\CN,q|_\CN)=\Br(\CL,q)$.
\qed
\endcorollary

\remark{Remark}
The notion of informative subspace still makes sense if the quadratic
form~$q$ is only defined on~$\CN$ (while the bilinear form is defined on
the whole space~$\CL$). In this case Corollary~\ref{inf-space} states
that $\Br q'=\Br q$ for any extension~$q'$ of~$q$ to~$L$.
\endremark

\claim\label{Br}
For any informative quadratic space $(\CL,q)$ one has\rom:
\roster
\item\local1
$\Br q=\dim(\CL/\CL\ort)\bmod2$\rom;
\item\local2
$\Br q= q(w)\bmod4$ for any characteristic element $w\in\CL$\rom;
\item\local3
$\Br(q+v)=\Br q-2q(v)$ for any $v\in\CL$
\rom(where $q+v$ stands for the quadratic form $x\mapsto
q(x)+2(vx)$\rom)\rom;
\item\local4
$\Br q=0$ if and only if $(\CL,q)$ is\/ \emph{null cobordant}, i.e., there
is a subspace $\CH\subset\CL$ such that $\CH\ort=\CH$ and $q|_\CH=0$.
\pni
\endroster
\endclaim

\subsection{Integral lattices}\label{lattices}
An \emph{\rom(integral\rom) lattice} is a free abelian group $L$ of finite
rank supplied with a symmetric bilinear form $b\:L\otimes L\to\Z$. As in
the case of finite forms, we often abbreviate $b(x,y)=xy$ and
$b(x,x)=x^2$. A lattice~$L$ is called \emph{even} if $x^2=0\bmod2$ for
all $x\in L$; otherwise, $L$ is called \emph{odd}. Let
$L\spcheck=\Hom(L,\Z)$ be the dual abelian group. $L$ is called
\emph{nondegenerate} (\emph{unimodular}) if the \emph{correlation
homomorphism} $L\to L\spcheck$, $x\mapsto b(x,\,\cdot\,)$, is monic
(respectively, one-to-one). The cokernel of the correlation homomorphism
is called the \emph{discriminant group} of~$L$ and denoted by $\discr L$
or~$\CL$.

The bilinear form on~$L$ extends to a $\Q$-valued bilinear form on
$L\otimes\Q$. If $L$ is nondegenerate, there is a natural identification
$L\spcheck=\{x\in L\otimes\Q\,|\,\text{$xy\in\Z$ for $y\in L$}\}$ and the
correlation homomorphism is the embedding $L\hookrightarrow L\spcheck$.
Hence, $\CL$ is finite and the form on~$L$ extends to a $\Q$-valued form
on $L\spcheck$, which, in turn, induces a finite symmetric bilinear form
$b\:\CL\otimes\CL\to\Q/\Z$. The order of~$\CL$ is equal to the absolute
value of the determinant $\det L$ of the Gram matrix of~$L$. If, in
addition, $L$ is even, $b$ is associated with the finite quadratic form
$q\:\CL\to\Q/2\Z$ defined via $q(x+L)=x^2\bmod2\Z$, $x\in L\spcheck$.

\claim[Van der Blij formula]\label{VdB}
For a nondegenerate even integral lattice~$L$ one has
$\Br\CL=\Gs(L)\bmod8$.
\pni
\endclaim

The following proposition is an immediate consequence of \ref{VdB}:

\claim[Corollary]\label{10.2.B}
Let~$L$ be an even nondegenerate integral lattice whose discriminant group
$\CL=\discr L$ is of period~$2$. Then\rom:
\roster
\item\local1
if $L$ is unimodular \rom(i.e., $\CL=0$\rom), then $\Gs(L)=0\bmod8$\rom;
\item\local2
if $\det L=\pm2$ \rom(i.e., $\CL=\ZZ$\rom), then $\Gs(L)=\pm 2\bmod8$\rom;
\item\local3
if $\det L=\pm4$ \rom(i.e., $\CL=\ZZ\oplus\ZZ$\rom) and the discriminant
form on $\CL$ is odd, then $\Gs(L)=0,\pm2\bmod8$\rom;
\item\local4
if the discriminant form on $\CL$ is even, then $\Gs(L)=0\bmod4$.\qed
\endroster
\endclaim

\subsection{Involutions on unimodular lattices}\label{inv.lattice}
Assume that $L$ is a unimodular lattice and $c\:L\to L$ an involutive
auto-isometry of~$L$. Let $L^{\pm}\subset L$ be the eigenlattices of $c$,
$q^\pm=\discr L^{\pm}$, and $J=L/(L^{+}\oplus L^{-})$. As groups, $J$ and
$q^\pm$ are isomorphic $\ZZ$-vector spaces. The eigenlattices are the
orthogonal complements of each other and one has $\rk L^++\rk L^-=\rk L$
and
$$
 \log_2\ls|\discr L^\pm|=\dim_{\ZZ}J
   =\dim_{\ZZ}L\otimes\ZZ-\dim_{\ZZ}\{L\otimes\ZZ\}^c,
$$
where $\{L\otimes\ZZ\}^c=\{x\in L\otimes\ZZ\,|\,cx=l\}$.

Next statement is straightforward:

\claim\label{10.4.C}
The form~$q^+$ is even if and only if the characteristic elements of~$L$
coincide with those of the twisted form $b^c\:L\otimes L\to\ZZ$,
$x\otimes y\mapsto b(x,cy)$.
\qed
\endclaim


\end



\ifx\rokhloaded\undefined\input rokh.def \fi

\begin

\def\Pin{\operatorname{Pin}}
\def\ind{\operatorname{ind}}
\def\Card{\operator{Card}}
\def\sint{\operator{int}}
\def\gl{\frak l}
\def\gm{\operator{\frak{gm}}}
\def\gM{\frak M}
\def\PD{\operator{D\null}}

\let\dd\partial

\appendix{The Rokhlin-Guillou-Marin form}\label{RGMform}

\noindent
In this appendix we briefly remind the construction of the
Rokhlin-Guillou-Marin form (see~\cite{GM}) of a characteristic surface in
a $4$-manifold.

Let $Y$ be an oriented closed smooth $4$-manifold and $U$ a
\emph{characteristic surface} in~$Y$, i.e., a smooth closed
$2$-submanifold with $[U]=\PD u_2(Y)$ in $H_2(Y)$. Denote by
$i\:U\hookrightarrow Y$ the inclusion and let $K=\Ker[i_*\:H_1(U)\to
H_1(Y)]$.  Then there is a natural quadratic extension $\gm\:K\to\Z_4$ of
the intersection index form, which is called the
\emph{Rokhlin-Guillou-Marin form} of~$(Y,U)$.  It can be defined as
follows. Pick a class $x\in K$ and realize it by a union~$\gl$ of disjoint
simple closed smooth loops in~$U$. It spans an immersed surface~$\gM$
in~$Y$, which can be chosen normal to~$U$ along $\gl=\dd\gM$ and
transversal to~$U$ at its inner points.  (Such a surface is called a
\emph{membrane}.) Consider a normal line field~$\xi$ on~$\gl$ tangent
to~$U$ and define the \emph{index} $\ind\gM\in\Z$ as the obstruction to
extending~$\xi$ to a normal line field on~$\gM$. (Since
$\Gt\gM\oplus\nu\gM$ is an oriented vector bundle, the obstruction is a
well defined integer. If $\gl$ is two-sided in~$U$, the index equals twice
the index defined via vector fields instead of line fields.) Then
$\gm(x)=\ind\gM+2\Card(\sint\gM\cap F)\bmod4$.

\theorem\label{2.4.1}
Let $Y$, $U$, and $(K,\gm)$ be as above. Then $(K,\gm)$ is an informative
subspace of $H_1(U)$ and
$$
2\Br\gm=\Gs(Y)-U\ii U\bmod16,
$$
where $U\ii U$ stands for the normal Euler number of~$U$ in~$Y$.
\endtheorem

For the proof in the case when $Y$ is simply connected (or, more
generally, when $K=H_1(U)$) see~\cite{GM}, for the general case see,
e.g.,~\cite{DIK}. Rokhlin considered the case when the characteristic
surface is orientable and used the $\operatorname{Arf}$-invariant: in
this case $4\operatorname{Arf} =\Br$.

\remark{Remark}
There is an alternative construction of the Rokhlin-Guillou-Marin form.
Since $U$ is characteristic, $Y\sminus U$ admits a $\Spin$-structure
which does not extend through any component of~$U$. Its restriction to
the boundary of a tubular neighborhood of~$U$ induces in a natural way a
$\Pin^-$-structure on~$U$ (cf\.~\cite{FinashinRGM}), which defines a
quadratic form~$q$ on~$H_1(U)$. It is not difficult to see that $q$ is
well defined up to adding elements of $\Im[i^*\:H^1(Y)\to H^1(U)]$ and,
hence, its restriction to~$K$ does not depend on the choice of a
$\Spin$-structure; it coincides with~$\gm$. (This gives an alternative
proof of \ref{2.4.1}, communicated to us by C.~Bailly: the quadratic forms
arising from $\Pin^-$-structures are isomorphic to each other since they
have the same Brown-invariant (see \ref{2forms}); these forms provide all
the quadratic extensions of~$\gm$ from~$K$ to $H_1(U)$ and, as it is
shown in~\cite{Bailly}, this implies that $K$ is informative.)
\endremark

\end



\ifx\rokhloaded\undefined\input rokh.def \fi

\def\mnote#1{}

\document

\end





\endgroup


\widestnumber\key{Kho1}

\Refs

\ref{Ak}
\by S.~Akbulut
\paper On quotients of complex surfaces under complex conjugation
\jour J. Reine Angew. Math.
\vol 447
\yr 1994
\pages 83--90
\endref\label{Akbulut}

\ref{A1}
\by V.~I.~Arnold
\paper On arrangement of ovals of real plane algebraic curves,
the involutions of four-dimensional smooth manifolds,
and the arithmetic of integral quadratic forms,
\jour Funkc. Anal. i Prilozhen.
\vol 5
\yr 1971
\issue 3
\pages 1--9
\endref\label{Arnold:break}

\ref{A2}
\by V.~I.~Arnold
\paper The index of a singular point of a vector field,
the Petrovkii-Oleinik inequalities and mixed Hodge structures
\jour Funkc. Anal. i Prilozhen.
\vol 12
\yr 1978
\issue 1
\pages 1--14
\endref\label{ArnoldPO}

\ref{A3}
\by V.~I.~Arnold
\paper The branched covering $\CP^2\to S^4$,
hyperbolocity and projective topology
\jour Sibirsk. Mat. Zh.
\yr 1988
\vol 29
\issue 5
\pages 36--47
\endref
\label{Arnold}

\ref{A4}
\by V.~I.~Arnold
\paper Topological content of the Maxwell theorem
on multipole representation of spherical functions
\jour Topol. Methods Nonlinear Anal.
\vol 7
\issue 2
\yr 1996
\pages 205--217
\endref
\label{ArOnMaxwell}

\ref{AO}
\by V.~I.~Arnold, O.~A.~Oleinik
\paper Topology of real algebraic manifolds
\jour Vestnik Mosk. Univ., Ser. I, Mat. i Mekh.
\yr 1979
\vol 6
\pages 7--17
\endref
\label{ArnOlei}

\ref{At}
\by M.~F.~Atiyah
\paper $K$-theory and reality
\jour Quart. J. Math. Oxford Ser. \rom{(2)}
\vol 17
\yr 1966
\pages 367--386
\endref\label{At:Kreal}

\ref{AD}
\by M.~F.~Atiyah, J.~L.~Dupont
\paper Vector fields with finite singularities
\jour Acta Math. Stockh.
\vol 128
\yr 1972
\pages 1--40
\endref\label{AD}

\ref{Ba}
\by C.~Bailly, A.~Vdovina
\paper Sous-espaces d\'eterminant l'invariant de Arf et un Th\'eor\`eme de
Rokhlin sur la signature
\toappear
\endref\label{Bailly}

\ref{Bi}
\by F.~Bihan
\book Constructions Combinatoires de Surfaces Alg\'ebriques
R\'eelles
\bookinfo Th\`ese
\yr 1998
\endref
\label{Bihan}

\ref{BPV}
\by W.~Barth, C.~Peters, A.~van~de~Ven
\book Compact complex surfaces
\bookinfo Ergebnisse der Mathemtik und ihrer
Grenzgebiete \rom{(3)}
\publ Springer-Verlag
\publaddr Berlin-New York
\yr 1984
\endref\label{BPV}

\ref{BH}
\by A.~Borel, A.~Haefliger
\paper La classe d'homologie fondamentale d'un espace analytique
\jour Bull. Soc. math. France
\vol 89
\yr 1961
\pages 461--513
\endref
\label{BH}

\ref{Br}
\by G.~E.~Bredon
\book Introduction to compact transformation groups
\publ Academic Press
\publaddr New York-London
\yr 1972
\endref
\label{Bredon}

\ref{C}
\by J.~Cerf
\paper Sur les diff\'eomorphismes de la sph\'ere de dimension trois,
$\Gamma_4=0$
\jour Lecture Notes Math.
\vol 53
\endref\label{Cerf}

\ref{CP}
\by C.~Ciliberto, C.~Pedrini
\paper Annibale Comessatti and real algebraic geometry
\jour   R.C. Circolo Mat. Palermo. Ser.~\rom{II}. Suppl. no.~36
\yr 1994
\pages  71--102
\endref\label{CilibertoPedrini}

\ref{Ch}
\by S.~S.~Chern
\book Complex manifolds
\bookinfo Russian translation
\publ Izd. Inostr. Lit.
\publaddr Moscow
\yr 1961
\endref\label{Chern}

\ref{Co1}
\by A.~Comessatti
\paper Sulla connessione delle superficie algebriche reali
\inbook Verhandl. Intern. Math. Kongr. Z\"urich
\vol 2
\pages 169
\endref\label{Comessatti}

\ref{Co2}
\by A.~Comessatti
\paper Fondamenti per la geometria sopra le superficie razionali dal punto
di vista reale
\jour Math. Annalen
\vol 43
\yr 1912
\pages 1--72
\endref\label{Comessatti:1}

\ref{Co3}
\by A.~Comessatti
\paper Sur la connessione delle superficie razionali reale
\jour Annali di Matematica
\vol 23
\yr 1914
\pages 215--285
\endref\label{Comessatti:2}

\ref{Co4}
\by A.~Comessatti
\paper Reele Fragen in der algebraischen Geometrie
\jour Jahresbericht d. Deut. Math. Vereinigung
\vol 41
\yr 1932
\pages 107--134
\endref\label{Comessatti:3}

\ref{CM}
\by P.~E.~Conner, E.~Y.~Miller
\book Equivariant self-intersection
\bookinfo Preprint
\yr 1979
\endref\label{CM}

\ref{DKh}
\by V.~I.~Danilov, A.~G.~Khovanskii
\paper Newton polygons and an algorithm for computing Hodge-Deligne numbers
\jour Izv. Akad. Nauk SSSR, Ser. Matem.
\vol 50
\yr 1986
\issue 5
\pages 925--945
\transl English transl.
\jour Math. USSR Izv.
\vol 29
\yr 1987
\pages 279--298
\endref\label{KhoDanilov}

\ref{Deg1}
\by A.~Degtyarev
\paper Stiefel orientations on a real algebraic variety
\inbook Proc. Rennes Symp. on Real Alg. Geom.
\bookinfo Lect. Notes in Math.
\vol 1524
\yr 1992
\pages 205--220
\endref\label{Degt:Stiefel}

\ref{Deg2}
\by A.~Degtyarev
\paper Cohomology approach to killing structures on Steenrod bundles
\inbook Geometry and Topology,~\rom{I}
\bookinfo Adv. Soviet Math.
\vol 18
\yr 1994
\pages 1--22
\endref\label{Degt:ks}

\ref{Deg3}
\by A.~Degtyarev
\paper On the Pontrjagin-Viro form
\jour Advances Math. Sciences, Rokhlin's Memorial volume
\toappear
\endref\label{Dg:PVF}

\ref{DIK}
\by A.~Degtyarev, I.~Itenberg, V.~Kharlamov
\book Real Enriques Surfaces
\toappear
\endref\label{DIK}

\ref{DK1}
\by A.~Degtyarev, V.~Kharlamov
\paper Halves of a real Enriques Surface
\jour Comment. Math. Helvetici
\vol 71
\yr 1996
\pages 628--663
\endref\label{DegtKh:2}

\ref{DK2}
\by A.~Degtyarev, V.~Kharlamov
\paper Topological classification of real Enriques surfaces
\jour Topology
\vol 35
\issue 3
\yr 1996
\pages 711--729
\endref\label{DegtKh:1}

\ref{DZ}
\by A.~Degtyarev, V.~I.~Zvonilov
\paper Rigid isotopy classification of real algebraic curves of
bidegree $(3,3)$ on quadrics
\jour Mat. Zametki
\vol 66
\issue 6
\yr 1999
\pages 810-815
\endref\label{DZ}

\ref{Don}
\by S.~Donaldson
\paper Yang-Mills invariants of smooth four-manifolds
\inbook Geometry of low-dimens\-ion\-al manifolds
\publ Cambridge Univ. Press
\vol  1
\yr 1990
\pages 5--40
\endref\label{Donald}

\ref{Ed}
\by A.~L.~Edmonds
\paper Orientability of fixed point sets
\jour Proc. Amer. Math. Soc.
\vol 82
\yr 1981
\issue 1
\pages 120--124
\endref
\label{Edmonds}

\ref{EKh}
\by Ya.~M.~Eliashberg, V.~Kharlamov
\paper On the number of complex points of a real surface in a complex one
\inbook Proceedings of Leningrad International
Topology conference, 1982, Leningrad
\pages 143--148
\endref\label{EKh}

\ref{Fi1}
\by T.~Fiedler
\paper Pencils of lines and the topology of real algebraic curves
\jour Izv. Akad. Nauk SSSR Ser. Mat.
\vol 46
\issue 4
\yr 1982
\pages 853--863
\endref
\label{Fie:cong}

\ref{Fi2}
\by T.~Fiedler
\paper New congruences in the topology of real plane algebraic curves
\jour Dokl. Akad. Nauk SSSR
\vol 270
\issue 1
\pages 56--58
\yr 1983
\endref
\label{Fie:pencil}

\ref{Fi3}\mnote{Kh:new}
\by T.~Fiedler
\paper Additional inequalities in the topology of real plane algebraic curves
\jour Izv. Akad. Nauk SSSR Ser. Mat.
\vol 49
\issue 4
\pages 874--883
\yr 1985
\endref
\label{Fie:ro}

\ref{F1}
\by S.~Finashin
\paper Rokhlin conjecture and quotients of complex surfaces by complex
conjugation
\jour J. Reine Angew. Math.
\vol 481
\yr 1996
\pages 55--71
\endref\label{Finashin1}

\ref{F2}
\by S.~Finashin
\paper Decomposability of quotients by complex conjugation for
 rational and Enriques surfaces
\jour Topology Appl.
\vol 79
\yr 1997
\issue 2
\pages 121--128
\endref\label{Finashin2}

\ref{F3}
\by S.~Finashin
\paper A $Pin^-$ cobordism invariant and a
generalization of the Rokhlin signature congruence
\jour Algebra i Analiz
\vol 2
\issue 4
\yr 1990
\pages 242--250
\endref\label{FinashinRGM}

\ref{FS}
\by S.~Finashin, E.~Shustin
\paper On imaginary plane curves and $\Spin$-quotients of complex surfaces
by complex conjugation
\jour AMS Translations, Ser.~2
\vol 180
\yr 1997
\pages 93--101
\endref\label{FinShust}

\ref{G1}
\by D.~A.~Gudkov
\paper Construction of a new series of $M$-curves
\jour Dokl. Akad. Nauk SSSR
\vol 200
\issue 6
\yr 1971
\pages 1269--1272
\endref\label{Gudkov71}

\ref{G2}
\by D.~A.~Gudkov
\paper The topology of real projective algebraic varieties
\jour Uspekhi. Matem. Nauk
\vol 29
\yr 1974
\issue 4
\pages 3--79
\endref
\label{Gudkov:survey}

\ref{G3}
\by D.~A.~Gudkov
\paper  Topology of algebraic curves on a hyperboloid
\jour   Uspekhi Mat. Nauk
\vol    34
\issue  6~(210)
\yr 1979
\pages  26--32 \lang   Russian
\moreref\nofrills English transl. in
\jour   Russian Math. Surveys
\vol    34
\yr     1979
\endref\label{Gudkov}

\ref{GS}
\by D.~A.~Gudkov, E.~I.~Shustin
\paper Classification of nonsingular eight-order curves on an ellipsoid
\inbook Methods of the Qualitative Theory of Differential
Equations \ed E.~A.~Leontovich-Andronova
\publ   Gor$'$kov. Gos. Univ.
\publaddr Gor$'$ki
\yr 1980
\pages  104--107 \lang   Russian
\endref\label{Gudkov-Shustin}

\ref{GU}
\by D.~A.~Gudkov, A.~K.~Usachev
\paper Nonsingular curves of small orders on a hyperboloid
\inbook Methods of the Qualitative Theory of Differential Equations \ed
E.~A.~Leontovich-Andronova
\publ   Gor$'$kov. Gos. Univ.
\publaddr Gor$'$ki
\yr 1980
\pages  96--103 \lang   Russian
\endref\label{Gudkov-Usachev}

\ref{GUt}
\by D.~A.~Gudkov, G.~A.~Utkin
\paper Topology of curves of order 6 and surfaces of order 4
\inbook Uch. Zapiski Gor$'$kov. Univ.
\vol 87
\yr 1969
\endref\label{Gudkov69}

\ref{GM}
\by L.~Guillou, A.~Marin
\paper Une extension d'in th\'eor\`eme de Rohlin sur la signature
\jour C. R. Acad. Sci. Paris
\yr 1977
\vol 285
\pages A95--A98
\endref
\label{GM}

\ref{Ha}
\by R.~Hartshorne
\book Algebraic Geometry
\bookinfo Graduate Text in Math.
\vol    49
\yr 1977
\publ   Springer-Verlag
\endref\label{Hartshorne}

\ref{H1}
\by F.~Hirzebruch
\paper
\"Uber vierdimensionale Riemannsche Fl\"achen mehrdeufiger analiti\-scher
Functionen von zwei komplexen Ver\"anderlichen
\jour Math. Ann.
\yr 1953
\vol    126
\issue  1
\pages  1--22 
\endref\label{Hirz53}

\ref{H2}
\by F.~Hirzebruch
\book Topological Methods in Algebraic Geometry
\publ Springer-Verlag
\yr   1978
\endref\label{Hirzebruch}

\ref{Ho}
\by E.~Horikawa
\paper On deformations of quintic surfaces
\jour Invent. Math.
\vol 31
\yr 1975
\pages 43--85
\endref\label{Horikawa}

\ref{IKh}
\by I.~Itenberg, V.~Kharlamov
\paper Towards the maximal number of components
of a nonsingular surface of degree 5 in $\RP^3$
\jour Amer. Math. Soc. Trans. \rom{(2)}
\vol 173
\yr 1996
\pages 111--118
\endref\label{ItenKh}

\ref{IV1}
\by I.~Itenberg, O.~Ya.~Viro
\paper Patchworking algebraic curves disproves the Ragsdale conjecture
\jour Math. Intelligencer
\vol 18
\yr 1996
\issue 4
\pages 19--28
\endref\label{VIten}

\ref{IV2}
\by I.~Itenberg, O.~Ya.~Viro
\paper Maximal real algebraic hypersurfaces of projective space
\jour preprint
\endref\label{ItenViro}

\ref{Ka}
\by I.~Kalinin
\paper Cohomology characteristics of real algebraic hypersurfaces
\jour Algebra i Analis
\vol    3
\yr     1991
\pages 313--332
\lang Russian
\endref\label{Kalinin}

\ref{Kh1}
\by V.~Kharlamov
\paper The maximal number of components of a 4th degree surface in $\RP^3$
\jour Funkc. Anal. i Prilozhen.
\vol 6
\yr 1972
\issue 4
\pages 101
\endref\label{Kh:10}

\ref{Kh2}
\by V.~Kharlamov
\paper A generalized Petrovskii inequality
\jour Funkc. Anal. i Prilozhen.
\vol    8
\issue  2
\yr     1974
\pages  50--56
\lang Russian
\transl\nofrills English transl. in
\jour Functional Anal. Appl.
\vol    8
\issue  2
\yr     1974
\pages  132--137
\endref\label{Kh:petr1}

\ref{Kh3}
\by V.~Kharlamov
\paper The topological type of nonsingular surfaces in $\Rp3$ of
degree~4
\jour Funkc. Anal. i Prilozhen.
\vol    10
\issue  4
\yr 1976
\pages  55--68
\lang Russian
\transl\nofrills English transl. in
\jour Functional Anal. Appl.
\vol 10
\yr 1976
\issue 4
\pages 295--305
\endref\label{Kh:K3}

\ref{Kh4}
\by V.~Kharlamov
\paper A generalized Petrovskii inequality. \rom{II}
\jour Funkc. Anal. i Prilozhen.
\vol    9
\issue  3
\yr     1975
\pages  93--94
\lang Russian
\transl\nofrills English transl. in
\jour Functional Anal. Appl.
\vol    9
\issue  3
\yr     1975
\pages  266--268
\endref\label{Kh:petr2}

\ref{Kh5}
\by V.~Kharlamov
\paper Additional congruences for the Euler characteristic
of even-dimensional real algebraic varieties
\jour Funkc. Anal. i Prilozhen.
\vol 9
\yr 1975
\issue 2
\pages 51--60
\endref\label{Kh75}

\ref{Kh6}
\by V.~Kharlamov
\paper Isotopy types of non-singular surfaces of degree 4 in $\RP^3$
\jour Funct. Anal. Appl.
\yr 1978
\vol 12
\pages 86--87
\endref\label{Kh:iso}

\ref{Kh7}
\by V.~Kharlamov
\paper On the classification of non-singular surfaces
of degree 4 in $\RP^3$ with respect to rigid isotopies
\jour Funct. Anal. Appl.
\yr 1984
\vol 18
\issue 1
\pages 49--56
\endref\label{Kh:rigid}

\ref{Kh8}
\by V.~Kharlamov
\paper Nonamphicheiral surfaces of degree 4 in $\RP^3$
\jour Lect. Notes Math.
\vol 1346
\yr 1988
\pages 349--356
\endref\label{Kh:chiral}

\ref{Kh9}
\by V.~Kharlamov
\paper On a number of components of an $M$-surface of degree 5 in $\RP^3$
\inbook Proc. of XVI soviet algebraic conference, Leningrad
\yr 1981
\pages 353--354
\endref\label{KhLen}

\ref{Kh10}
\by V.~Kharlamov
\paper Estimates of Betti numbers in topology of real algebraic surfaces
\inbook Proc. Internat. Topological Conference, Oberwolfach
\yr 1987
\pages 12--13
\endref\label{KhOber}

\ref{Kho1}
\by A.~G.~Khovansky
\paper The index of a polynomial vector field
\jour Funkc. Anal. i Prilozhen.
\vol 13
\yr 1979
\issue 1
\pages 49--58
\endref\label{KhoPO}

\ref{Kho2}
\by A.~G.~Khovansky
\paper Boundary indices of polynomial $1$-forms with homogeneous components
\jour Algebra i Analiz
\vol 10
\yr 1998
\issue 3
\pages 193--222
\endref\label{KhoPO2}

\ref{Kl}
\by F.~Klein
\paper Ueber Fl\"achen dritter Ordnung
\jour Math. Ann.
\vol   6
\yr    1873
\pages 551--581
\endref\label{Klein}

\ref{KS}
\by K.~Kodaira, D.~C.~Spencer
\paper Groups of complex line bundles over compact K\"ahler
varieties. Divisor class groups on algebraic varieties.
\jour   Proc. Nat. Acad. Sci. U.S.A.
\vol    39
\yr 1953
\pages  868--877
\endref\label{Kodaira-Spencer}

\ref{Ko}
\by A.~Korchagin
\paper Classification of real algebraic curves on toric varieties,
Preliminary report
\inbook 911-th AMS Meeting, Program
\bookinfo Louisiana State Univ., Baton Rouge, Apr.~19--21
\yr 1996
\page 397
\endref\label{Korchagin}

\ref{Kr1}
\by V.~A.~Krasnov
\paper The Harnack-Thom inequality for mappings of real algebraic varieties
\jour Izv. Akad. Nauk SSSR Ser. Mat.
\vol 47
\yr 1983
\issue 2
\pages 268--297
\endref
\label{Kr:GM}

\ref{Kr2}
\by V.~A.~Krasnov
\paper On the equivariant Grothendieck cohomology of
a real algebraic variety and its application
\jour Izv. Ross. Akad. Nauk Ser. Matem.
\vol 58
\yr 1994
\issue 3
\pages 36--52
\endref
\label{Kr:94}

\ref{Ku}
\by N.~Kuiper
\paper The quotient space of $\CP^2$ by complex conjugation is
 the $4$-sphere
\jour Math. Ann.
\yr 1974
\vol 208
\pages 175--177
\endref\label{Kuiper}

\ref{La}
\by Hon-Fei Lai
\paper Characteristic classes of real manifolds immersed in
complex manifolds
\jour Trans. A.M.S.
\vol 172
\yr 1972
\pages 1--33
\endref\label{Lai}

\ref{LP}
\by F.~Laudenbach, V.~Poenaru
\paper A note on four-dimensional handlebodies
\jour Bull. Soc. math. France
\yr 1972
\vol 100
\pages 337--344
\endref\label{LP}

\ref{Let}
\by M.~Letizia
\paper Quotients by complex conjugation of nonsingular quadrics and
cubics in $\CP^3$ defined over $\R$
\jour Pacif. J. Math.
\yr 1984
\vol 110
\issue 2
\pages 307--314
\endref
\label{Leti}

\ref{Ma1}
\by A.~Marin
\paper $\CP^2/\sigma$ ou Kuiper et Massey au pays des coniques
\inbook A la recherche de la topologie perdue
\eds L.~Guillou, A.~Marin
\publ Birkh\"auser
\yr 1986
\endref\label{MarinonKuiperMassey}

\ref{Ma2}
\by A.~Marin
\paper Quelques remarques sur les courbes r\'eelles alg\'ebriques planes
\inbook A la recherche de la topologie perdue
\bookinfo Russian transl.
\lang Russian
\publ Mir
\yr 1989
\pages 162--172
\endref
\label{GM.zamecha}

\ref{Mas}
\by W.~Massey
\paper The quotient space of the complex projective
plane under conjugation is a $4$-sphere
\jour Geom. Dedicata
\yr 1973
\vol 2
\pages 371--374
\endref\label{Massey}

\ref{Mat1}
\by S.~Matsuoka
\paper Nonsingular algebraic curves in $\RP^1\times\RP^1$
\jour Trans. Amer. Math. Soc.
\vol 324
\issue 1
\yr 1991
\pages 87--107
\endref\label{Mat:91}

\ref{Mat2}
\by S.~Matsuoka
\paper Congruences for $M$- and $(M-1)$-curves
with odd branches on a hyperboloid
\jour Bull. London Math. Soc.
\yr 1992
\vol 24
\pages 61--67
\endref\label{Mat:92}

\ref{Mat3}
\by S.~Saito \rom(formerly S.~Matsuoka\rom)
\paper Classification of involutions of lattices with conditions and real
algebraic curves on a hyperboloid
\jour preprint
\pages 1--21
\endref\label{Matsuoka}

\ref{MS}
\by J.~W.~Milnor, J.~D.~Stasheff
\book Characteristic classes
\bookinfo Annals of Math. Studies, 76
\publ Princeton Univ. Press and Univ. Tokyo Press
\publaddr Princeton, New Jersey
\yr 1974
\endref\label{Milnor}

\ref{Mik1}
\by G.~Mikhalkin
\paper The complex separation of real surfaces and extensions of
Rokhlin congruence
\jour   Invent. Math.
\vol    118
\yr 1994
\pages  197--222
\endref\label{Mikhalkin}

\ref{Mik2}
\by G.~Mikhalkin
\paper Adjunction inequality for real algebraic curves
\jour   Math. Res. Lett.
\vol   4
\yr 1997
\issue 1
\pages  45--52
\endref\label{Mikh2}

\ref{Mik3}
\by G.~Mikhalkin
\paper Congruences for real algebraic curves on an ellipsoid
\inbook  Adv. Sovient Math.
\vol    18
\yr 1994
\pages  223--233
\endref\label{Mikh3}

\ref{Mik4}
\by G.~Mikhalkin
\paper Topology of curves of degree $6$ on cubic surfaces
in $\RP^3$
\jour J. Algebraic Geom.
\vol 7
\yr 1998
\issue 2
\pages 219--237
\endref\label{Mikh:cubic}

\ref{MH}
\by J.~Milnor, D.~Husemoller
\book Symmetric bilinear forms
\publ Springer-Verlag
\publaddr Berlin-Heidelberg-New York
\yr 1973
\endref
\label{MilnorHusemoller}

\ref{Na}
\by S.~M.~Natanzon
\paper Klein surfaces
\jour Uspekhi Mat. Nauk
\vol 45
\yr 1990
\issue 6
\pages 47--90
\endref\label{Natanzon:survey}

\ref{Nik1}
\by V.~V.~Nikulin
\paper Integral symmetric bilinear forms and some of their applications
\jour Izv\. Akad\. Nauk SSSR
\vol    43
\yr     1979
\issue  1
\pages  117--177 \lang Russian \transl\nofrills English transl. in
\jour Math\. USSR-Izv.
\vol    14
\yr     1980
\pages  103--167
\endref\label{Nikulin:forms}

\ref{Nik2}
\by V.~V.~Nikulin
\paper Involutions of integral quadratic forms
and their application to real algebraic geometry
\jour Izv. Acad. Nauk SSSR
\yr 1983
\vol 47
\pages 109--188 \lang Russian
\transl\nofrills English transl. in
\jour Math. USSR Izv.
\vol 22
\yr 1984
\endref\label{Nikulin:condition}

\ref{Och}
\by Ochanine
\paper Signature modulo 16, invariants de Kervaire generalis\'es et
nombres charact\'eristiques dans la $K$-th\'eorie re\'ele
\jour Mem. Soc. Math. France \rom(N.S.\rom)
\yr 1981
\issue  5
\endref\label{Ochanine}

\ref{Ol}
\by O.~A.~Oleinik
\paper On the topology of real algebraic curves on an algebraic
surface
\jour Matem. Sb.
\vol    29
\yr 1951
\pages  133--156 \lang   Russian
\endref\label{Oleinik}

\ref{Or1}
\by S.~Orevkov
\paper Link theory and oval arrangements of real algebraic curves
\jour Topology
\yr 1999
\vol 38
\issue 4
\pages 779--810
\endref\label{Or}

\ref{Or2}
\by S.~Orevkov
\paper Riemann existence theorem and construction of real algebraic curves
\toappear
\endref
\label{Or:RET}

\ref{Pet}
\by I.~G.~Petrovskiy
\paper On the topology of real plane algebraic curves
\jour Ann. Math
\vol    39
\issue  1
\yr 1938
\pages  189--209
\endref\label{Petrovsky}

\ref{PO}
\by I.~G.~Petrovskiy, O.~A.~Oleinik
\paper On the topology of real algebraic surfaces
\jour Isv. Akad. Nauk SSSR, Ser. mat.
\vol    13
\yr 1949
\pages  389--402 \lang   Russian
\endref\label{PetrovskyOleinik}

\ref{Ri}
\by J.-J.~Risler
\paper Construction d'hypersurfaces r\'eelles \rom(d'apre\`es Viro\rom)
\jour Ast\'ersique
\vol 216
\yr 1993
\pages 69--86
\endref\label{Risler}

\ref{R1}
\by V.~A.~Rokhlin
\paper Proof of a conjecture of Gudkov
\jour Funkc. Anal. i Prilozhen.
\yr 1972
\vol 6
\issue 2
\pages 62--64
\endref\label{Ro}

\ref{R2}
\by V.~A.~Rokhlin
\paper Complex orientation of real algebraic curves
\jour Funkc. Anal. i Prilozhen.
\yr 1974
\vol 8
\issue 4
\pages 71--75
\endref\label{Ro1}

\ref{R3}
\by V.~A.~Rokhlin
\paper Complex topological characteristics of real algebraic curves
\jour Uspekhi Mat. Nauk
\yr 1978
\vol 33
\issue 5
\pages 77--89
\endref\label{Ro2}

\ref{R4}
\by V.~A.~Rokhlin
\paper Congruences modulo $16$ in Hilbert's sixteenth problem.
\jour Funkc. Anal. i Prilozhen.
\yr 1972
\vol 6
\issue 4
\pages 58--64
\endref\label{Ro3}

\ref{R5}
\by V.~A.~Rokhlin
\paper Congruences modulo $16$ in Hilbert's sixteenth problem. II
\jour Funkc. Anal. i Prilozhen.
\yr 1973
\vol 7
\issue 2
\pages 91--92
\endref\label{Ro4}

\ref{R6}
\by V.~A.~Rokhlin
\paper New inequalities in the topology of real plane algebraic curves
\jour Funkc. Anal. i Prilozhen.
\yr 1980
\vol 14
\issue 1
\pages 37--43
\endref\label{Ro5}

\ref{S}
\by W.~Scharlau
\book Quadratic and Hermitian forms
\bookinfo Grundlehren der mathematische Wissenschaften, 270
\publ Springer-Verlag
\publaddr Berlin-New York
\yr 1985
\endref
\label{Scharlau}

\ref{Sch}
\by L.~Schl\"afli
\paper On the distribution of surfaces of the third order into species,
in reference to the absence or presence of singular points,
and the reality of their lines
\jour Phil. Trans. Roy. Soc. London
\yr 1863
\vol 153
\pages 195--241
\endref\label{Schlafli}

\ref{Seg}
\by B.~Segre
\book The Non-singular Cubic Surfaces
\publ Oxford University Press
\publaddr Oxford
\yr 1942
\endref
\label{Segre}

\ref{Sil}
\by R.~Silhol
\book Real Algebraic Surfaces
\bookinfo Lecture Notes in Math.
\vol    1392
\yr 1989
\publ   Springer-Verlag
\endref\label{Silhol}

\ref{St}
\by J.~Steenbrink
\paper Intersection form for quasi-homogeneous singularities
\jour Comp. Math.
\vol 34
\issue 2
\yr 1977
\pages 211--223
\endref
\label{Steen}

\ref{SE}
\by N.~E.~Steenrod, D.~B.~A.~Epstein
\book Cohomology operations
\publ Princeton University Press
\publaddr Princeton
\yr     1962
\endref\label{SteenrodEpstein}

\ref{Sz}
\by Z.~Szab\'o
\paper Simply-connected irreducible $4$-manifolds
with no symplectic structure
\jour Invent. Math.
\vol 132
\yr 1998
\pages 457--466
\endref\label{Szabo}

\ref{Th1}
\by R.~Thom
\paper Une th\'eorie intrins\`eque des puissance de Steenrod
\jour Colloque de Topologie de Strasbourg
\yr 1951
\issue\nofrills expos\'e N.~6.
\endref\label{ThomStras}

\ref{Th2}
\by R.~Thom
\paper Espaces fibr\'es en sph\`eres et carr\'es de Steenrod
\jour Ann. \'Ec. Norm.
\vol 3
\issue 2
\yr 1952
\pages 109--182
\endref\label{ThomAEN}

\ref{Th3}
\by R.~Thom
\paper Sur l'homologie des vari\'et\'es
alg\'ebriques r\'eelles
\inbook A symposium in honour of Marston Morse
\publ Princeton Univ. Press
\yr 1965
\pages 255--265
\endref\label{ThomSmith}

\ref{Va1}
\by A.~N.~Varchenko
\paper Theorems on the topological equisingularity of families
of algebraic varieties and families of polynomial mappings
\jour Izv. Acad. Nauk SSSR
\vol 36
\yr 1972
\pages 957--1019
\endref\label{Varchenko}

\ref{Va2}
\by A.~N.~Varchenko
\paper Asymptotic Hodge structure on vanishing
cohomology
\jour Izv. Akad. nauk SSSR Ser. Mat.
\vol 45
\yr 1981
\issue 3
\pages 540--591
\endref\label{VarHodge}

\ref{Va3}
\by A.~N.~Varchenko\mnote{Kh:to add this ref
into r-tools to PetrOlei? then to add also Loezer?}
\paper Local residue and the intersection form
in vanishing cohomology
\jour Izv. Akad. Nauk SSSR, Ser. Matem.
\vol 49
\yr 1985
\issue 1
\pages 32--54
\endref\label{VarPO}


\ref{V1}
\by O.~Ya.~Viro
\paper Real algebraic plane curves: constructions
with controlled topology
\jour Leningrad Math. J.
\vol 1
\yr 1990
\pages 1059--1134
\endref\label{Viro:survey}

\ref{V2}
\by O.~Ya.~Viro
\paper Construction of $M$-surfaces
\jour Funct. Anal. Appl.
\vol 13
\issue 3
\yr 1979
\endref\label{ViroM}

\ref{V3}
\by O.~Ya.~Viro
\paper Achievements in the topology of real algebraic
varieties in the last six years
\jour Uspekhi Mat. Nauk
\vol 41
\yr 1986
\issue 3
\pages 45--67
\endref\label{Viro.congres}

\ref{V4}
\by O.~Ya.~Viro
\paper Complex orientations of real algebraic surfaces
\jour Adv. Soviet Math.
\yr 1994
\pages 261--284
\vol 18
\endref
\label{Viro.???}

\ref{Wl}
\by C.~T.~C.~Wall
\paper Quadratic form in finite groups and related topics
\jour Topology
\vol 2
\yr 1964
\pages 281--298
\endref\label{Wall}

\ref{Wc}
\by A.~Wallace
\paper Linear sections of algebraic varieties
\jour Indiana Univ. Math. J.
\vol 20
\yr 1971
\pages 1153--1162
\endref\label{Wallace}

\ref{Wa}
\by Sh.~Wang
\paper A vanishing theorem for
Seiberg-Witten invariants
\yr 1995
\jour Math. Res. Lett.
\vol 2
\issue 3
\pages 305--310
\endref\label{Wang}

\ref{We1}
\by J.-Y.~Welschinger
\paper Orientations complexes des $J$-courbes r\'eelles
\jour Pr\'epublications de \newline l'Univ.~Louis~Pasteur, Strasbourg
\yr 1999
\endref\label{W1}

\ref{We2}
\by J.-Y.~Welschinger
\paper $J$-courbes r\'eelles \`a nids profonds
sur les surfaces r\'egl\'ees
\jour Pr\'epu\-blications de l'Univ.~Louis~Pasteur, Strasbourg
\yr 1999
\endref\label{W2}

\ref{Wi}
\by G.~Wilson
\paper Hilbert's sixteenth problem
\jour Topology
\vol 17
\pages 53--73
\endref\label{Wilson}

\ref{Ze}
\by H.~G.~Zeuthen
\paper Etudes des propri\'et\'es de situation
des surfaces cubiques
\yr 1875
\jour Math. Ann
\vol 8
\pages 1--30
\endref\label{Zeuthen}

\ref{Zv1}
\by V.~I.~Zvonilov
\paper Harlamov's inequalities and Petrovsky-Oleinik inequalities
\jour Funkc. Anal. i Prilozhen.
\vol 9
\yr 1975
\issue 2
\pages 69--70
\endref
\label{ZvonPetr}

\ref{Zv2}
\by V.~I.~Zvonilov
\paper Complex topological characteristics of real algebraic curves
on surfaces
\jour   Funkc. Anal. i Prilozhen.
\vol    16
\issue  3
\yr 1982
\pages  56--57 \lang   Russian
\moreref\nofrills English transl. in
\jour   Functional Anal. Appl.
\vol    16
\issue  3
\yr     1982
\endref\label{Zvonilov:RF}

\ref{Zv3}
\by V.~I.~Zvonilov
\paper Complex topological invariants of real algebraic curves
on a hyperboloid and on an ellipsoid
\jour   Algebra i Analiz
\vol    3
\issue  5
\yr 1991
\lang   Russian
\moreref\nofrills English transl. in
\jour   St. Petersburg Math. J.
\vol    3
\issue  5
\yr     1992
\pages  1023--1042
\endref\label{Zvonilov}

\ref{Zv4}
\by V.~I.~Zvonilov
\paper Stratified spaces of real algebraic curves
 of bidegree $(m,1)$ and $(m,2)$ on a hyperboloid
\jour Amer. Math. Soc. Transl. \rom{(2)}
\yr 1996
\vol 173
\pages 253--264
\endref\label{Zv4}


\endRefs
\enddocument